

\documentstyle{amsppt}
\magnification=\magstep1
\hsize=6 truein
\hcorrection{.375in}
\vsize=8.5 truein
\parindent=20pt
\baselineskip=14pt
\TagsOnRight
\NoBlackBoxes
\footline{\hss\tenrm\folio\hss}

\centerline{ Invariant Measures for}
\centerline{Unitary Forms of Kac-Moody Groups}

\vskip1truein

\centerline{Doug Pickrell}
\centerline{Mathematics Department}
\centerline{University of Arizona}
\centerline{Tucson, Arizona 85721 }

\vskip1truein

\flushpar Abstract.  The purpose of this paper is to describe
some conjectures and results on the existence and uniqueness
of invariant measures on formal completions of Kac-Moody
groups and associated homogeneous spaces.  The basic invariant
measure is a natural generalization of Haar measure for the
unitary form of a finite type Kac-Moody group, and its
projection to flag spaces is a generalization of the normalized
invariant volume element.  The other ``invariant measures''
are actually measures having values in line bundles over these
spaces; these bundle-valued measures heuristically arise from
coupling the basic invariant measure to hermitian structures
on associated line bundles, but in the examples that we
consider, they are in fact singular with respect to the basic
invariant measure.  The main examples considered are infinite
classical groups (i.e.  affine Kac-Moody groups of infinite
rank), central extensions of loop groups, and the Virasoro group
(which is not technically a Kac-Moody group).

\vfill
\eject

\centerline{\S 0. General Introduction}

\bigskip

In this introduction we will begin by describing some general
conjectures about the existence and uniqueness of invariant
measures for unitary forms of symmetrizable Kac-Moody
groups.  In the body of the paper we prove the existence parts
of the conjectures for affine algebras (of possibly infinite
rank).

The structure necessary to define and work with Kac-Moody
algebras and groups is extensive.  I will use the book [Kac1]
and the survey article [Kac2] as general references for the
basics of the theory.

\bigskip

\flushpar\S 0.1. Some Conjectures.

Let $A$ denote an indecomposable generalized $n\times n$ Cartan matrix
which is symmetrizable.  Let $G(A)$ denote the Kac-Moody group
associated to $A$, in the sense introduced by Kac and Peterson
([KP1]), and let $K(A)$ denote the unitary real form of $G(A)$.  If
$A$ is of finite type, i.e.  if $n$ is finite and $A$ is a classical
Cartan matrix, then in particular $K(A)$ is a compact topological
group, so that there exists a unique $K(A)$-invariant probability
measure supported on $K(A)$, the Haar measure of $K(A)$.
Conversely the existence of a $K(A)$-invariant measure class
supported on $K(A)$ implies that $K(A)$ is locally compact, by
Mackey and Weil's converse to Haar's theorem, hence that $A$ is
of finite type.

To go beyond the finite type case we are guided by the fact
that the most elemental functions which we would like to be
able to integrate, namely the matrix coefficients of integrable
highest weight representations of $G(A)$, naturally define regular
functions on a space which we will call the formal completion
of $G(A)$, and which we will denote either by $\Cal G$ or $G(A)_{
formal}$.
This space is a complex manifold which is modelled on the
space
$$\frak n^{-}_{formal}\times \frak h_{formal}\times \frak n^{+}_{
formal},$$
where $\frak n^{\pm}_{formal}$ is the direct product of the root spaces
corresponding to positive roots (resp.,  negative roots), and
$\frak h_{formal}$ is the direct product of the coroot spaces (of course
if $n<\infty$, then $\frak h_{formal}=\frak h$).  The transition functions
arise
from continuously extending the translational action of $G(A)$ on
itself, expressed in terms of the coordinates provided by the
top stratum of the Birkhoff decomposition.  The two-sided
translational action of $G(A)$ on itself extends continuously to a
holomorphic two-sided action of $G(A)$ on $\Cal G$.  If $A$ is of finite
type, then $\Cal G=G(A)$; otherwise, the space $\Cal G$ is not a group, and
there does not exist an involution which fixes a completion of
$K$(A).  The description of this formal completion is the subject
of \S 1 of Part I.

{}From the point of view of this paper, the rationale for
considering such a large space enveloping $K(A)$ is that the
process of constructing a measure on an infinite dimensional
space inevitably involves some regularization scheme, and
regularization tends to spread out the support of the measure.
In our context the support definitely tends to spread out in
nonunitary directions from $K(A)$ in $\Cal G$.  In the long run we
would like to precisely describe in geometric terms the
support of any measure that we succeed in constructing.
However this issue is distinct from and secondary to the issue
of existence.

Let $C$ denote the center of the group $G$ (if $n$ is finite, then
$G=G(A)$; if $n$ is infinite, then $G$ is slightly larger than $G
(A)$
(see subsection (1.4) of Part I)).  Given a holomorphic character
$k$ of the identity component, $C_0$, we can consider the space of
matrix coefficients which have level $k$, i.e.  the space of
regular functions on $\Cal G$ satisfying
$$\phi (g\cdot\lambda )=k(\lambda )\phi (g),\tag 0.1.1$$
for $g\in \Cal G$ and $\lambda\in C_0$.  This space can be identified with the
space of regular sections of a line bundle,
$$\Cal L^{*}_k\to \Cal G/C_0,\tag 0.1.2$$
associated to the character $k$.  For this space to be
nonvacuous, $k$ must satisfy a positivity condition; specifically,
$k$ must be the restriction of a character defined by a dominant
integral functional (the trivial character is allowed).  We
assume that this condition is satisfied.

If $n$ is finite, then it is known that $C_0$ is nontrivial precisely
when $detA=0$.  In this case $A$ is of affine type, meaning that
$G/C_0$ is a (possibly twisted) polynomial loop group, $C_0=\Bbb C^{
*}$,
and we can regard $k$ as a nonnegative integer.   In the case
that $detA\ne 0$, the regular sections of $\Cal L_k^{*}$ are just the regular
functions on $\Cal G$.

\proclaim{(0.1.3) Conjecture} There exists a measure $\mu =\mu^{\vert
\Cal L_k\vert^2}$
on $\Cal G/C_0$ having values in $\vert \Cal L_k\vert^2\equiv \Cal L_
k\otimes\bar {\Cal L}_k$ such that (i) $\mu$ is
$K(A)$-biinvariant, and (ii) for each regular section $\sigma$ of $
\Cal L_k^{*}$,
$(\sigma\otimes\bar{\sigma })d\mu$ is a finite positive measure.  If $
n<\infty$, and if $A$ is
not of finite type, then the measure is unique up to a constant
multiple.  \endproclaim

\flushpar(0.1.4)Remarks. (1) To say that $\mu$ is a measure having values
in $\vert \Cal L\vert^2$ means that for each local section $s$ of $
\Cal L^{*}$, $(s\otimes\bar {s})d\mu$ is a
measure on the base in the ordinary sense, and if $s'$ is another
such section, then $(s'\otimes\bar {s}')d\mu$, restricted to the set $
\{s\ne 0\}$, is
absolutely continuous with respect to $(s\otimes\bar {s})d\mu$, and the
Radon-Nikodym derivative is given by $\vert\frac {s'}s\vert^2$.

(2) Using $\mu$, we can define a Hilbert space structure on a
subspace of the space of sections of $\Cal L^{*}$ by
$$\langle s_1,s_2\rangle =\int (s_1\otimes\bar {s}_2)d\mu .$$
Since $\mu$ is $K(A)$-biinvariant, there is a natural unitary action of
$K(A)\times K(A)$ on this Hilbert space of sections.  The subspace of
square-integrable holomorphic sections will be invariant, and
this will provide an analytic setting for Kac and Peterson's
generalization of the Peter-Weyl theorem.  Conjecturally these
actions will extend continuously to the completion of $K(A)$
constructed in [KP2] (the strong operator completion of $K(A)$ in
a faithful highest weight module), which we will refer to as
the Kac-Peterson completion of $K(A)$.

(3) If (i) and (ii) do not suffice to characterize the measure,
then we expect that the diagonal distribution (which is a rough
analogue of the spherical transform) will faithfully map the
space of measures satisfying (i) and (ii) to a space of measures
on the Cartan subalgebra.  The preferred measure $\mu$ will
correspond to a regularization of Harish-Chandra's {\bf c}-function
(see (0.2.8) below).

\smallskip

In this paper we will prove the existence part of (0.1.3) in all
affine cases, i.e.  when $G(A)$ is either an infinite classical
group, or a central extension of a (possibly twisted) loop group.
In the first case we will be able to produce relatively explicit
expressions for the measures.  In the second case we will see
that there is a relation between $\mu$ and Brownian motion in the
limit as temperature tends to infinity (the observation that
Wiener measure in the presence of high temperature is
approximately invariant is due to the Malliavins; see [MM2]).
In this second case there are some conjectural formulas for
the measures.

The analogue of (0.1.3) for homogeneous spaces is the following

\proclaim{(0.1.5) Conjecture} Let $\Cal L^{*}\to \Cal F$ denote a positive line
bundle, where $\Cal F$ is the formal completion of the flag space
$G(A)/B^{+}$.  Then up to a multiple there exists a unique measure
$\mu$ on $\Cal F$ having values in $\vert \Cal L\vert^2$ such that (i) $
\mu$ is $K(A)$-invariant,
and (ii) $\vert\sigma\vert^2d\mu$ is a finite positive measure, for each
regular
section $\sigma$ of $\Cal L^{*}$.\endproclaim

This second conjecture sounds more plausible than the first
because the formal flag space is closer to being a homogeneous
space for $K(A)$ than is $\Cal G$.  However, what should be true is
that the nonunitary directions in $\Cal G$ can be recovered from
asymptotic information on the flag space.

\proclaim{(0.1.6)Conjecture}If $n<\infty$, and if $A$ is not of finite
type, then the natural maps
$$\matrix \Cal L_k/T&\rightarrow&\Cal L_k/\Cal B^{+}\\\downarrow&&
\downarrow\\\Cal G/TC_0&\rightarrow&\Cal F&\endmatrix $$
induce a measure algebra isomorphism of the respective
$K(A)$-invariant measures.\endproclaim

\bigskip

\flushpar\S 0.2. Basic intuition.

The flag space
$$\Cal F=K(A)/T=G(A)/B^{+}\tag 0.2.1$$
is not in general smooth.  In the informal discussion that
follows, for ease of exposition, we will overlook this point.

Because the Cartan matrix $A$ is symmetrizable, there is a
essentially unique nondegenerate symmetric bilinear form on
$\frak g(A)$.  In a standard way, using the Cartan involution defining
$\frak k(A)$, this induces a $K(A)$-invariant Kahler form $\omega$ on $
\Cal F$ (the
fact that this form is positive was proven by Garland in the
affine case, and Kac and Peterson in general; see [KP2]).

Now suppose that $\Cal L^{*}\to \Cal F$ is a positive line bundle corresponding
to a dominant integrable weight $\Lambda$.  There is a essentially
unique $K(A)$-invariant hermitian structure on $\Cal L^{*}$.  Let $
\sigma =\sigma_{\Lambda}$
denote the essentially unique regular section of $\Cal L^{*}\to \Cal F$ which
is
fixed by $N^{-}$.  Given a cell $C_w\subset \Cal F$, the positive measure
$$\vert\sigma\vert^{2s}\omega^{l(w)}\vert_{\bar {C}_w},\tag 0.2.2$$
on the Schubert variety $\bar {C}_w$, is finite, hence it can be
renormalized; let
$$d\mu_s^{(w)}=\frac 1E\vert\sigma\vert^{2s}\omega^{l(w)}\vert_{\bar {
C}_w}\tag 0.2.3$$
denote the corresponding probability measure (here $s$ is any
real number $\ge 0$).

If $w$ is expressed as a product of fundamental reflections,
$$w=r_{i_l}\cdot ...\cdot r_{i_1},\tag 0.2.4$$
then the corresponding Schubert variety $\bar {C}_w$ is invariant under
the action of the copy of $SU(2,\Bbb C),$ $K_{(i_l)}$, corresponding to the
simple root $\alpha_{i_l}$.  Since the symplectic form $\omega$ is
$K(A)$-invariant, it follows that
$$\vert\sigma\vert^{2s}d\mu =\lim_{w\to\infty}d\mu_s^{(w)}\tag 0.2.5$$
defines a $K(A)$-quasi-invariant measure with the appropriate
Radon-Nikodym derivative, provided that the limit exists.
The basic question then is ``why does the limit exist?''.

For simplicity suppose that $n$ is finite.  Because of the formal
structure of the space $\Cal F$ ( specifically, because Kolmogorov's
theorem can be applied in affine coordinate patches), it turns
out that it essentially suffices to analyze the measures
$$proj_{*}(d\mu_s^{(w)})\in Prob(X_M),\tag 0.2.6$$
for each $M$, where
$$X_M=(N_M^{-}\backslash G(A)/B^{+})_{ss}/equivalence\tag 0.2.7$$
is the projective variety obtained by applying geometric
invariant theory to the semistable points of the double coset
space, and where $N^{-}_M$ is the subgroup of $N^{-}$ corresponding to
the Lie subalgebra of root spaces of height $\le -M$.

The important point is that the space $X_M$ is compact, hence
the measures (0.2.6) have limit points, as $w$ tends to infinity in
an appropriate sense.   Unfortunately, for this method to
work, we need to know that these Schubert approximations
have a limit point which is simultaneously invariant under all
of the $K_{(j)}$.  For the infinite classical groups there is a unique
limit, but I do not know whether there is a unique limit even
in the case of loop groups, where another method is known to
produce an invariant measure.

Plausible or not, this same kind of heuristic argument applies
to the construction of measures on the space $\Cal G/C_0$, provided
that $C_0$ is nontrivial.  One replaces the Schubert variety $\bar {
C}_w$
by its inverse image in $K(A)$, which is a torus bundle over $\bar {
C}_w$.
The condition on $C_0$ is important because in this case the
diagonal distribution, i.e.  the pivots for the Bruhat
decomposition of $\Cal G/C_0$, are expressed in terms of $\underline {
ratios}$ of
fundamental matrix coefficients; these ratios do not require
regularization.

This approximation scheme definitely does not apply to the
affine space $\Cal G$, for the diagonal distribution is expressed in
terms of fundamental matrix coefficients, and these do not
have limits.  If the Cartan matrix is noninvertible, then we
expect that there will not exist an invariant measure on $\Cal G$
itself (we can prove this for the infinite classical groups).  If
the Cartan matrix $A$ is invertible, then a priori it is
conceivable that there may be some way to scale the torus
bundle in the preceding paragraph to obtain a nontrivial limit.
Since so little is known about nonaffine Kac-Moody groups, it
is difficult to test this directly, but it may be possible to
explore the idea in the following indirect way.

Suppose first that $A$ is of finite type.  We can write a generic
$g\in G(A)$ as
$$g=l\cdot diag\cdot u\tag 0.2.8$$
where $l\in N^{-}$, $diag\in H$, and $u\in N^{+}$.  In turn we can write
$diag=ma$, corresponding to $H=T\cdot exp(\frak h_{\Bbb R})$.  Harish-Chandra
proved that
$$\int_{K(A)}e^{-i\lambda loga(k)}dm(k)=\bold c(\rho -i\lambda )=
\prod_{\alpha >0}(1-i\frac {\langle\alpha ,\lambda\rangle}{\langle
\alpha ,\rho\rangle})^{-1},\tag 0.2.9$$
where $\rho$ is the sum of the positive roots, and $\langle\cdot
,\cdot\rangle$ is a
nondegenerate invariant form.  This is the Fourier transform
of the diagonal distribution for the Haar measure of $K(A)$.
(More generally we have
$$\frac 1E\int_{K(A)}\vert\sigma_{\Lambda}(k)\vert^{2s}e^{-i\lambda
loga(k)}dm(k)=\bold c(\rho +2s\Lambda -i\lambda )\tag 0.2.10$$
where $E$ is a normalization constant).

For the infinite classical groups we will see that the diagonal
distribution for our analogue of Haar measure is given by
essentially this same formula.  We believe the same is roughly
true in the loop group case, although as of this writing, this is
merely a conjecture.  It clearly is of interest, independent of
considerations of measure, to determine whether there is a
canonical analytical interpretation of a regularized {\bf c}-function
for a general Kac-Moody group.

\smallskip

\flushpar(0.2.11)Remark. The formula (0.2.9) is closely related to
the Weyl dimension formula. For if $\Lambda /2$ is dominant integral
and $\sigma_{\Lambda}$ is the corresponding fundamental matrix coefficient
normalized by $\sigma_{\Lambda}(1)=1$, then
$$\aligned
\int_{K(A)}a^{\Lambda}dm(g)&=\int\vert\sigma_{\Lambda /2}(g)\vert^
2dm(g)\\=dimn(\pi_{\Lambda /2})^{-1}&=(\prod_{\alpha >0}\frac {\langle
\rho +\Lambda ,\alpha\rangle}{\langle\rho ,\alpha\rangle})^{-1}\endaligned
$$
I believe that it has been observed previously that the
standard modules have finite formal dimensions in the affine
case.  One naturally suspects that these formal dimensions will
occur in the discrete part of a Plancherel formula.

\bigskip

\flushpar\S 0.3. Contents of these notes.

The contents of these notes is fairly well described by the
Table of Contents below.  We will just add a few comments
about part I, which concerns what little we can presently say
about Kac-Moody groups in general.  In section 1 there is a
description of the formal completions of $G(A)$ and the
associated flag space.  In section 2 we formulate the
Kac-Peterson generalization of the Peter-Weyl theorem in the
language of section 1.  We also indicate a criterion for a given
Kac-Moody algebra to have an associated complex Lie group.
This is all relatively straightforward.  The measure-theoretic
results are in section 3.  The first is that the stable points
have full measure for any $K(A)$-quasi-invariant measure.  The
second is a compactness result for certain approximation
schemes which we will use.  This, together with the
asymptotic invariance of Wiener measure (which we prove in
Part III), immediately implies the existence of an invariant
measure for flag spaces in affine cases.  The third is a
heuristic proof of uniqueness of measures on $\Cal F$.

\bigskip

\centerline{\S 0.4. Table of Contents.}

\bigskip

\flushpar General Introduction.

\smallskip

\flushpar Part I. General Results.

\S 1. The Formal Completions of $G(A)$ and $G(A)/B^{+}$.

\S 2. Matrix Coefficients.

\S 3. Measures on the Flag Space $\Cal F$.

\bigskip

\flushpar Part II. Infinite Classical Groups.

\S 0. Introduction for Part II.

\S 1. Measures on Flag Spaces.

\S 2. The Case $\frak g=sl(\infty ,\Bbb C)$.

\S 3. The Case $\frak g=sl(2\infty ,\Bbb C).$

\S 4. The Other Cases.

\bigskip

\flushpar Part III. Loop Groups.

\S 0. Introduction for Part III.

\S 1. Extensions of Loop Groups.

\S 2. Hyperfunction, Formal and Measurable Loop Spaces.

\S 3. Existence of the Measures $\nu_{\beta ,k}$, $\beta >0$.

\S 4. Existence of Invariant Measures.

\bigskip

\centerline{Part I. General Theory.}

\bigskip

\centerline{\S 1. The Formal Completions of $G(A)$ and $G(A)/B$}

\bigskip

\flushpar\S 1.1 The Kac-Moody group $G(A)$.

Given an n x n generalized Cartan matrix $A$, the associated
Kac-Moody algebra $\frak g=\frak g'(A)$ is a certain complex Lie algebra
generated by 3n letters $e_i,$ $f_i,$ and $h_i$ (1$\le i<n+1)$.  The
possibility $n=\infty$ is allowed.  These letters satisfy the
Chevalley-Serre relations \newline
\centerline{$[h_i,h_j]=0,$ $[e_i,f_i$ $]=h_i$, $[e_i$ $,f_j$ $]=0
,$     $i\ne j$}
\centerline{$[h_i$ $,e_j]=a_{ij_{}}e_j$, $[h_i,f_j]=-a_{ij}f_j,$}
\centerline{(ad $e_i)^{1-a_{ij}}e_j=0,$ (ad $f_i)^{1-a_{ij}}f_j=0
,$     $i\ne j.$   }
We will always assume that $A$ is symmetrizable.  In this case
it is known that these relations actually define the Lie algebra
$\frak g$.

The Lie algebra $\frak g$ has a triangular decomposition\newline
\centerline{$\frak g=\frak n^{-}\oplus \frak h\oplus \frak n^{+}$}
where $\frak n^{-}$ is generated by the $\{f_j\},$ $ $ $\frak h$ by the $
\{h_j\},$ and $\frak n^{+}$ by the $\{e_j\}.$
As a Lie algebra $\frak g$ is graded by its root space decomposition; if
$n=\infty$, then we will view $\frak g_0=\frak h$ as the direct sum $
\oplus \Bbb Ch_j$, so that the
graded pieces will all be finite dimensional.  It also graded as a Lie
algebra by height, where the heights of $f_j,$ $h_j,$ and $e_j$ are -1, 0, and
$+1$, respectively. We will write
$$heights(x)\ge k\quad (resp,\quad\le k)\tag 1.1.1$$
if each of the nonzero homogeneous components of $x$ has height $
\ge k$
(resp, $\le k$)

For each $1\le j<n+1$, the Lie subalgebra $\frak g_{(j)}=$ $\Bbb C
f_j$ $\oplus$ $\Bbb Ch_j$ $\oplus$ $\Bbb Ce_j$ is
isomorphic to $sl(2,\Bbb C).$ In each integrable module of $\frak g${\it ,} the
action by
$\frak g_{(j)}$ can be exponentiated to an action by a copy $G_{(
j)}$ of $SL(2,\Bbb C).$ The
group $G(A)$, the minimal Kac-Moody group associated to $A$, is the
group generated by these n copies of $SL(2,\Bbb C),$ subject to the
relations imposed by considering all integrable modules of $\frak g$.  Each
$G_{(j)}$ is faithfully represented in $G(A)$.  Below we will use the fact that
$\oplus L(\Lambda_i)$, the sum of the fundamental highest weight modules, is a
faithful $G(A)-$module.

The group $G(A)$ is not a Lie group (except in the finite type
case).  In the introduction we mentioned the measure-theoretic
motivation for introducing the formal completion for $G(A)$.  In
algebraic geometry formal completions form a bridge between
analysis and algebra; see Remark (1.5.2) below and [Hart].  In
\S 2.4 below we will comment on completions of $G(A)$ from a
different perspective.

\bigskip

\flushpar\S 1.2 Formal completions, proalgebraic complex manifolds.

Given a countable direct sum of finite dimensional complex vector
spaces, $V=\oplus V_k,$ the formal completion of $V$, denoted $V_{
formal},$ is the
topological vector space $\prod V_k$, where the topology is the product
topology.

\proclaim{(1.2.1) Lemma} Suppose that $U$ is a connected open
subset of $V_{formal}$.  If $f:U\to \Bbb C$ is a holomorphic function, then $
f$
depends on finitely many variables; that is, there is a domain
$U_1\subset V_{k_1}\times ...\times V_{k_n}$ and a holomorphic function $
F$ on $U_1$ such that
$U$ is contained in the inverse image of $U_1$ in $\prod V_k$, and
$f((v_k))=F(v_{k_1},..,v_{k_n})$.  \endproclaim

\flushpar Proof.  We can assume that the product is indexed by the
natural numbers.  The differential of $f$ can be regarded as a map
$D:U\to\oplus V_k^{*}$.  Let $D_k$ denote the kth component of $D$.  Then
$$U=\bigcup_N\bigcap_{k>N}\{D_k=0\}.$$
Since $U$ is of second category, this implies that for some $N$
$$\bigcap_{k>N}\{D_k=0\}$$
contains an open subset of $U$.  Because $D$ is holomorphic, it follows
that $D_k$ is identically zero on $U$ for all but finitely many $
k$.  This
implies that $f$ depends upon finitely many variables.//

\smallskip

In connection with Remark (1.5.2) below, we note that this has
the following

\proclaim{(1.2.2)Corollary}Suppose that S is a connected finite
dimensional complex manifold. Then
$$H^q(V_{formal}\times S;\Cal O)=lim\,\,H^q(V_1\times ...\times V_
N\times S;\Cal O).$$
In particular if S is a Stein manifold, then
$$H^q(V_{formal}\times S;\Cal O)=0,\,\,\,q>0.$$
\endproclaim

We will use the following terminology.  A function $f$ as in the
lemma is a polynomial (resp.,  a rational function) if $F$ is a
polynomial (resp.,  a rational function).  A function $U\to\prod
V_k$ is
polynomial (resp.,  rational) if each component is polynomial (resp.,
rational).  A complex manifold modelled on $V_{formal}$, together with a
distinguished system of charts, is proalgebraic if each distinguished
chart is modelled on a space of the form $\prod V_k\setminus D$, where $
D$ is the
zero set of finitely many polynomials, and if the transition functions
are rational.  A function on such a manifold is regular if it is
represented by a polynomial in each distinguished chart.

A second countable, Hausdorff manifold modelled on $V_{formal}$
admits partitions of unity.  This is because for each open
subset $U$ of $V_{formal}$, there exists a smooth function (of
finitely many variables) which vanishes on the complement of
$U$.  Hence the DeRham Theorem holds for such a manifold.

\bigskip

\flushpar\S 1.3. The groups $\Cal N^{+}$ and $\Cal G^{+}$.

Because $\frak n^{+}$ with the root space gradation is a graded profinite
nilpotent Lie algebra, the formal completion $\frak n^{+}_{formal}$ is a
proalgebraic profinite nilpotent Lie algebra.  By definition
$\Cal N^{+}_{formal}$, which we will usually abbreviate to $\Cal N^{
+}$, is the
corresponding simply connected proalgebraic Lie group.  It can
be realized explicitly in the following way:  the exponential
map is an isomorphism of proalgebraic manifolds,
$$exp:\frak n^{+}_{formal}\to \frak N^{+},\tag 1.3.1$$
and in these coordinates the multiplication is given by the
Baker-Campbell-Hausdorff formula, which is clearly a polynomial
map.

Now suppose that $L$ is a direct sum of integrable highest
weight modules, e.g.  $L=\oplus L(\Lambda_i).$ Let $M$ denote a finite
collection of positive roots, and let $\frak n^{+}_M$ denote the largest
graded subalgebra of $\frak n^{+}$ which intersects $\frak g_{\alpha}$
trivially, for each
$\alpha\in M$.  Then $\frak n^{+}/\frak n^{+}_M$ is finite dimensional, and the
space $
L^{\frak n^{+}_M}$
is a direct sum of finite dimensional spaces.  Since
$L=\sum_ML^{\frak n^{+}_M}$ and $\frak n^{+}_{}/\frak n^{+}_M=\frak n^{
+}_{formal}/\frak n^{+}_{M,formal}$, it follows that
the action of $\frak n^{+}$on $L$ extends to an action
$$\frak n^{+}_{formal}\times L\to L.\tag 1.3.2$$
In turn this action can be exponentiated to a group action
$$\frak N^{+}\times L\to L.\tag 1.3.3$$

Note that $L$ is a faithful $\frak n^{+}$-module if and only if $
\frak n^{+}/\frak n^{+}_M$ acts
faithfully on $L^{\frak n^{+}_M}$ for each $M$.  The same is true with $
\Cal N^{+}$ in
place of $\frak n^{+}$.  Since $\frak n^{+}/\frak n^{+}_M$ is the Lie algebra
of $
\Cal N^{+}/\Cal N^{+}_M$, and
since this group acts unipotently on $L^{\frak n^{+}_M},$ it follows that $
\frak n^{+}$
acts faithfully on L if and only if $\Cal N^{+}$ acts faithfully on $
L$.

For a general $L$ we always have $\Cal N^{+}\cdot G(A)=G(A)\cdot
\Cal N^{+}$.  To see
this it suffices to check that
$$\frak N^{+}\cdot H\cdot G_{(j)}=G_{(j)}\cdot \frak N^{+}\cdot H
,\tag 1.3.4$$
where {\it H} corresponds to {\it h.} This follows from the fact that
$\frak n^{+}\oplus \frak h\oplus \Bbb Ce_j$ is a parabolic subalgebra of $
\frak g${\it .}

Now suppose that $n$, the size of the Cartan matrix, is finite.  The
preceding paragraph justifies the following definition:
$$\Cal G^{+}=\frak N^{+}\cdot G(A)=G(A)\cdot \frak N^{+}.\tag 1.3.5$$
Similarly we define $\Cal G^{-}$. Note that $\Cal G^{+}$ acts on each
integrable highest
weight module of {\it g,} whereas $\Cal G^{-}$ acts in each integrable lowest
module of $\frak g$, in each case extending the action of $G(A)$.

If $n$ is infinite, then we must modify this definition slightly. The
change which we must incorporate is to replace $\frak h$ by $\frak h_{
formal}$. Let $\Cal H$
denote the corresponding infinite complex torus; precisely,
$$\Cal H=\prod_{1\le j<\infty}H_{(j)},\tag 1.3.6$$
where $H_{(j)}$ is the torus in $G_{(j)}$.  In the preceding paragraph we then
replace $G(A)$ by the group
$$G=G(A)\cdot \Cal H=\Cal H\cdot G(A).\tag 1.3.7$$
Of course if $n$ is finite, then $G=G(A)$.

Let $B^{+}$ denote the upper Borel subgroup of $G$, and let $\Cal B^{
+}=\Cal H\cdot \Cal N^{+}$
denote the upper Borel subgroup for $\Cal G^{+}$.  We will denote the common
lower Borel subgroup of {\it G} and $\Cal G^{+}$ by $B^{-}$.  Because $
\Cal G^{+}/\Cal B^{+}=G/B^{+},$ the
Birkhoff decomposition and unique factorization for the top stratum
of $G$, \newline
\centerline{$G=\bigsqcup B^{-}\cdot w\cdot B^{+}$,}
\centerline{$B^{-}\cdot 1\cdot B^{+}=N^{-}\times \Cal H\times N^{
+}$,}
induce a Birkhoff decomposition and unique factorization for the top
stratum of $\Cal G^{+}$:\newline
\centerline{$\Cal G^{+}=\bigsqcup$ $ $ $B^{-}\cdot w\cdot \Cal B^{
+}$,  }
\centerline{$B^{-}\cdot 1\cdot \Cal B^{+}=N^{-}\times \Cal H\times
\Cal N^{+}$,}
where the disjoint union is over all elements of the Weyl group
$W=N(H)/H$.  Similarly there is a Birkhoff decomposition and unique
factorization for the top stratum of $\Cal G^{_{-}}$.

\bigskip

\flushpar\S 1.4. The space $\Cal G=G(A)_{formal}$.

We can now define the formal completion of the Kac-Moody group
$G(A)$ as the quotient \newline
\centerline{$\Cal G=G(A)_{formal}=\Cal G^{-}\times_G\Cal G^{+}$,}
where $G$ acts on $\Cal G^{-}$ by multiplication on the right, and on $
\Cal G^{+}$ by
multiplication on the left. For $(g,h)\in \Cal G^{-}\times \Cal G^{
+}$, $[g,h]$ will denote its
image in $\Cal G$. There is a natural inclusion\newline
\centerline{$G\to \Cal G:g\to [g,1]=[1,g]$.}
Note that $\Cal G^{-}$ acts on $\Cal G$ from the left, and $\Cal G^{
+}$ acts on $\Cal G$ from the right.

\proclaim{(1.4.1) Proposition}  (a) $\Cal G=\bigsqcup \Cal B^{-}\cdot
w\cdot \Cal B^{+}$, where the disjoint
union is over all elements of the Weyl group $W=N(H)/H$.

(b) Each element of the top stratum $\Cal B^{-}\cdot 1\cdot \Cal B^{
+}$ can be written
uniquely in the form $l\cdot diag\cdot u$, where $l\in \Cal N^{-}$, $
diag\in \Cal H$, $u\in \Cal N^{+}$.
More generally
$$\Cal N^{-}_w\times \{w\}\times \Cal B^{+}\cong \Cal B^{-}\cdot
w\cdot \Cal B^{+}\cong \Cal B^{-}\times \{w\}\times \Cal N^{+}_w,$$
where $\Cal N^{-}_w$ $=\Cal N^{-}\cap w\Cal N^{-}w^{-1}$, and $\Cal N^{
+}_w=w^{-1}\Cal N^{+}w\cap \Cal N^{+}$.

(c) $\Cal G$ has the structure of a $\Cal G^{-}\times \Cal G^{+}$-homogeneous
proalgebraic
complex manifold, modelled on the space $\frak n^{-}_{formal}\times
\Cal H\times \frak n^{+}_{formal}$.

(d) $w\cdot (\Cal B^{-}\cdot 1\cdot \Cal B^{+})\cong (N^{-}\cap w
N^{+}w^{-1})\times (\Cal B^{-}\cdot w\cdot \Cal B^{+})$.

(e) The inclusion $G(A)\hookrightarrow \Cal G$ is a homotopy equivalence.
\endproclaim

\flushpar Proof.  Parts (a), (b) and (d) follow immediately from the
corresponding facts for $\Cal G^{-}$ and $\Cal G^{+}$.

To prove (c) recall that $\Cal P^{+}_{(j)}=G_{(j)}\cdot \Cal H\cdot
\Cal N^{+}$ is a parabolic subgroup of
$\Cal G^{+}$. This group has a Levi decomposition given by\newline
\centerline{$\Cal P^{+}_{(j)}=G_{(j)}\propto \Cal R^{+}_{(j)}$,}
where the radical $\Cal R^{+}_{(j)}$ has Lie algebra
$$Lie(\Cal R_{(j)}^{+})=(\sum_{k\ne j}\Bbb Ch_k\quad\oplus\quad\sum_{\matrix
\alpha >0\\\alpha\ne\alpha_j\endmatrix }\frak g_{\alpha})_{formal}
.$$
Now suppose that $x^{-}\in \frak n^{-},x^{+}\in \frak n^{+},$ and $
h\in \Cal H$. Let $w$ denote the $g_{\alpha_j}$
component of $x^{+}$, so that $exp(-we_j)\cdot exp(x^{+})\in \Cal R^{
+}_{(j)}$. Then for
$g=\left(\matrix a&b\\c&d\endmatrix \right)\in G_{(j)}$,
$$exp(x^{-})\cdot h\cdot exp(x^{+})\cdot g=exp(x^{-})\cdot g_1\cdot
r,$$
where $g_1=h\cdot exp(we_j)\cdot g\in G_{(j)},$ and
$$r=Ad_{g^{-1}}(exp(-we_j)\cdot exp(x^{+}))\in \Cal R^{+}_{(j)}.$$
Thus $exp(x^{-})\cdot h\cdot exp(x^{+})\cdot g$ will belong to the top stratum
if and
only if $g_1$ belongs to the top stratum, which is the case precisely
when $a+cw\ne 0$.  From this it is clear that each root space
component of log({\it r}) will be a rational function of finitely many root
space components of $x^{+}$, where the corresponding divisor is given by
$a+cw=0.$ Thus the right action by g on $\Cal G$, expressed in terms of
the coordinate $\frak n^{-}_{formal}\times \Cal H\times \frak n^{
+}_{formal}$ for the top stratum, defines a
rational map.  Since the subgroups $G_{(j)}$ generate {\it G}, it follows that
the same is true for any $g\in G$.  Since the right action by a $
g\in \Cal N^{+}$
defines a polynomial map, it follows that the right action by $\Cal G^{
+}$ is
rational in the above coordinate.

Part (e) is proved by the same argument used to prove (8.6.6) in
[PS].//

\bigskip

\flushpar\S 1.5. The formal flag space $\Cal F=\Cal G/\Cal B^{+}$.

We will refer to the space
$$\Cal F=\Cal G/\Cal B^{+}=\Cal G^{-}/B^{+}$$
as the formal completion of the flag space $G/B^{+}$, or the
formal flag space.  The following essentially follows from
(1.4.1).

\proclaim{(1.5.1) Proposition} (a) $\Cal F=\bigsqcup\,\,\Sigma_w$, where $
\Sigma_w$ is the
$\Cal B^{-}$-orbit of $[w]\in \Cal F$, and the disjoint union is over $
w\in W$.

(b) $\Cal N^{-}\cong\Sigma_1$. More generally, $\Cal N^{-}_w\cong
\Sigma_w$.

(c) $\Cal F$ is a $\Cal G^{-}$-homogeneous proalgebraic manifold, modelled on $
\frak n^{-}_{formal}$.

(d) $w\cdot\Sigma_1\cong (N^{-}\cap wN^{+}w^{-1})\times\Sigma_w$.

(e) The inclusion $G/B^{+}\hookrightarrow \Cal F$ is a homotopy equivalence.
\endproclaim

\flushpar(1.5.2) Remark. Together with (1.2.2), this shows that
the manifold $\Cal F$ satisfies the conditions I-VI in section 14.5 of
[PS].  The point of mentioning this is that, together with the
discussion in [PS], it illustrates how formal completions meld
analysis and algebra.  This is usually referred to as the ``gaga''
principle of algebraic geometry (see [Hart]).

\bigskip

\flushpar\S 1.6. The infinitesimal action of $\frak g$ on $\Cal F$.

The group $G(A)$ acts holomorphically on the space $\Cal F$.  Since
$G(A)$ is not in general a Lie group, it does not automatically
follow that there is an infinitesimal action of $\frak g$ on $\Cal F$.
However enough Lie algebra elements can be exponentiated, so
that the action can be pieced together.  The same applies to
the left and right actions of $G(A)$ on $\Cal G$.  We will write the
actions out explicitly.

Given $\xi\in \frak g$ such that there is a corresponding one parameter
group $\{e^{t\xi}\}\subset G(A)$ (e.g.  if $\xi =h_i,f_i,e_i$, or, more
generally, a real
root vector), there is a corresponding vector field $\hat{\xi}$ on $
\Cal F$
which is defined by
$$\hat{\xi}\vert_q=\frac d{dt}\vert_{t=0}e^{t\xi}\cdot q.\tag 1.6.1$$

To compute a formula for this vector field, we first
parameterize the top stratum of $\Cal F$ by $\Cal N^{-}$.  Since $
\Cal N^{-}$ is a Lie
group, we can use left translation to trivialize the tangent
bundle, and there is a resulting identification
$$Functions(\Cal N^{-},\frak n^{-}_{formal})\to Vect(\Cal N^{-}).\tag 1.6.2$$
We write $g=l\cdot diag\cdot u$ for the Bruhat factorization of a $
g$ in
the top stratum of $\Cal G$.  We then have
$$\aligned
\hat{\xi}\vert_{l_0}&=\frac d{dt}\vert_{t=0}(l^{-1}_0l(e^{t\xi}l_
0))=proj_{\frak n^{-}_{formal}}Ad_{l_0}^{-1}(\xi )\endaligned
\tag 1.6.3$$
for $l_0\in \Cal N^{-}$.  This formula makes sense, because if
$heights(\xi )\le M$, then it follows that $heights(Ad_{l_0}^{-1}
\xi )\le M$.  From
this formula we can see that there is an infinitesimal action
$\frak g\to Vect(\Cal F)$, and (1.6.4) determines $\hat{\xi}$ for a general $
\xi\in \frak g$.

In \S 3.2 below we will need another formula for $\hat{\xi}$. Recall that
$$exp:\frak n^{-}_{formal}\to \Cal N^{-}:x\to l_0=exp(x)$$
is an isomorphism of manifolds. In terms of the linear space
$\frak n^{-}_{formal}$, we have
$$\aligned
\hat{\xi}\vert_x&=\frac d{dt}\vert_{t=0}(log(l(e^{t\xi}e^x)))\\&=
dlog\vert_{e^x}\circ dL_{e^x}(\frac d{dt}\vert_{t=0}l(e^{te^{-adx}
(\xi )}))\\&=Td(adx)(proj_{\frak n^{-}_{formal}}(e^{-adx}(\xi )))\endaligned
\tag 1.6.5$$
where $Td(y)$ denotes the Todd series (for $\vert y\vert <2\pi$,
$Td(y)=\frac y{1-e^{-y}})$.

It is also a simple matter to write down the infinitesimal left
and right actions of $\frak g$ on $\Cal G$.  We parameterize the top stratum
by $\Cal N^{-}\times \Cal H\times \Cal N^{+}$, and we identify vector fields on
$
\Cal N^{\pm}$ using the left
trivialization.  Then for $\xi\in \frak g$, the vector field corresponding to
the left action is given on the top stratum by
$$\hat{\xi}\vert_{(l_0,diag_0,u_0)}=\frac d{dt}\vert_{t=0}(l^{-1}_
0l(e^{t\xi}l_0),diag(e^{t\xi}l_0),u_0^{-1}u(e^{t\xi}g_0))=$$
$$=(proj_{\frak n^{-}_{formal}}(Ad_{l_0}^{-1}\xi ),proj_{\frak h}
(Ad_{l_0}^{-1}\xi ),Ad^{-1}_{u_0}(proj_{\frak n^{+}_{formal}}(Ad_{
ldiag}^{-1}\xi )))\tag 1.6.6$$
where $g_0=l_0\cdot diag_0\cdot u_0$.

\bigskip

\centerline{\S 2. Matrix Coefficients}

\bigskip

\flushpar\S 2.1. The Kac-Peterson version of the Peter-Weyl theorem.

Let $(\pi_{\Lambda},$ $L(\Lambda )$) denote the integrable highest weight
representation
corresponding to a dominant integral functional $\Lambda$, and let $
(\pi^{*}_{\Lambda},L^{*}(\Lambda ))$
denote the dual integrable lowest weight representation ($L^{*}(\Lambda
)$ is the
sum of the weight spaces of the topological dual of $L(\Lambda )$). Given
$v^{*}\otimes v\in L^{*}(\Lambda )\otimes L(\Lambda ),$ the matrix coefficient
defined on $
\Cal G^{-}\times \Cal G^{+}$ by
$$f_{v^{*}\otimes v}(g,h)=(\pi_{\Lambda}(h)\cdot v,\pi^{*}_{\Lambda}
(g^{-1})\cdot v^{*})\tag 2.1.1$$
descends to a well-defined regular function on $\Cal G$. The map
$$\Phi :L^{*}(_{}\Lambda )\otimes L(\Lambda )\to \Cal O_{reg}(\Cal G
):v^{*}\otimes v\to f_{v^{*}\otimes v}\tag 2.1.2$$
is a $\Cal G^{-}\times \Cal G^{+}$-equivariant map. The Kac-Peterson version of
the Peter-Weyl Theorem can be formulated as the following

\proclaim{(2.1.3) Proposition}  The map $\Phi :\sum\oplus L^{*}(\Lambda
)\otimes L(\Lambda )\to \Cal O_{reg}(\Cal G)$
is a $\Cal G^{-}\times \Cal G^{+}$-equivariant isomorphism.  \endproclaim

Suppose that $k$ is a holomorphic character of $C_0$, the identity
component of the center of $G(A)$.  We will say that $\Lambda$ has
level $k$ if
$$\pi_{\Lambda}:C_0\to \Bbb C^{*}:\lambda\to\lambda^k.\tag 2.1.4$$
Also, let $\Cal L_k^{*}\to \Cal G/C_0$ denote the induced bundle, where
$$\Cal O(\Cal L_k^{*}\to \Cal G/C_0)\cong \{f\in \Cal O(\Cal G)\vert
f(g\cdot\lambda )=\lambda^kf(g),\forall\lambda\in C_0\}.\tag 2.1.5$$

\proclaim{(2.1.6)Corollary}Suppose that $n$, the size of the Cartan
matrix, is finite. Then the induced map
$$\Phi_k:\sum_{level(\Lambda )=k}\oplus L^{*}(\Lambda )\otimes L(
\Lambda )\to \Cal O(\Cal L^{*}_k\to \Cal G/C_0)$$
is a $\Cal G^{-}\times \Cal G^{+}$-equivariant isomorphism. \endproclaim

\flushpar Proof.  This will follow from (2.1.3), once we know
that any holomorphic section of $\Cal L_k^{*}\to \Cal G/C_0$ is actually
regular.
By (1.2.1) such a section $\sigma$ must depend upon finitely many
variables, in the sense that for a sufficiently large $M$,
the corresponding function $f$, as in (2.1.5), is actually a
function on the finite dimensional double coset space
$$N^{-}_M\backslash G(A)/N^{+}_M=\frak N^{-}_M\backslash \Cal G/\frak N^{
+}_M,\tag 2.1.7$$
where $N^{+}_M$ is the subgroup corresponding to positive roots
dominating $M$, as in \S 1.3.
We can interpret $\sigma$ as a holomorphic section of a bundle
over a finite dimensional complex projective variety
$$(N^{-}_M\backslash G(A)/C_0N^{+}_M)_{ss}/equivalence,\tag 2.1.8$$
which contains
$$(N^{-}_M\backslash N^{-})\times (H/C_0)\times (N^{+}/N^{+}_M)\tag 2.1.9$$
as a dense open set.  Hence $\sigma$ must actually be regular.   In
more detail, $(\cdot )_{ss}$ is the set of semistable points, i.e.  the
points $q$ for which there exists a regular function $\sigma$ on $
G(A)$ of
level $k$ with $\sigma (q)\ne 0$ and with $\sigma (lgu)=\sigma (g
)$ for $l\in N^{-}_M$ and
$u\in N^{+}_M$.  Equivalence of a pair of semistable points means that
the pair is not separated by a regular section satisfying these
conditions.  A projective embedding of (2.1.8) is obtained simply
by mapping a point $q$ in (2.1.8) to the line of the evaluation
functional on the space of regular $\sigma$ of level $k$ stabilized by
$N^{-}_M$ and $N^{+}_M$.  The fact that (2.1.9) injects into this variety is
essentially a restatement of the fact that there is a unique
Bruhat factorization for the points of the top stratum of
$G(A)/Ker(k)$, and the components of the factorization are
expressible in terms of rational functions.  If we write $g$ in
the top stratum of $\Cal G$ as $g=l\cdot diag\cdot u$,
$$l=exp(\sum_{\alpha >0}x_{-\alpha}),\quad diag=\prod\sigma_j(g)^{
h_j},\quad u=exp(\sum_{\alpha >0}x_{\alpha}),$$
then the $x_{\pm\alpha}$ are meromorphic functions, expressible in terms
of ratios of regular sections, relative to the choice of basis
for each root space.  For $height(\alpha )\le M$, the sections
corresponding to $x_{\pm\alpha}$ are invariant with respect to $\Cal N^{
-}_M$ and $\Cal N^{+}_M$
(since for example multiplying $g$ on the left by an $l_1\in \Cal N^{
-}_M$ does
not change the $\frak g_{-\alpha}$ component of $log(l)$).  Since when we pass
to the quotient $\Cal G/Ker(k)$, the Bruhat factorization remains
unique, it follows that there are enough of these sections of
level $k$ to separate these points.  //

\bigskip

\flushpar\S 2.2. The Kac-Peterson version of the Borel-Weil theorem.

Given a dominant integral functional $\Lambda$, define characters
$\chi_{\Lambda}^{\pm}$ of $\Cal B^{\pm}$ by
$$\chi_{\Lambda}^{\pm}(exp(h)\Cal N^{\pm})=exp(\pm\Lambda (h)).\tag 2.2.1$$
If $v_{\Lambda}$, resp.  $v^{*}_{\Lambda}$, is a choice of highest, resp.
lowest, weight
vector for $\pi_{\Lambda}$, resp.  $\pi_{\Lambda}^{*}$, then
$$\chi_{\Lambda}^{+}(b^{+})v_{\Lambda}=\pi_{\Lambda}(b^{+})\cdot
v_{\Lambda},\quad resp.\quad\chi_{\Lambda}^{-}(b^{-})v^{*}_{\Lambda}
=\pi_{\Lambda}^{*}(b^{-})\cdot v^{*}_{\Lambda},\tag 2.2.2$$
for $b^{\pm}\in \Cal B^{\pm}$.

Let $\Cal L_{\chi^{+}_{\Lambda}}\to (\Cal G/\Cal B^{+}=\Cal G^{-}
/B^{+})$ and $\Cal L_{\chi_{\Lambda}^{-}}\to (\Cal B^{-}\backslash
\Cal G=B^{-}\backslash \Cal G^{+})$ denote the
associated bundles with isotropy representations $\chi_{\Lambda}^{
\pm}$,
respectively, i.e.
$$\Cal L_{\chi^{+}_{\Lambda}}=\Cal G^{-}\times_{\chi_{\Lambda}^{+}}
\Bbb C\quad and\quad \Cal L_{\chi^{-}_{\Lambda}}=\Bbb C\times_{\chi_{
\Lambda}^{-}}\Cal G^{+},\tag 2.2.3$$
where $(g\cdot b^{+},\lambda )\sim (g,\chi^{+}_{\Lambda}(b^{+})\cdot
\lambda )$ and $(\chi_{\Lambda}^{-}(b^{-})^{-1}\lambda ,b^{-}\cdot
g)$,
respectively.

\proclaim{(2.2.4)Corollary}As $\Cal G^{\pm}$-modules
$$\Cal O(\Cal L^{*}_{\chi^{-}_{\Lambda}}\to B^{-}\backslash \Cal G^{
+})\cong L(\Lambda )\quad and\quad \Cal O(\Cal L^{*}_{\chi^{+}_{\Lambda}}
\to \Cal G^{-}/B^{+})\cong L^{*}(\Lambda ),$$
respectively.\endproclaim

\flushpar Proof. By the Peter-Weyl theorem
$$\Cal O(\Cal L^{*}_{\chi_{\Lambda}^{+}})\tilde {=}\{f\in \Cal O_{
reg}(\Cal G)\vert f(g\cdot b^{+})=\chi^{+}_{\Lambda}(b^{+})\cdot
f(g),\forall g\in \Cal G,b^{+}\in \Cal B^{+}\}$$
$$\cong L^{*}(\Lambda )\otimes \Bbb Cv_{\Lambda}\cong L^{*}(\Lambda
),$$
as $\Cal G^{-}$-modules. Similarly
$$\Cal O(\Cal L^{*}_{\chi^{-}_{\Lambda}})\cong \{f\in \Cal O_{reg}
(\Cal G)\vert f((b^{-})^{-1}\cdot g)=\chi_{\Lambda}^{-}(b^{-})\cdot
f(g),\forall g\in \Cal G,b^{-}\in \Cal B^{-}\}$$
$$\cong \Bbb Cv_{\Lambda}^{*}\otimes L(\Lambda )\cong L(\Lambda )
,$$
as $\Cal G^{+}$-modules.//

\bigskip

\flushpar\S 2.3. Linear Borel embeddings of $\Cal G$.

Let $L$ denote a direct sum of integrable highest weight modules.
Choose a basis for $L$, say $\{v_i\}$, which is compatible with the direct
sum and weight space decompositions of $L$, and let $\{v_i^{*}\}$ denote the
dual basis.  Let $\Bbb M$ denote the space of matrices indexed by $
\{i\}$, which
we identify with $\prod \Bbb Cv_i\otimes v_j^{*}$.

\proclaim{(2.3.1) Proposition} (a) The map
$$\imath :\Cal G\mapsto \Bbb M:g\mapsto (g_{ij}),$$
where $g_{ij}=f_{v_i^{*}\otimes v_j}(g)$, is regular.

(b) If $L$ is the direct sum of the fundamental highest weight
modules, then $\imath$ is a Borel isomorphism of $\Cal G$ onto its image
$\imath (\Cal G)$.\endproclaim

\flushpar Proof.  Part (a) is obvious.  The injectivity of the
map can be proved by mimicking the method used by Kac and
Peterson (to prove that $G(A)$ is an algebraic group in the sense
of Shafarevich) in [KP3].  The fact that the map is a Borel
isomorphism follows automatically, since it is an injective
continuous map into a standard Borel space.//

\smallskip

The simplest examples show that the map $\imath$ is not in general a
homeomorphism; see (2.1.3) of Part II.

\bigskip

\flushpar\S 2.4. A digression on completions of $K(A)$ and $G(A)$.

In the case that $A$ is affine, it is well-known how to complete
$\frak g$ so that there exists a Lie group completion of $G(A)$ which
has $\frak g$ as a dense subalgebra of its Lie algebra.  This can be
done in many ways, so that in this case, there exist many
interesting group-theoretic completions of $G(A)$ contained in
$G(A)_{formal}$.  In this subsection we will briefly indicate a
practical criterion for existence of such Lie completions in
nonaffine cases.  The idea is to observe that the Toda
equations, as formulated by Kostant, make sense for a general
Kac-Moody algebra, and the existence of a Lie completion
implies that the solutions are necessarily meromorphic.

Let $\langle\cdot ,\cdot\rangle$ denote the essentially unique nondegenerate
symmetric bilinear form on $\frak g$.  Assume that there exists a
complex Lie group $\bar {G}$ which has a Lie algebra $\bar {\frak g}$
containing $
\frak g$ as
a dense subalgebra.  Without further comment, we will denote
the completion of $\frak n^{-}$ in $\bar {\frak g}$ by $\bar {\frak n}^{
-}$, and so on.  We will assume
that $\bar {G}$ has the following properties.  First, the form $\langle
\cdot .\cdot\rangle$
extends continuously in each variable to $\bar {\frak g}$ (this is certainly
true for all the completions constructed in [RC]).  Second,
$$exp:\bar {\frak g}\to\bar {G}$$
is well-defined and holomorphic. Third, $\bar {G}$ has a stratification
parameterized by the elements of the Weyl group $W=N(H)/H$,
namely
$$\bar {G}=\bigsqcup_{w\in W}\bar{\Sigma}_w,\quad\bar{\Sigma}_w=\bar {
B}^{-}w\bar {B}^{+}.$$
And finally for the top stratum $\bar{\Sigma}_1$, there is a diffeomorphism
$$\bar{\Sigma}_1=\bar {N}^{-}\times H\times\bar {N}^{+},$$
and the corresponding Gaussian factorization is given by
$$g=l\cdot diag\cdot u$$
where
$$diag(g)=\prod\sigma_j(g)^{h_j},$$
and the $\sigma_j$ are the fundamental matrix coefficients.

These assumptions are obviously valid for the $C^{\infty}$ completion
of an affine Lie algebra.

\smallskip

We now recall the Kostant formulation of the generalized Toda
equations.

Since the form $\langle\cdot ,\cdot\rangle$ extends continuously to $
\bar {\frak g}$, there
is an induced map with dense image,
$$\bar {\frak b}^{-}\to (\bar {\frak b}^{+})^{*}:x\to\langle x,\cdot
\rangle .\tag 2.4.1$$
In a formal sense this induces a Poisson structure on $\bar {\frak b}^{
-}$.
Define
$$\epsilon =\sum e_j,\quad M=\epsilon +\bar {\frak b}^{-}.\tag 2.4.2$$
Via the identification
$$M=\epsilon +\bar {\frak b}^{-}\cong\bar {\frak b}^{-},\tag 2.4.3$$
there is an induced Poisson structure on $M$.

Suppose that $x(t)$ is a curve in $M$. The complex Kostant-Toda
equations are given by
$$\dot {x}(t)=[x(t),proj_{\bar {\frak n}^{-}}(x(t))],\tag 2.4.4$$
where $proj_{\bar {\frak n}^{-}}$ has kernel $\bar {\frak b}^{+}$.  The real
Kostant-Toda equations
are the equations obtained by restricting to the normal real
form.  Since we are interested in analyticity properties of the
solutions, we will stick exclusively with the complex case.
The basic facts about these equations for finite type $A$ carry
over directly to this more general setting, namely

\proclaim{(2.4.5)Proposition} (a) The Kostant-Toda equations are
Hamiltonian for $M$. More precisely, define
$$\Cal H:M\to \Bbb C:x\to\frac 12\langle x,x\rangle ;$$
then for any smooth function $f:M\to \Bbb C$,
$$\frac d{dt}f(x(t))=\{\Cal H,f\}\vert_{x(t)}$$
for any solution $x(\cdot )$ of (2.4.4).

(b) Given an initial condition $x_0$ with $e^{x_0}\in\bar{\Sigma}_
1$,
$$x(t)=l(t)^{-1}x_0l(t)\tag 2.4.6$$
where
$$e^{tx_0}=l(t)diag(t)u(t)$$
solves (2.4.4), up to the time when $e^{tx_0}$ exits $\bar{\Sigma}_
1$.

(c) Suppose that $x_0\in\epsilon +\frak b^{-}$, and suppose $heig
hts(x_0)\ge -n$ (in the
sense that each homogeneous piece of $x_0$ has $height\ge -n$). Then
$$heights(x(t))\ge -n$$
for $x(\cdot )$ in (b). In particular, suppose that
$$x_0=\epsilon +\sum (a_jh_j+b_jf_j);$$
then the equations (2.4.4) are equivalent to
$$\aligned
\dot {a}_j&=b_j\\\dot {b}_j&=\sum_ia_ia_{ij}b_j\endaligned
\tag 2.4.7$$
\endproclaim

\flushpar Brief proof of (b) and (c). Directly from the definition
(2.4.6), we see that
$$\dot {x}=[x,l^{-1}\dot {l}].$$
Let $b=diag\cdot u$. Since
$$lb=e^{tx_0}$$
it follows that
$$\dot {l}b+l\dot {b}=e^{tx_0}x_0\implies l^{-1}\dot {l}+\dot {b}
b^{-1}=bx_0b^{-1}=x$$
Thus we see that
$$proj_{\frak n^{-}}x=l^{-1}\dot {l}\quad and\quad proj_{b^{+}}x=
\dot {b}b^{-1}.$$
This proves (b).

In terms of the solution, part (c) follows immediately from the
observation that
$$x(t)=l^{-1}x_0l=bx_0b^{-1},$$
since
$$heights(bx_0b^{-1})\ge -n.$$
In terms of the equations themselves,
$$\dot {x}=[x,proj_{\frak n^{-}}x]=[proj_{\frak h}x+\epsilon ,pro
j_{\frak n^{-}}x].$$
The right hand side clearly has $heights\ge -n$, so all the
homogeneous components of $x$ of $height<-n$ are constant, hence
identically zero.

This completes the proof.//

\smallskip

The equations (2.4.7) were studied by Zhou in his thesis in the
simplest hyperbolic cases (see [Z]).  He observed that the
solutions are not meromorphic.  This implies that in these
cases there does not exist a Lie completion of $G(A)$ (having the
properties we ascribed to $\bar {G}$ above).  From an analytical point
of view, this suggests that these algebras are perhaps
comparable to the Virasoro algebra, in the sense that there
exists a reasonable global unitary real form (see [RC]), but
there does not exist a corresponding global complexification.

\bigskip

\centerline{\S 3. Measures on the Formal Flag Space}

\bigskip

\noindent\S 3.1.  The top stratum has full measure.

The following result draws a purely technical conclusion, but
the proof involves a key idea.

\proclaim{(3.1.1) Proposition} Suppose that $\mu$ is a
$K$-quasi-invariant probability measure on the formal completion
of the flag space, $\Cal F$.  Then $\mu (\Sigma_1)=1$, i.e.  $\mu$ is supported
on the
top stratum.  \endproclaim

\flushpar Proof. Fix $1\le j<n+1$. We will exploit the following diagram:
$$\matrix \Cal G/\Cal B^{+}&\mapsto&\Cal R_{(j)}^{-}\backslash \Cal G
/\Cal B^{+}\\\uparrow&&\uparrow\\\Sigma_1\cup (r_j\cdot\Sigma_1)&
\mapsto&G_{(j)}/B_{(j)}^{+}\\\uparrow&&\uparrow\\\Sigma_1\tilde {
=}\Cal N^{-}&\mapsto&N_{(j)}^{-}.\endmatrix \tag 3.1.2$$
The top arrow is $G_{(j)}-$equivariant, because $\Cal R^{-}_{(j)}$, the radical
of the
parabolic subgroup $\Cal P^{-}_{(j)}=G_{(j)}\cdot \Cal B^{-}$, is normalized by
$
G_{(j)}$.  All of the
up arrows are injective.

The projection of $\mu$ to the space $\Cal R_{(j)}^{-}\backslash
\Cal F$ is $K_{(j)}$-quasi-invariant;
therefore the restriction of this projection to the sphere $G_{(j
)}/B_{(j)}^{+}$
belongs to the unique invariant measure class.  It follows that
$\Sigma_1=r_j\cdot\Sigma_1$, modulo sets of $\mu$-measure zero, hence that $
\Sigma_1=w\cdot\Sigma_1$,
modulo sets of $\mu$-measure zero, $\forall w\in W$.  This implies that $
\Sigma_1$ has full
measure.//

\bigskip

\flushpar\S 3.2. On uniqueness of $K$-invariant measures on $\Cal F$.

Let $\Lambda$ denote a dominant integral functional, and let $\Cal L^{
*}=$
$\Cal L^{*}_{\chi^{+}_{\Lambda}}\to \Cal F$ denote the corresponding positive
line bundle.  Recall
that as $\Cal G^{-}$-modules
$$\Cal O(\Cal L^{*})\tilde {=}L^{*}(\Lambda ).$$
Let $\sigma =\sigma_{\Lambda}$ denote a regular section which is of lowest
weight.

\proclaim{(3.2.1) Conjecture} For each $k\ge 0$, there is at most
one measure $\mu^{\vert \Cal L\vert^{2k}}$ on $\Cal F$ having values in $
(\Cal L\otimes\bar {\Cal L})^k$ such that (i)
$\mu^{\vert \Cal L\vert^{2k}}$ is $K$-invariant, and (ii) $d\mu_k
=$ $(\sigma\otimes\bar{\sigma })^kd\mu^{\vert \Cal L\vert^{2k}}$ is a
probability measure.  \endproclaim

\proclaim{(3.2.2) Remarks} (i) This conjecture can be reformulated in
the following way:  there is at most one $K$-quasi-invariant
probability measure $\mu_k$ on $\Sigma_1$ satisfying
$$\frac {d\mu_k(g^{-1}\cdot x)}{d\mu_k(x)}=\vert\frac {g_{*}\sigma}{
\sigma}\vert^{2k}(x)=\vert\frac {g\cdot\sigma (g^{-1}\cdot x)}{\sigma
(x)}\vert^{2k},$$
for all $g\in K$.  This follows from (3.1.1) above.

(ii) If $g=\left(\matrix a&b\\c&d\endmatrix \right)\in K_{(j)}$, and $
x=exp(zf_j+...)\cdot \Cal B^{+}$, where
$zf_j+...\in n^{-}_{formal},$ then (viewing $\sigma$ now as a function on $
\Cal G$)
$$\vert\frac {\sigma (g\cdot x)}{\sigma (x)}\vert^{}=\vert\sigma
(g\cdot exp(zf_j+...)\vert =\vert a+bz\vert^{\Lambda (h_j)};$$
for $\Cal N^{-}\cdot g\cdot exp(zf_j+...)\cdot \Cal N^{+}=\Cal N^{
-}\cdot exp(log(a+bz)h_j)\cdot \Cal N^{+}$.

This implies that for any $g\in K$, the Radon-Nikodym derivative
in (i) is a function of finitely many variables.

\endproclaim

\flushpar Idea of the Proof.  Let $\mu$ denote a $K$-invariant measure
with values in $\vert \Cal L\vert^{2k}$ such that $d\mu_k=\vert\sigma
\vert^{2k}d\mu$ is a probability
measure.  We will view $\mu_k$ either as a measure on $_{}n^{-}_{
formal}$ or
$\Cal N^{-}$, via the isomorphisms $exp:_{}n^{-}_{formal}\mapsto
\Cal N^{-}$, and $\Cal N^{-}\tilde {=}\Sigma_1$,
depending upon which is most convenient.   The diagram (3.1.2)
implies that the image of $\mu$ in each of the spheres $G_{(j)}/B^{
+}_{(j)}$ is
uniquely determined; in terms of the parameterization
$z\to exp(zf_j)\cdot B_{(j)}^{+}$ of $\Bbb P^1\tilde {=}G_{(j)}/B^{
+}_{(j)}$, the image of $\mu_k$ is given
by the formula
$$\frac 1E(1+\vert z\vert^2)^{-2-k\Lambda (h_j)}\frac {dz\wedge d
\bar {z}}{2i}.\tag 3.2.3$$

Suppose that $\nu$ is another measure as in the Proposition.
Without loss of generality, we can assume that $\nu <<\mu$.  For
the convex combination $\frac 12(\nu +\mu )$ is another measure as in the
Proposition, and this latter measure dominates $\nu$ and $\mu$; hence,
if necessary, we could replace $\mu$ by $\frac 12(\nu +\mu )$.  Thus we can
assume that $d\nu =\rho d\mu$, where $\rho$ is a nonnegative $K(A
)$-invariant
scalar function with $\int\rho d\mu_k=1$.  By replacing $\rho$ by an
appropriate multiple of $\rho^{1/2}$, we can suppose that $\rho$ is in
$L^2(d\mu_k)$.

By the first paragraph we know that
$$\int f\rho d\mu_k=\int fd\mu_k,\tag 3.2.4$$
for all $L^2$ functions of the form
$$f(x)=\sum_jF_j(x_{-\alpha_j}),\quad for\quad x=\prod_{\alpha >0}
x_{-\alpha}\in \frak n^{-}_{formal}.\tag 3.2.5$$

The idea now is to use infinitesimal methods to show that we
can generate enough functions to force $\rho$ to be 1.  For
simplicity we will consider the case $k=0$.  Recall from \S 1.6
that $\frak g$ acts by holomorphic vector fields on $\Cal F$.  For each $
M$,
we consider the natural inclusion
$$\Cal S(\Cal N^{-}/\Cal N^{-}_M)\subset L^2(\Cal N^{-},d\mu_0),\tag 3.2.6$$
where $\Cal S$ denotes the Schwarz space (and note that for a
nilpotent group, there is no ambiguity about the meaning of a
Schwarz class function).  The induced action of $\frak b^{-}$ on functions
on $\Cal N^{-}$ restricts to the natural action of $\frak b^{-}/\frak n^{
-}_M$ on $\Cal S(\Cal N^{-}/\Cal N^{-}_M)$.
The explicit formula (1.6.4) shows that $\hat {e}_i$ maps $\Cal S
(\Cal N^{-}/\Cal N^{-}_M)$ to
$\Cal S(\Cal N^{-}/\Cal N^{-}_{M+1})$.  Thus there is an induced action
$$\frak g\times (\Cal S_{\infty}\equiv\lim_{M\to\infty}\Cal S(\Cal N^{
-}/\Cal N^{-}_M))\to \Cal S_{\infty}.\tag 3.2.7$$

If $\xi\in \frak k$ and there exists a one-parameter group corresponding to
$\xi$ in $K$, then we can assert that for $f\in \Cal S_{\infty}$ and $
f\perp\rho -1$,
$\hat{\xi}\cdot f\perp\rho -1$.  Since such vector fields generate $
\frak k$, it follows that
$$\Cal S_{\infty}\cap (\rho -1)^{\perp}\tag 3.2.8$$
is invariant under the natural action of $\Cal U(\frak g)$, the universal
enveloping algebra. Thus it suffices to show that the Schwarz
class functions of the form (3.2.5) are cyclic for the action of $
\Cal U(\frak g)$
on $\Cal S_{\infty}\subset L^2$. This seems very plausible to me, but I cannot
prove it.

\bigskip

\flushpar\S 3.3. A compactness result.

Suppose that $X$ is a complete separable metric space. Let
$BC(X)$ denote the space of bounded continuous functions with
the uniform topology, and equip $Prob(X)\subset BC(X)^{*}_{\le 1}$ with the
$weak^{*}$-topology. Prohorov's theorem characterizes compactness
in this space (see section 6 of [Bill1]):

\proclaim{(3.3.1)Proposition}$\Pi\subset Prob(X)$ is relatively compact if
and only if $\Pi$ is tight, i.e.  for all $\epsilon >0$, there exists a
compact set $K_{\epsilon}\subset X$ such that $\nu K_{\epsilon}>1
-\epsilon$, for all
$\nu\in\Pi$.\endproclaim

We will be applying this to formal completions, in which case
our space $X$ is equivalent to $\Bbb R^{\infty}$.  In this case $
\Pi\subset Prob(\Bbb R^{\infty})$
is tight if and only if for each finite dimensional projection
$p:\Bbb R^{\infty}\to \Bbb R^n$, $p_{*}\Pi\subset Prob(\Bbb R^n)$ is tight, if
and only if for each
coordinate function
$$p_j:\Bbb R^{\infty}\to \Bbb R:x\to x_j,$$
$(p_j)_{*}\Pi\subset Prob(\Bbb R)$ is tight.

Suppose that $\{\nu_l\}$ is a sequence of measures with values in the
bundle $\vert \Cal L_{\Lambda}\vert^{2k}\to \Cal F$.  We will consider
sequences satisfying the
condition
$$\vert\int\Phi d(g_{*}\nu_l)-\int\Phi d\nu_l\vert\to 0\quad as\quad
l\to\infty\tag 3.3.2$$
for each $g\in K(A)$ and for each density $\Phi$ of the form
$$\Phi =f\cdot\vert\sigma_{\Lambda}\vert^{2k},\tag 3.3.3$$
where $f$ is a bounded Borel function of finitely many variables.

If $k=0$, so that $\nu_l$ is an ordinary measure, then the condition
(3.3.2) implies that for $g\in K(A)$ and each $f=\chi_E$,
$$\vert\nu_l(g^{-1}E)-\nu_l(E)\vert\to 0\quad as\quad l\to\infty
.\tag 3.3.4$$
In this case the meaning of the condition is clear.  Suppose
that $k>0$, and let $\mu_{l,k}$ denote the measure
$$\mu_{l,k}=\vert\sigma_{\Lambda}\vert^{2k}d\nu_l.\tag 3.3.5$$
Then we have
$$dg_{*}\mu_{l,k}=\vert g_{*}\sigma_{\Lambda}\vert^{2k}dg_{*}\nu_
l.\tag 3.3.6$$
If we take $f=\chi_{gE}$, then (3.3.2) reads
$$\vert\int_E\vert\frac {g_{*}\sigma_{\Lambda}}{\sigma_{\Lambda}}
\vert^{2k}d\mu_{l,k}-\mu_{l,k}(gE)\vert\to 0\quad as\quad l\to\infty
.\tag 3.3.7$$

The example that motivates (3.3.2) is the family of Wiener
measures on a (possibly twisted) loop group.  Here $l$ is
interpreted as the temperature.  The family of Wiener
measures is in fact approximately invariant in a much stronger
sense, namely they satisfy a uniform version of (3.3.4); see
\S 4.1 of Part III.

\proclaim{(3.3.8)Proposition} Suppose that the measures $\nu_l$
(having values in the line bundle $\vert \Cal L\vert^{2k}$) are normalized so
that
$$d\mu_{l,k}=\vert\sigma_{\Lambda}\vert^{2k}d\nu_l$$
is a probability measure for each $l$.  If the sequence $\{\nu_l\}$
satisfies the condition (3.3.2), then $\{\mu_{l,k}\}$ is a tight set of
probability measures on $\Cal F$.\endproclaim

\flushpar Proof.  For each $j$ let $\Cal P_{(j)}^{+}$ denote the parabolic
subgroup $G_{(j)}\cdot \Cal B^{+}$.  We will use the parameter $x
\in \frak n_{formal}^{-}$ to
parameterize the top stratum of $\Cal F=\Cal G/\Cal B^{+}$, i.e.
$$x\to exp(x)\cdot \Cal B^{+}\in \Cal G/\Cal B^{+}.\tag 3.3.9$$
Similarly we will use the parameter $x_{(j)}\in \frak n_{(j)}^{-}$ to
parameterize
$\Cal G/\Cal P_{(j)}^{+}$, where
$$\frak n^{-}_{(j)}=\frak n^{-}_{formal}\ominus \frak g_{-\alpha_
j}\tag 3.3.10$$
We are principally interested in the $x$ distributions of the
measures $\mu_{l,k}$.  However it is convenient to also consider the
$x_{(j)}$ distributions, because $K_{(j)}$ acts linearly on this variable,
through its natural action on the partial flag space $\Cal G/\Cal P^{
+}_{(j)}$.  In
particular the $K_{(j)}$ action will take a bounded function of $
x_{(j)}$
to another bounded function.

Consider the $G_{(j)}$ equivariant projection in (3.1.2).  The
condition (3.3.2) implies that the $l\to\infty$ limit of the image of
$\mu_{l,k}$ in the sphere $G_{(j)}/B^{+}_{(j)}$ is given by (3.2.3).  This
implies
that the projection of $\{\mu_{l,k}\}$ into
$$\prod_{0<height(\alpha )\le m}\frak g_{-\alpha}\subset \frak n^{
-}\tag 3.3.11$$
is a tight family, in the case $m=1$.  We will now show that
this is true for any $m$ by induction.

Assume that tightness has been established for $m-1$.

Suppose that $\alpha$ is a root and $0<height(\alpha )<m$.  With respect
to the ordinary linear adjoint action, consider the smallest
$\frak g_{(j)}$-invariant subspace containing a $v\in \frak g_{-\alpha}$.  As
we vary $
j$, $\alpha$
and $v$, the span of these subspaces will contain all root spaces
of $height=-m$, because $n^{-}$ is generated by $\{f_i\}$.  Thus it
suffices to establish tightness in the part of $\frak g_{(j)}\cdot
v$ of
$height=-m$. We can assume that $\alpha\ne\alpha_j$, for otherwise there is
nothing to prove.

Now $x$ and $x_{(j)}$ are related by
$$exp(x)=exp(x_{(j)})exp(zf_j)\tag 3.3.12$$
where $zf_j$ is the $\frak g_{-\alpha_j}$ component of $x$. It follows from
this that
for a positive root $\beta$,
$$proj_{\frak g_{-\beta}}(x)=proj_{\frak g_{-\beta}}(x_{(j)})+p,\tag 3.3.13$$
where $p$ is a polynomial in $z$ and the components of $x_{(j)}$ in
root spaces $\frak g_{-\gamma}$ with $0>height(-\gamma )>height(-
\beta )$. Similarly
$$proj_{\frak g_{-\beta}}(x_{(j)})=proj_{\frak g_{-\beta}}(x)+q\tag 3.3.14$$
where $q$ is a polynomial in the components of $x$ in root spaces
$\frak g_{-\gamma}$ with $0>height(-\gamma )>height(-\beta )$. The basic idea
of the proof
is to exploit these identities, which will allow us to work in a
setting where $K_{(j)}$ acts linearly.

It follows from (3.3.13) and the induction hypothesis that for
the $x_{(j)}$ distribution of the measures $\mu_{l,k}$, the projection into
$$(\prod_{0<height(\gamma )<m}\frak g_{-\gamma})\ominus \frak g_{
-\alpha_j}\subset \frak n^{-}_{(j)}\tag 3.3.15$$
is tight.  In particular the projection (relative to $\langle\cdot
,\cdot\rangle$) of the
$x_{(j)}$ distribution into $\Bbb Cv$ is tight.  We now claim that the
condition (3.3.2) implies that for the $x_{(j)}$ distribution of the
measures $\mu_{l,k}$, the projection into
$$span\{K_{(j)}\cdot v\}=\frak g_{(j)}\cdot v\tag 3.3.16$$
is tight (and recall here that the action of $K_{(j)}$ on $\Cal G
/\Cal P^{+}_{(j)}$,
expressed in the coordinate $x_{(j)}$, coincides with its $\underline {
linear}$
adjoint action on $\frak n^{-}_{(j)}$, so that there is no ambiguity about the
meaning of the dot in (3.3.16)).  It will follow from this claim
and (3.3.14) that for the $x$ distribution of the measures $\mu_{
l,k}$, the
projection into the part of $\frak g_{(j)}\cdot v$ of $height=-m$ is tight.
This
will complete the induction.

We first consider the case $k=0$, so that $\nu_l=\mu_{l,k}$.  Fix
$g\in K_{(j)}$.  We know that given $\epsilon >0$ there exists a compact
subset $C_{\epsilon}$ of $\Bbb Cv$ such that
$$\nu_l(\{proj_{\Bbb Cv}(x_{(j)})\in C_{\epsilon}\})>1-\epsilon ,
\quad\forall l.\tag 3.3.17$$
Now
$$\aligned
\nu_l(\{proj_{\Bbb Cg\cdot v}(x_{(j)})\in g\cdot C_{\epsilon}\})=
g_{*}\nu_l(\{proj_{\Bbb Cv}(x_{(j)})\in C_{\epsilon}\})\endaligned
\tag 3.3.18$$
It follows from the condition (3.3.2) that (3.3.18) is larger than
$1-\epsilon$ for almost all $l$.  Thus the image of the set $\{\nu_
l\}$ under the
map $x\to proj_{\Bbb Cg\cdot v}(x_{(j)})$ is tight, for each $g\in
K_{(j)}$.  This implies
our claim in this case, since a finite number of the $g\cdot v$ will
span $K_{(j)}\cdot v$. This completes the proof in the case $k=0$.

Now suppose that $k>0$, and fix $g\in K_{(j)}$ as before.  By
induction we know that the image of $\{\mu_{l,k}\}$ under the map
$$x\to proj_{\Bbb Cf_j+\Bbb Cv}(x)\tag 3.3.19$$
is tight.

\proclaim{(3.3.20)Lemma} The image of $\{\vert\frac {g^{-1}_{*}\sigma_{
\Lambda}}{\sigma_{\Lambda}}\vert^{2k}d\mu_{l,k}\}$ under the
map (3.3.19) is tight.
\endproclaim

\flushpar Proof of (3.3.20).  In (ii) of Remark (3.2.2), we
calculated that
$$\vert\frac {g^{-1}_{*}\sigma_{\Lambda}}{\sigma_{\Lambda}}\vert^{
2k}(x)=\vert a+bz\vert^{2k},\tag 3.3.21$$
where as before $zf_j$ is the $\frak g_{-\alpha_j}$ component of $
x$, and $g$
corresponds to $\left(\matrix a&b\\c&d\endmatrix \right)\in SU(2,
\Bbb C)$. Thus
$$(proj_{\Bbb Cf_j+\Bbb Cv})_{*}(\vert\frac {g_{*}\sigma_{\Lambda}}{
\sigma_{\Lambda}}\vert^{2k}\mu_{l,k})=\vert a+bz\vert^{2k}(proj_{
\Bbb Cf_j+\Bbb Cv})_{*}\mu_{l,k}.\tag 3.3.22$$

We thus have the following data:
$$dm_l(z,w)=(proj_{\Bbb Cf_j+\Bbb Cv})_{*}\mu_{l,k}$$
is a tight sequence of probability measures on $\Bbb Cf_j+\Bbb Cv$;
$$\delta (z)=\vert a+bz\vert^{2k}$$
is a fixed nonnegative function on $\Bbb Cf_j$; and
$$\delta (z)\int_wdm_l(z,w)$$
is a tight sequence of probability measures on $\Bbb Cf_j$.  The last
bit of data implies that given $\epsilon >0$, there is an $M>0$ such that
$$\int_{\vert z\vert >M}\delta (z)\int_wdm_l(z,w)<\epsilon /2,\quad
\forall l.\tag 3.3.23$$
Since $\delta$ is bounded on $\{\vert z\vert\le M\}$, and $dm_l$ is tight,
there is an
$N>0$ such that
$$\int_{\vert z\vert <M}\delta (z)\int_{\vert w\vert >N}dm_l(z,w)
<\epsilon /2,\quad\forall l.\tag 3.3.24$$
It follows that $\delta dm_l$ will have measure exceeding $1-\epsilon$ on
$\{\vert z\vert <M,\vert w\vert <N\}$, for all $l$. This implies (3.3.20). //

\smallskip

We can now complete the proof of (3.3.8) by arguing as in the
case $k=0$.  Given $\epsilon >0$, the lemma implies that there exists a
compact set $C_{\epsilon}\subset \Bbb Cv$ such that for
$$E=\{x:proj_{\Bbb Cv}(x_{(j)})\in C_{\epsilon}\}\tag 3.3.25$$
$$\int_E\vert\frac {g_{*}^{-1}\sigma_{\Lambda}}{\sigma_{\Lambda}}
\vert^{2k}d\mu_{l,k}>1-\epsilon .\tag 3.3.26$$
Then
$$\mu_{l,k}(\{proj_{\Bbb Cg\cdot v}(x_{(j)})\in g\cdot C_{\epsilon}
\})=g_{*}\mu_{l,k}(E),\tag 3.3.27$$
and this is approximately equal to the left hand side of (3.3.26)
for $l$ large, by condition (3.3.2), as spelled out by (3.3.7).  Thus
(3.3.27) will exceed $1-\epsilon$ for large $l$.  This completes the proof
for $k>0$.  //

\smallskip

\proclaim{(3.3.24)Corollary} Suppose that the Cartan matrix is of
finite rank and of affine type. Then there exists a $K(A)$-invariant
probability measure on $\Cal F$.
\endproclaim

This follows from the fact that the conditioned Wiener
measures satisfy the condition (3.3.2), as the temperature tends
to infinity; see \S 4.1 of Part III.

\bigskip
\centerline{Part II. Infinite Classical Groups.}

\bigskip

\centerline{\S 0. Introduction and Notation.}

\bigskip

\noindent\S 0.1. The Examples, and the Kac-Moody Completions.

The irreducible affine infinite rank Kac-Moody algebras are
precisely the classical inductive limits $sl(\infty ,\Bbb C)$, $s
l(2\infty ,\Bbb C)$,
$o(2\infty ,\Bbb C)$, $o(2\infty +1,\Bbb C)$, and $sp(\infty ,\Bbb C
)$ (see \S 7.7 of [Kac]).  We will
write $\frak g=\frak g(\infty )=lim\,\,\frak g(n)$ for such an inductive limit.
 We have
$$G(A)=G(\infty )=\lim_{n\to\infty}G(n),$$
where $G(n)$ is the unique connected, simply connected complex
Lie group corresponding to the Lie algebra $\frak g(n)$.

The Kac-Moody strong operator closure of the unitary real forms are
well-known. In the case of $sl(\infty ,\Bbb C)$, the defining representation is
of
highest weight type, so that the Kac-Moody strong operator closure is
just the ordinary strong operator closure of $SU(\infty )$, namely $
U(H)$, the
unitary group of a separable infinite dimensional complex Hilbert
space. This group is a central extension of the projective unitary
group by its center,
$$0\to \Bbb T\to U(H)\to \Bbb PU(H)\to 0.$$

In the other cases the defining representation is not of highest
weight.  In these other cases the Kac-Moody strong operator
closures, as groups, are equal to the universal central extensions of
the identity components of the Hilbert-Schmidt restricted unitary,
orthogonal (for both type B and D), and unitary symplectic groups,
respectively.  We will denote these extensions by
$$0\to \Bbb T\to\tilde {U}_{(2)}\to U_{(2)}\to 0,$$
$$0\to \Bbb T\to\tilde {O}_{(2)}\to O_{(2)}\to 0,\tag 0.1.1$$
and
$$0\to \Bbb T\to\tilde {S}p_{(2)}\to Sp_{(2)}\to 0,$$
respectively. These extensions are thoroughly studied as Banach Lie
groups in [PS], where the extensions are realized using the canonical
determinant line bundle over Toeplitz operators. We only need to add
the following

\bigskip

\noindent(0.1.2)Remarks.  (a) The Kac-Moody strong operator
topology on the above groups is weaker than the Banach Lie
topologies studied in [PS].  The former topology involves the
strong operator topology on Toeplitz operators rather than the
norm topology.  So far as representation and measure theory
are concerned, the difference between the topologies is
inconsequential (see section 5 of [Pi2]).

(b) The group $U(H)$ is the isometry group of a finite rank,
infinite dimensional Riemannian symmetric space ($Gr(n,H),$ for
any $n$), while the restricted groups above are isometry groups
of certain infinite rank Riemannian symmetric spaces
(Grassmannians of the appropriate type).  The collection of
automorphism groups of separable irreducible Riemannian
symmetric spaces is, to my knowledge, the only class of
infinite dimensional groups for which there is a substantial and
(apparently) tractable theory of general unitary
representations.  This theory, due largely to Ol'shanskii,
generalizes aspects of Harish-Chandra's work on representations
of finite dimensional semisimple groups in a remarkable way
(see [Ol] for a general account).

\bigskip

\flushpar\S 0.2. Contents.

In the next section we will give a simple conceptual proof of
existence and uniqueness of measures on the formal flag space
for the infinite classical groups. Unfortunately we have very
little to say about the structure of these measures, or
harmonic analysis.

The remaining sections take up the analysis of measures on
the formal completion $G(\infty )_{formal}$, on a case by case basis.
This is necessitated by the curious fact that the results are
not uniform.  There is a single anomalous example, namely
$SL(\infty )$.  The odd features are that (1) there exists an invariant
measure on $SL(\infty )_{formal}$, despite the fact that $G$ has a
nontrivial center, and (2) the center does not act transitively
on the measures.  In all the other cases there is a simple
existence and uniqueness statement for measures on
$G(\infty )_{formal}/C_0$, and there does not exist an invariant measure
on the formal completion itself.

In the case of $SU(\infty )$, there is a simple discrete version of the
Plancherel theorem, in complete analogy with the situation for
compact groups. Ol'shanskii has conjectured that the other
unitary groups considered here are also type I. This suggests
the existence of interesting Plancherel formulae for these
other groups as well. However this has not been worked out
completely even in the case of symmetric quotient spaces (but
see the recent announcement [KOV] for an interesting
analogue).

\bigskip

\centerline{\S 1. Measures on the Formal Flag Space.}

\bigskip

\flushpar\S 1.1. Existence and Uniqueness.

We will use the same notation as in (3.2) of Part I.

\proclaim{(1.1.1) Proposition}  For each $s\ge 0$, there exists a unique
measure $\mu^{\vert \Cal L\vert^{-2s}}$ on $\Cal F$ having values in $
\vert \Cal L\vert^{-2s}$ such that (i) $\mu^{\vert \Cal L\vert^{-
2s}}$ is
$K(\infty )-$invariant, and (ii) $d\mu_s=\vert\sigma\vert^{2s}d\mu^{
\vert \Cal L\vert^{-2s}}$ is a probability measure.
\endproclaim

\noindent Proof.  Let $\Delta^{-}(n)$ and $\Pi (n)$ denote the sets of negative
and
simple roots for $\frak g(n)$, respectively.  Also let $\Pi (n)^c
=\Pi\setminus\Pi (n),$ and
$\Delta^{-}(n)^c=\Delta^{-}\setminus\Delta^{-}(n)$.  Define
$$\Cal N^{-}(\Delta^{-}(n)^c)=exp(\,\,\,(\sum_{\alpha\in\Delta^{-}
(n)^c}\,\,\,\frak g_{\alpha})_{formal}),$$
and
$$\Cal R^{-}=(\prod_{\alpha_j\in\Pi (n)^c}H_{(j)})\,\,\propto\,\,
\Cal N^{-}(\Delta^{-}(n)^c).$$
The group $\Cal R^{-}$ is the radical of the lower parabolic subgroup of $
\Cal G^{-}$
corresponding to $\Pi (n)$. In particular it is normalized by $G(
n)$.

Consider the following commutative diagram, analogous to (3.1.2)
of Part I:
$$\matrix \Cal L_{\Lambda}\to \Cal G/\Cal B^{+}&\mapsto&\Cal R^{-}
\backslash \Cal L_{\Lambda}\to \Cal R^{-}\backslash \Cal G/\Cal B^{
+}\\\uparrow&&\uparrow\\\Cal L_{\Lambda}\to\bigcup_{w\in W(n)}w\cdot
\Sigma_1&\mapsto&\Cal L_{\Lambda}\to G(n)/B^{+}(n)\\\uparrow&&\uparrow\\
\frak n_{formal}^{-}\times \Bbb C&\mapsto&\frak n^{-}(n)\times \Bbb C\endmatrix
\!\tag 1.1.2$$
In the last line we have used several identifications.  First we used
the section $\sigma_{\Lambda}$, and the isomorphism $\Cal N^{-}\cong
\Sigma_1$, to identify the
restriction of $\Cal L_{\Lambda}$ to $\Sigma_1$ with $\Cal N^{-}\times
\Bbb C$.  Secondly, we used the
isomorphism $exp:\frak n^{-}_{formal}\to \Cal N^{-}$.  The projection $
\frak n^{-}_{formal}\to \frak n^{-}(n)$ is
simply the linear projection along the root spaces.  The top arrow is
$G(n)$-equivariant, because $\Cal R^{-}$ is normalized by $G(n)$.  All of the
up
arrows are injective.

The uniqueness of $\mu^{\vert \Cal L\vert^{-2s}}$ follows from (1.1.2), since $
K(n)$-invariance
determines its image in $\Cal L_{\Lambda}\to G(n)/B^{+}(n)$.  The existence of
$
\mu^{\vert \Cal L\vert^{-2s}}$ also
follows from (1.1.2), together with the Kolmogorov extension
theorem.//

\flushpar Remark.  The above argument clearly applies to
generalized flag spaces as well, in particular to the
Grassmannian homogeneous spaces.  It is the same argument as
in [Pi1], except that the conceptual simplicity of the argument
was completely obscured in [Pi1] by the introduction of
coordinates. See [Pi1] for several explicit representations of the
Grassmannian measure.

\bigskip

\noindent\S 1.2. Comparison with Product Measures.

Measures on infinite dimensional spaces tend to be small
perturbations of product measures.  Here I would like to suggest how
this might be true for $d\mu_0$, the invariant probability on $\Cal F$.  I do
not
know how to generalize this to the case of measures with
values in a line bundle.

For each simple root $\alpha$, let $P^{+}(\alpha )$ and $\Cal P^{
+}(\alpha )$ denote the maximal
parabolic subgroups of $G$ and $\Cal G^{+}$, respectively, corresponding to $
\{\alpha \}$.
Let $\Cal X(\alpha )$ denote the formal completion of the minimal flag space
$G/P^{+}(\alpha )$,
$$\Cal X(\alpha )=\Cal G/\Cal P^{+}(\alpha )=\Cal G^{-}/P^{+}(\alpha
).$$
The measure $d\mu_0$ projects to the unique invariant probability on
$\Cal X(\alpha )$.

For the examples at hand, the space $G/P^{+}(\alpha )$ is one of the
Grassmannians $U(\infty )/(U(n)\times U(\infty ))$, $U(2\infty )/
(U(\infty )\times U(\infty ))$,
$SO(2\infty )/U(\infty )$, and $Sp(\infty )/U(\infty )$.  Each of these spaces
is a
(incomplete) Riemannian symmetric space.  For spaces of this
type it is known (abstractly) how to express the invariant
measure on the formal completion as a product measure, using
polar coordinates, i.e.  the ergodic decomposition of the
measure relative to the stability subgroup (see [Pi4]).

There is a natural map
$$\Cal F\to\prod \Cal X(\alpha_j).$$
Is the image of $d\mu_0$ related in an elementary way to the
product of the invariant measures?

\bigskip

\centerline{\S 2. The Case $\frak g=sl(\infty ,\Bbb C)$.}

\bigskip

\noindent\S 2.1. Borel embedding of $G_{formal}$ into $\infty\times
\infty$ matrices.

Let $\epsilon_1,\epsilon_2...$ denote an ordered orthonormal basis for a
countably
infinite dimensional complex Hilbert space H.  We will identify $
\Bbb C^n$
with the span of $\{\epsilon_j\vert j\le n\}.$ Then $\frak g=sl(\infty
,\Bbb C)=lim$ $sl(n,\Bbb C)$, and $\frak h$
consists of the diagonal matrices, $\frak n^{+}$ consists of the upper
triangular
matrices, and $\frak n^{-}$ consists of the lower triangular matrices in $
\frak g$.

Let $\Bbb M=_{}\Bbb M_{\infty ,\infty}$ denote the space of complex matrices
indexed by
$\Bbb N\times \Bbb N$.  Given $g\in \Bbb M$ and $n\in \Bbb N$, we will write $
g=\left(\matrix A&B\\C&D\endmatrix \right)$, where
$A=A^{(n)}\in \Bbb M_{n,n}$, the space of $n\times n$ matrices.  We then have
matrix
representations
$$\Cal H\tilde {=}\{diag(z_j)\in \Bbb M\vert\lambda_j\ne 0,\forall
j\},\tag 2.1.1$$
where $\prod\sigma_j^{h_j}\in \Cal H$ corresponds to $diag(z_j)$, $
z_1=\sigma_1,$ $z_{n+1}=\sigma_{n+1}/\sigma_n,$
$n\ge 1$,
$$G=\{g\in \Bbb M\vert g_{ij}=0,\forall i\ne j\ni i+j\gg 0;detA_n
\ne 0,\forall n\gg 0\},$$
$$\Cal N^{+}\tilde {=}\{g\in \Bbb M\vert g_{ij}=0,\forall i,j\ni
i>j;g_{ii}=1,\forall i\},\tag 2.1.2$$
$$\Cal G^{+}\tilde {=}\{g\in \Bbb M\vert g_{ij}=0,\forall i>j\ni
i+j\gg 0;detA_n\ne 0,\forall n\gg 0\}.$$
Similarly $\Cal N^{-}$ and $\Cal G^{-}$ are represented by essentially lower
triangular
matrices.  Note that after $SL(\infty ,C)$ is augmented, the condition
$det=1$ never reappears. Note also that $\Cal G^{-}$ acts from the left on $
\Bbb M$, and
$\Cal G^{+}$acts from the right, by matrix multiplication.

\bigskip

\proclaim{(2.1.3) Proposition}  (a) The map $\imath :\Cal G\to \Bbb M
:[g,h]\to g\cdot h$, given by
matrix multiplication, is regular and $\Cal G^{-}\times \Cal
G^{+}$-equivariant.

(b) The image of $\imath$ is given by
$$\imath (\Cal G)=\{g\in \Bbb M\vert detA_n(g)\ne 0,\forall n\ni
n\gg 0\}=liminf\{g\vert detA_n\ne 0\}.$$

(c) $\imath$ is a Borel isomorphism of $\Cal G$ onto $\imath (\Cal G
)$; $\imath$ is not an open map onto
$\imath (\Cal G)$, hence $\imath (\Cal G)\in G_{\delta\sigma}(\Bbb M
)\setminus G_{\delta}(\Bbb M)$; $\imath$ is a homeomorphism of the top
stratum $\Cal N^{-}\times \Cal H\times \Cal N^{+}$ onto its image.

(d) $\imath^{*}:\Cal O_{reg}(\Bbb M)\to \Cal O_{reg}(\Cal G)$ is an isomorphism
of $
\Cal G^{-}\times \Cal G^{+}$-modules.
\endproclaim

\bigskip

\noindent Proof.  Part (a) is obvious.  Part (b) and the injectivity of
$\imath$ in (c) follow from the basic facts about Gaussian elimination.

The first statement of (c) follows from the general fact that an
injective map of Polish spaces is a Borel isomorphism of its domain
onto its image. If $r_n\in\imath (\Cal G)$ is the permutation matrix which
switches
$\epsilon_n$ and $\epsilon_{n+1}$, and fixes all other $\epsilon_
j$, then $r_n$ is not in the image of the
top stratum, but nonetheless $r_n$ converges to the identity matrix in
$\Bbb M$. This shows that $\imath$ is not open. On the top stratum $
\imath$ is inverted
by the Gaussian elimination process, hence it is open. This proves c).

Part d) follows from taking the direct limit of the isomorphism

$$(\Cal O_{reg}(\Bbb M_{n,n})/<det=1>)\cong \Cal O_{reg}(SL(n,\Bbb C
)).//$$

\bigskip\

\noindent\S 2.2. Classification of $SU(\infty ,\Bbb C)$-biinvariant measures on
$
\Cal G$.

The classification of invariant probability measures on the
linear space $\Bbb M$ is known.  The basic examples are the
following.  First, let $\nu_{\beta}$ denote the Gaussian measure on $
\Bbb M$
corresponding to inverse temperature $\beta$, where $0\le\beta <\infty$.
Formally,
$$d\nu_{\beta}(g)=E^{-1}exp(-\beta tr(g^{*}g)/2)dLeb(g);$$
rigorously, if $\beta >0$, then the Fourier transform of $\nu_{\beta}$ is given
by
$$\hat{\nu}_{\beta}(x)=\int exp(-itr(x^{*}g))d\nu_{\beta}(g)=exp(
-tr(x^{*}x)/(2\beta )),\tag 2.2.1$$
where $x\in \Bbb M^{*}$, $\nu_0=\delta_0.$ Clearly $\nu_{\beta}$ is $
U(\infty )$-biinvariant.

Secondly, suppose that $d_1\ge d_2\ge ...\ge 0$.  If $\sum d_j^2<
\infty ,$ then the map
$$\phi :\Bbb M\times \Bbb M\mapsto \Bbb M:X,Y\mapsto X\cdot diag(
d_j)\cdot Y,$$
is defined for a.e. $(X,Y)\in (\Bbb M\times \Bbb M,\nu_1\times\nu_
1)$. This is so because for each
$(i,j),$ $(X_{ik}d_kY_{kj})_{0<k<\infty}$ is a sequence of independent random
variables,
each with expectation zero, such that the corresponding variances are
summable, for
$$\sum d_k^2\cdot var(X_{ik}Y_{kj})=constant\sum d_k^2<\infty ;$$
hence $\sum X_{ik}d_kY_{kj}$ is conditionally summable. Let $\nu_{
\vec {d}}$ denote the image
of $_{}\nu_1\times\nu_1$ with respect to the map $\phi$ above. The measure $
\nu_{\vec {d}}$ is also
clearly $U(\infty )$-biinvariant. A straightforward calculation shows that
$$\hat{\nu}_{\vec {d}}(x)=\prod_{i=1}^n\prod_{j=1}^{\infty}(1+d_j^
2u_i^2)^{-1},\tag 2.2.2$$
where $x$ is an $n\times n$ matrix, and $spectrum(\vert x\vert )=
\{u_1,..,u_n\}$.

The following is proven in [Pi3].

\proclaim{(2.2.3)Proposition} Each ergodic $SU(\infty )\times SU(
\infty )$-invariant
probability on $\Bbb M$ is of the form $\nu_{\beta}*\nu_{\vec {d}}$, where $
\beta$ and $\vec {d}$ are uniquely
determined by the conditions $\beta\ge 0$ and $\vec {d}=(d_j)$, where $
d_1\ge d_2\ge ..\ge 0$,
and $\sum d_j^2<\infty .$ \endproclaim

It is a simple matter to determine which of these measures is
supported on $\Cal G$. If $d_{n+1}=0$, then $\nu_{\vec {d}}$ is supported on
the complement
of $\Cal G$. For the image of $\nu$ $_{\vec {d}}$ in $\Bbb M_{N,N}$ equals the
image of the Gaussian
relative to the map
$$\phi_N:\Bbb M_{N,N}\times \Bbb M_{N,N}\mapsto \Bbb M_{N,N}:X,Y\mapsto
X\cdot diag(d_1,..,d_n,0,..,0)\cdot Y,$$
where $n<N$.  Thus $detA_N=0$ with $\nu_{\vec {d}}$ measure one, whenever $
n<N$.
In all other cases, for all $N$, the image of $\nu_{\beta}*\nu_{\vec {
d}}$ in $\Bbb M_{N,N}$ is
represented by a real analytic density relative to Lebesgue measure.
This can be checked directly from the explicit form for the Fourier
transform. This proves the following

\proclaim{(2.2.4)Corollary}The ergodic $SU(\infty )\times SU(\infty
)$-invariant
probability measures on $\Cal G$ are of the form $\nu_{\beta}*\nu_{
\vec {d}}$, where either $\beta >0$,
or if $\beta =0$, then $d_j>0$, $\forall j.$ \endproclaim

This completely determines the space of unitarily biinvariant
probability measures on $\Cal G$, but the approach in [Pi3] is not in
the spirit of this paper. The better approach is suggested by
Remark (3) of (0.1.4) in the introduction. If $\Cal G$ was a finite
dimensional complex group, then there would be an Iwasawa (or
$KAN$) decomposition, and the Abel transform would define
isomorphisms
$$\Cal M(K\backslash \Cal G/K)\to \Cal M(A)^W\quad and\quad Prob(
K\backslash \Cal G/K)\to Prob(A)^W,\tag 2.2.5$$
where $\Cal H=A\cdot T$ is the polar decomposition and $W$ is the Weyl
group. But in our infinite dimensional context, there does not
exist an Iwasawa decomposition. A reasonable substitute is the
diagonal distribution
$$a_{*}:\Cal M(\Cal G)^{SU(\infty )\times SU(\infty )}\to \Cal M(
A):\nu\to a_{*}\nu ;\tag 2.2.6$$
here $a$ denotes the $A$-component of $diag\in \Cal H$. It will follow
from what we say in \S 4.4 below that this map is injective.
However at this point we do not have sufficient information
about the image of $a_{*}$ to give another proof of (2.2.4).

\bigskip

\noindent\S 2.3.  The analogue of Haar measure.

In this section we will argue that the Gaussians $\nu_{\beta}$ are the
natural analogues of Haar measure for $SU(\infty )$, viewed as the
unitary real form of a Kac-Moody group.  The following result
indicates that practically any reasonable condition in addition
to biinvariance singles these measures out. In the next
subsection we will see that one can also prescribe the diagonal
distribution via a regularization of the {\bf c}-function, as suggested
in the introduction.

\proclaim{(2.3.1)Proposition} (a) If $\nu_{}$ is a $SU(\infty )$-biinvariant
probability measure on $\Cal G$ (or more generally $\Bbb M$), then $
\nu$ is
ergodic for the left action of $SU(\infty )$ if and only if $\nu$ is
Gaussian.

(b) Every probability measure on $\Bbb M$ which is invariant and
ergodic for the left action of $SU(\infty )$ is an equivariant image of
$\nu_{\beta}$ $(\beta >0)$.

(c) The invariant distribution on $(n/\beta )^{1/2}SU(n)$ converges
weakly to $\nu_{\beta}$ as $n\to\infty$.

(d) As a representation of $U(\infty )\times U(\infty )$,
$$L^2(\Cal G,d\nu_{\beta})\tilde {=}\sum\oplus\rho^{*}\otimes\rho$$
$$(L^2\Cal O(\Cal G,\!\nu_{\beta})\cong\sum\oplus\rho^{*}\otimes\rho
,\quad resp)$$
where the sum is over all irreducible separable unitary
representations of $U(H)$ (holomorphic representations, resp).
\endproclaim

Parts (a) and (b) follow from (2.2.3), (c) is a straightforward
calculation (see [Pi4]), and (d) is essentially due to I.E. Segal
(see [Pi3]).

\bigskip

\noindent\S 2.4.  Calculation of diagonal distributions.

Suppose that $g$ is in the top stratum of $\Cal G$. We have
$$g=l\cdot diag(z_j)\cdot u,\tag 2.4.1$$
where $l\in \Cal N^{-},u\in \Cal N^{+},$ $\prod\sigma_j(g)^{h_j}\in
\Cal H$ corresponds to $diag(z_j)$ as in
(2.1.1), and
$$\sigma_j(g)=detA^{(j)}(g).\tag 2.4.2$$
The diagonal distribution of $\nu_{\beta}$, $diag_{*}\nu_{\beta}$, is invariant
with
respect to the unitary part of $\Cal H$.  Let $a=diag(\vert z_j\vert
)$ denote the
$A$-component, regarded as a function on the top stratum,
and suppose that
$$\lambda\in exp(\frak h_{\Bbb R}\hat {)}=\sum \Bbb R\epsilon^{*}_
j.\tag 2.4.3$$
We now use part (c) of (2.3.1) and the formula (0.2.9) of
Harish-Chandra to compute the Fourier transform for the
diagonal distribution of $\nu_{\beta}$ (we will take $\beta =1$ for
simplicity).
We have
$$\aligned
(a_{*}\nu_1\hat {)}(\lambda )&=\int_{\Cal G}e^{-i\lambda loga(g)}
d\nu_1(g)=\lim_{n\to\infty}\int_{SU(n)}e^{-i\lambda loga(n^{-1/2}
g)}dm(g)\\&=\lim_{n\to\infty}exp(\frac i2logn\sum\lambda_j)\prod_{
1\le j<k\le n}(1+\frac i2\frac {\lambda_j-\lambda_k}{j-k})^{-1}\\&
=(\prod_{1\le j<k}(1+\frac i2\frac {\lambda_j-\lambda_k}{j-k})e^{
-\frac i2\frac {\lambda_j}{j-k}})^{-1}\lim_nexp(\frac i2\sum\lambda_
j(logn-\sum_{j+1}^n\frac 1{k-j}))\\&=exp(\frac i2\gamma\sum\lambda_
j)(\prod_{1\le j<k}(1+\frac i2\frac {\lambda_j-\lambda_k}{j-k})ex
p(-\frac i2\frac {\lambda_j}{j-k}))^{-1}\endaligned
\tag 2.4.4$$
where $\gamma$ denotes Euler's constant.  This is clearly a
regularization of the formal {\bf c}-function for $sl(\infty ,\Bbb C
)$, although
the Lie-theoretic significance is not so transparent.

It is interesting to check this formula by calculating the
individual distributions of the $\sigma_j$ in a different way.  It
suffices to calculate $\vert\sigma_n\vert^2_{*}\nu_1\in Prob(\Bbb R^{
+})$.  For $s\in \Bbb R=(\Bbb R^{+}\hat {)}$, the
Fourier transform is given by
$$(\vert\sigma_n\vert^2_{*}\nu_1\hat {)}(s)=\frac 1E\int_{GL(n,\Bbb C
)}det(g^{*}g)^{-is}e^{-\frac 12tr(g^{*}g)}dLeb(g)\tag 2.4.5$$
By using standard integral formulae for the Cartan
decomposition of $GL(n,\Bbb C)$, and a Selberg integral calculation, one
finds
$$\aligned
(\vert\sigma_n\vert^2_{*}\nu_G\hat {)}(s)=&\frac 1E\int_0^{\infty}
..\int_0^{\infty}(u_1..u_n)^{-is}\prod_{i<j}(u_i-u_j)^2e^{-\frac
12\sum u_j}\prod du_j\\=&\prod_1^n\frac {\Gamma (j-is)}{\Gamma (j
)}.\endaligned
\tag 2.4.6$$
to which (2.4.4) reduces upon inserting $\lambda_j=2s$ for $1\le
j\le n$, and
$\lambda_j=0$ otherwise.

\bigskip
\centerline{\S 3. The Case $\frak g=sl(2\infty ,\Bbb C).$}

\bigskip

\flushpar\S 3.1. On Realizing $\Cal G$ $^{\pm}$ and $\Cal G$.

Let $\epsilon_j,$ $j\in \Bbb Z,$ denote an ordered orthonormal basis for a
complex
Hilbert space H.  We will identify $\Bbb C^{2n}$ with the span of
$\{\epsilon_j\vert -n\le j<n\}$. Then $\frak g=sl(2\infty ,\Bbb C
)=lim$ $sl(2n,\Bbb C)$, $\frak h$ consists of the
diagonal matrices, $\frak n^{+}$ consists of the upper triangular matrices, and
$\frak n^{-}$ consists of the lower triangular matrices in $\frak g$, and
$h_i=\epsilon_i\otimes\epsilon_i^{*}-\epsilon_{i+1}\otimes\epsilon_{
i+1}^{*}$, $i\in \Bbb Z$.

Let $\Bbb M=\Bbb M_{2\infty ,2\infty}$ denote the space of matrices indexed by
$
\Bbb Z\times \Bbb Z$.
Using the basis $\{\epsilon_j\}$, we can faithfully represent $\frak g$ and
$G(\infty )=SL(2\infty ,\Bbb C)$ in $\Bbb M$.  We can extend this
representation to the
groups $G$ and $\Cal G^{\pm}$, but the extension is not faithful.  The basic
problem is that the extended map
$$\Cal H=\prod\,\,H_{(j)}\mapsto \Bbb M:\,\,\prod\,\,a_j^{h_j}\mapsto
diag(a_j/a_{j-1}),\tag 3.1.1$$
has kernel $\Bbb C^{*}$, where $\Bbb C^{*}$ is diagonally embedded in $
\Cal H$.  The upshot is
that we can only realize $G/\Bbb C^{*}$ and $\Cal G^{\pm}/\Bbb C^{
*}$ in $\Bbb M$, namely
$$G/\Bbb C^{*}\cong \{g\in \Bbb M\vert g_{ij}=0,\forall i\ne j,\,\,
\vert i\vert +\vert j\vert\gg 0;\,\,det((g_{ij})_{-n\le i,j\le n}
)\ne 0,\,\,n\gg 0\},$$
$$\Cal G^{+}/\Bbb C^{*}\cong \{g\in \Bbb M\vert g_{ij}=0,\forall
i>j,\,\,\vert i\vert +\vert j\vert\gg 0;\,\,det((g_{ij})_{-n\le i
,j\le n})\ne 0,\,n\gg 0\},$$
i.e. $\Cal G^{+}/\Bbb C^{*}$ is realized by invertible, asymptotically upper
triangular
matrices; similarly $\Cal G^{-}/\Bbb C^{*}$ is realized by invertible
asymptotically lower
triangular matrices.

\bigskip

\flushpar Remarks (3.1.2).  (a) $_{}\Cal N^{\pm}$ are faithfully represented by
unipotent upper and lower triangular matrices in $\Bbb M$, respectively.
Thus there is a concrete way to think about the top stratum of $\Cal G$, in
terms of the canonical factorization $\Cal N^{-}\cdot \Cal H\cdot
\Cal N^{+}$.

(b) The extension $0\to \Bbb C^{*}\to \Cal G^{+}\to \Cal G^{+}/\Bbb C^{
*}\to 0$ can be realized directly using the
same device as in [PS] (see section 6.6). One simply replaces the
group $GL(H_{+})$ in [PS] by the group of invertible matrices indexed by $
\Bbb N^c$
which are asymptotically upper triangular.

(c) If $a\in \Bbb C^{*}$ and $g\in SL(2\infty ,\Bbb C)\hookrightarrow
\Cal G$, then
$$\sigma_n(\prod a^{h_j}\cdot g)=a\cdot det(A^{(n)}(g)),$$
where $A^{(n)}(g)=(g_{ij})_{i,j\le n}$.

(d) The sum of the fundamental highest weight representations is
the wedge representation, $\Lambda^{\infty}(\oplus \Bbb C\epsilon_
j)$ (see [Kac] or [PS]), so
that there is an injection
$$\Cal G\to \Bbb M_{\Lambda^{\infty}}.$$
I do not know whether there is a useful description of the image.

\bigskip

\flushpar\S 3.2.  Measures on $\Cal G/\Bbb C^{*}$; Statement of the Main
Result.

Let $\Cal G\times_{\Bbb C^{*}}\Bbb C$ denote the complex line bundle associated
to the
$\Bbb C^{*}$-bundle $\Cal G\to \Cal G/\Bbb C^{*}$.  We will write $
\vert \Cal G\vert$ in place of the clumsy
expression $\vert \Cal G\times_{\Bbb C^{*}}\Bbb C\vert$, and similarly we will
denote the dual line
bundle by $\Cal G^{*}$.  Note that
$$\Cal O_{reg}((\Cal G^{*})^{\otimes k})\cong \{f\in \Cal O_{reg}
(\Cal G)\vert f(g\cdot\lambda )=\lambda^kf(g),\lambda\in \Bbb C^{
*}\},$$
i.e.  the space of matrix coefficients corresponding to
representations of level $k$ for the action of the center $\Bbb C^{
*}.$

\bigskip

\proclaim{(3.2.1) Proposition} For each $s>-1$, there exists a
measure $\mu =\mu^{\vert \Cal G\vert^{2s}}$ on $\Cal G/\Bbb C^{*}$ having
values in $
\vert \Cal G\vert^{2s}$ such
that (i) $\mu$ is $SU(2\infty ,\Bbb C)$-biinvariant, and (ii) the measure $
\vert\sigma\vert^{2s}d\mu$
is finite for every regular section $\sigma$ of $\Cal G^{*}\to \Cal G
/\Bbb C^{*}$.
\endproclaim

We will essentially give three proofs of this result (at least
for $s=0$).  The first proof is straightforward, in that we will
write down explicit formulas for the projections of $d\mu$ onto
double coset spaces, and we will directly check their
consistency; we will do this in the next subsection.  The
second proof will be given in \S 3.6.  The third proof will be
given in the next section; it applies uniformly to all the
classical groups.

\bigskip

\flushpar\S 3.3. Explicit Expressions for $\vert\sigma_r\vert^{2s}
d\mu$.

Given $n$, let $\Cal R^{+}$ and $\Cal R^{-}$ denote the radicals of the upper
and lower
parabolic subgroups corresponding to $\Pi (n)$, the set of simple roots for
$G(n)=SL(2n,\Bbb C)$; i.e.
$$\Cal R^{+}=\,(\,\prod_{\alpha_j\in\Pi (n)^c}H_{(j)})\,\,\propto\,\,
\Cal N^{+}(\Delta^{+}(n)^c),$$
and $\Cal R^{-}$ is defined as in (1.1.2). Consider the commutative diagram
$$\matrix \Cal G/\Bbb C^{*}&\mapsto&\Cal R^{-}\backslash (\Cal G/
\Bbb C^{*})/\Cal R^{+}\\\uparrow&&\uparrow\\\bigcup_{w\in W(n)}w\cdot
(\Cal N^{-}\cdot \Cal H/\Bbb C^{*}\cdot \Cal N^{+})&\mapsto&GL(2n
,\Bbb C)\\\uparrow&&\uparrow\\\frak n^{-}_{formal}\times \Cal H/\Bbb C^{
*}\times \frak n^{+}_{formal}&\mapsto&\frak n^{-}(n)\times H^{GL}
(n)\times \frak n^{+}(n)\\\endmatrix \tag 3.3.1$$
The top arrow is $SL(2n,\Bbb C)\times SL(2n,\Bbb C)$-equivariant, the bottom
arrow is
given by linear projection along root spaces, and all the up arrows
are injective. Note that the double coset space is essentially
$GL(2n)$, rather than $SL(2n)$, because when we project $\Cal H$ into $
\Bbb M$
by (3.1.1), the determinant condition is not preserved.

The following Proposition asserts the existence of the measure
$\vert\sigma_r\vert^{2s}d\mu$.  In what follows it will be implicitly
understood that
$N>n>\vert r\vert .$

\bigskip

\proclaim{(3.3.2) Proposition} (a) For $s>-1$ the measure on $GL(
2n,\Bbb C)$
defined by
$$det\vert A^{(r)}(g)\vert^{2s}det(1+g^{*}g)^{-4n-s}dLeb(g),$$
is finite.  After normalization these measures, denoted $d\mu_s^{
(n)}$, define a
probability measure, $d\mu_s=\vert\sigma_r\vert^{2s}d\mu$, on $\Cal G
/\Bbb C^{*}$ via the system of
projections (3.3.1).

(b) The measure $d\mu_s$ is S$U(2\infty ,\Bbb C)\times SU(2\infty
,\Bbb C)$-quasi-invariant, and the
Radon-Nikodym derivative is given by
$$\frac {d\mu_s(g_L\cdot g\cdot g_R)}{d\mu_s(g)}=\vert\frac {\sigma_
r(g_L\cdot g\cdot g_R)}{\sigma_r(g)}\vert^{2s}.$$
\endproclaim

\noindent Remarks (3.3.3).  The measure appearing in (3.3.2) (with
$s=0$) can be realized as a quotient of a Gaussian.  This is
suggestive from our point of view, because ultimately we
would like to think of (3.2.1) as an instance of a projective
analogue of the central limit theorem.  Here we record several
facts that we will use about the measures appearing in (3.3.2);
see [Pi1] for details.

 (a) The measure
$$det(1+\vert Z\vert^2)^{-2M}dLeb(Z),$$
properly normalized, is the unique $U(2M,\Bbb C)$-invariant
probability on $Gr(M,\Bbb C^{2M})$, expressed in terms of the graph
coordinate, $\Cal L(\Bbb C^M)\to Gr(M,\Bbb C^{2M}):Z\to graph(Z)$.  In
particular it is
invariant under $Z\to Z^{-1}$.

(b) If $m<M$, then the orthogonal projection, $\Cal L(\Bbb C^M)\to
\Cal L(\Bbb C^m):Z\to z$,
maps the measure in (a) to a multiple of
$$det(1+\vert z\vert^2)^{-2m}dLeb(z).$$
This coherence property is equivalent to the assertion that there is
a unique invariant probability on the formal completion of the
infinite rank Grassmannian (a fact which we reproved in $(4.2))$.

(c) The measure $det(1+\vert Z\vert^2)^{-2M-s}dLeb(Z)$ is finite precisely
when $s>-1$; let $\mu_s$ denote the normalized version of this
measure.  If $g=\left(\matrix a&b\\c&d\endmatrix \right)\in U(2M,
\Bbb C)$, where $a$ is $M\times M$, then the
linear fractional action of $g$ on $Z$ has the property that
$$g_{*}(\mu_s)=det\vert a(g^{-1})+b(g^{-1})Z\vert^{2s}d\mu_s;$$
see (3.10) and (3.12) of [Pi1].

\bigskip

\noindent Proof of Proposition (3.3.2).  Let $A=A^{(r)}.$ By Remarks (b)
and (c), to prove finiteness, it suffices to check that
$$det\vert A\vert^{2s}det(1+\vert A\vert^2)^{-2(n-r)-s}dLeb(A)$$
is finite precisely when $s>-1$. This follows from Remark (c), by
changing Z to $A=Z^{-1}.$

To prove the second part of (a), we must show that the projection
$$GL(2N,\Bbb C)\,\,\,\,\,\to\,\,\,\,GL(2n,\Bbb C),$$
defined almost everywhere (more specifically, on the top stratum) by
(3.3.1), maps $d\mu_s^{(N)}$ to $d\mu_s^{(n)}$.

Let $m=N-n$, and write $g\in GL(2N,\Bbb C)$ as
$$g=\left(\matrix a_{11}&a_{12}&a_{13}\\\ a_{21}&a_{22}&a_{23}\\\
a_{31}&a_{32}&a_{33}\\\endmatrix \right),$$
where $a_{11}$ is $m\times m,$ and $a_{22}\in \Cal L(\Bbb C^{2n})$.  The
projection (3.3.1)  is given
by
$$g\quad\to\quad a_{22}-a_{21}a_{11}^{-1}a_{12};$$
it can be factored as the composition of the following four maps:
$$Pr_1:g\to h=\left(\matrix a_{11}&a_{12}\\\ a_{21}&a_{22}\\\endmatrix \right
),$$
$$I_1:h\to k\equiv\left(\matrix b_{11}&b_{12}\\\ b_{21}&b_{22}\\\endmatrix
\right
)=h^{-1}=\left(\matrix *&*\\*&(a_{22}-a_{21}a_{11}^{-1}a_{12})^{-
1}\\\endmatrix \right),$$
$$Pr_2:k\to b_{22},$$
and
$$I_2:b_{22}\to b_{22}^{-1}=a_{22}-a_{21}a_{11}^{-1}a_{12}.$$
We will consider these maps in turn.

(1) Since $det\vert A(g)\vert^2$ $=det\vert A(h)\vert^2$, it follows from (b)
of (3.3.3)  that the
image of d$\mu_s^{(N)}$ under $Pr_1$ is, up to normalization,
$$det\vert A(h)\vert^{2s}det(1+\vert h\vert^2)^{-2(m+2n)-s}dLeb(h
).\tag 3.3.4$$

(2) By (a) of (3.3.3), the effect of inverting h in (3.3.4) is, up to
normalization,
$$det\vert A(k^{-1})\vert^{2s}det\vert k\vert^{2s}det(1+\vert k\vert^
2)^{-2(m+2n)-s}dLeb(k)\tag 3.3.5$$

(3) This is the main step. If
$$k=(b_{ij})=\left(\matrix 1&b_{12}b_{22}^{-1}\\\ 0&1\\\endmatrix \right
)\left(\matrix b_{11}-b_{12}b_{22}^{-1}b_{21}&0\\0&b_{22}\\\endmatrix \right
)\left(\matrix 1&0\\\ b_{22}^{-1}b_{21}&1\\\endmatrix \right),$$
then
$$k^{-1}=\left(\matrix 1&0\\-b_{22}^{-1}b_{21}&1\\\endmatrix \right
)\left(\matrix (b_{11}-b_{12}b_{22}^{-1}b_{21})^{-1}&0\\0&b_{22}^{
-1}\\\endmatrix \right)\left(\matrix 1&-b_{12}b_{22}^{-1}\\0&1\\\endmatrix
\right
).$$
Therefore
$$\frac {detA(k^{-1})}{det(k^{-1})}=\frac {det(b_{11}-b_{12}b_{22}^{
-1}b_{21})^{-1}detA(b_{22}^{-1})}{det(b_{11}-b_{12}b_{22}^{-1}b_{
21})^{-1}det(b_{22}^{-1})}=\frac {detA(b_{22}^{-1})}{det(b_{22}^{
-1})}.$$
Thus the map $k\to b=b_{22}$ projects (3.3.5) to a multiple of
$$det\vert A(b^{-1})\vert^{2s}det\vert b\vert^{2s}det(1+\vert b\vert^
2)^{-4n-s}dLeb(b).\tag 3.3.6$$

(4) Inserting $b^{-1}$ into (II.3.11) now yields $d\mu_s^{(n)}$, by (2). This
proves the
second statement in part (a).

Part (b) follows from (c) of (3.3.3)  and the $SL(2N,\Bbb C)$-equivariance of
the maps in (3.3.1).//

\bigskip

\noindent\S 3.5. Some Projections of $d\mu^{\vert \Cal G\vert^{2s}}$

Fix $n\in \Bbb Z$.  Let $\Cal N^{\pm}_{\infty}$ and $\Cal H_{\infty}$ denote
the formal completions of
the unipotent and Cartan subgroups for $SL(\infty )$, and consider the
model for these groups described in \S 2.1, relative to the basis
$\epsilon_{n+1},\quad\epsilon_{n+2},...$ .  There are natural maps
$$\Cal N^{-}\cdot \Cal H/\Bbb C^{*}.\Cal N^{+}\to \Cal N_{\infty}^{
-}\cdot \Cal H_{\infty}\cdot \Cal N^{+}_{\infty}\to \Bbb M_{\infty}
:(L,D,U)\to (l,d,u)\to\delta =ldu,$$
where $(l_{ij})=(L_{ij})_{i,j>n}$, $d=diag(D_{n+1},..)$, and $(u_{
ij})=(U_{ij})_{i,j>n}$, and
$\delta =ldu$ denotes the product (note that the image of the second map
is described in (b) of (2.1.1)).

Given $N$ such that $n<N$, the projection of $\vert\sigma_n\vert^{
2s}d\mu$ relative to
the map $\delta\to\delta_N$, where $\delta_N=(\delta_{ij})_{n<i,j
\le N}$, is given by
$$\frac 1Edet(1+\vert\delta_N\vert^2)^{-2(N-n)}dLeb(\delta_N),$$
by (3.3.2).  Thus the projection of $\vert\sigma_n\vert^{2s}d\mu$ into $
\Bbb M_{\infty}$ has the
same form as the invariant probability for the infinite rank
Grassmannian, expressed in terms of graph coordinates.

Now consider models for $\Cal N^{\pm}_{\infty}$ and $\Cal H_{\infty}$ relative
to the ordered basis
$\epsilon_n,\epsilon_{n-1},..$ . There are natural maps

$$\Cal N^{-}\cdot \Cal H/\Bbb C^{*}\cdot \Cal N^{+}\to \Cal N_{\infty}^{
+}\cdot \Cal H_{\infty}\cdot \Cal N_{\infty}^{-}\to \Bbb M_{\infty}
:(L,D,U)\to (u^{-1^{*}},d^{-1^{*}},l^{-1^{*}})\to\alpha =ldu,$$
where $(u^{-1^{*}}_{ij})=(L_{ij})_{i,j\le n}$, $d^{-1^{*}}=diag(D_
n,D_{n-1},..)$, and
$(l^{-1^{*}})=(U_{ij})_{i,j\le n}$.  It follows from the way in which the local
expression for $\vert\sigma_n\vert^{2s}d\mu$ transforms under inversion (see
(2) of
the proof of (3.3.2)) that the image under the above map again
has the same form as the invariant probability for the infinite
rank Grassmannian, expressed in terms of graph coordinates.

The significance of these calculations will be explored in the
next subsection.

\bigskip

\flushpar\S 3.6. Finite Dimensional Approximations to $d\mu^{\vert
\Cal G\vert^{2s}}$.

We first consider the case $s=0$, in which case we write
$\mu_0=\mu^{\vert \Cal G\vert^0}$.

\proclaim{(3.6.1)Proposition} Let $\delta_{SU(2N)}$ denote the invariant
probability measure on $SU(2N)\hookrightarrow \Cal G/\Bbb C^{*}$. Then
$$\delta_{SU(2N)}\to\mu_0\quad as\quad N\to\infty$$
weakly with respect to bounded continuous functions on the top
stratum (hence also with respect to $BC(\Cal G/\Bbb C^{*})$).
\endproclaim

\flushpar(3.6.2) Remark.  The important point here is that it is
not necessary to scale the position of $SU(2N)$ inside $\Cal G/\Bbb C^{
*}$, in
contrast to part (c) of (2.3.1), where it was necessary to scale
the position of $SU(N)$ inside of $SL(\infty )_{formal}$. The role of
scaling will emerge quite clearly in the course of the proof.

\smallskip

\flushpar Proof of (3.6.1).  Given $N,M\ge 0$, identify $\Bbb C^N
\oplus \Bbb C^M$
with
$$\sum_{-N<j\le 0}\Bbb C\epsilon_j\quad\oplus\quad\sum_{0<j\le M}
\Bbb C\epsilon_j.\tag 3.6.4$$
 Let $\delta_{N,M}$ denote the invariant probability on the space
$$(\frac MN)^{1/2}SU(N+M)\hookrightarrow GL(N+M)\hookrightarrow \Cal G
/\Bbb C^{*}.\tag 3.6.5$$
We will prove the more general fact that
$$\lim_{N,M\to\infty}\delta_{N,M}=d\mu_0,\tag 3.6.6$$
in the weak sense with respect to bounded continuous
functions on $\Cal G/\Bbb C^{*}$.  Our proof of the existence of the limit
will depend upon some facts that we will establish in \S 4.4.
The main point here will be to prove that the iterated limits
equal the invariant measure.

We first discuss existence of the limit. Consider the top
stratum
$$\Cal N^{-}\cdot \Cal H/\Bbb C^{*}\cdot \Cal N^{+}\subset \Cal G
/\Bbb C^{*}.\tag 3.6.7$$
We have already established weak convergence for functions that
depend only upon $\Cal N^{\pm}$ (this follows from the fact that the
projections of the measures to $\Cal N^{\pm}$ are coherent).  The key is
to analyze the diagonal distribution.

Now there is an identification
$$\prod_{-\infty}^{+\infty}\Bbb C^{*}\to \Cal H/\Bbb C^{*}:(z_j)\to
diag(z_j).\tag 3.6.8$$
We know that the diagonal distribution will be invariant with
respect to the unitary part of this group.  Write $a=(a_j)$,
where $a_j=\vert z_j\vert$, and regard $a$ as a map on the top stratum, by
composing with the projection to $\Cal H/\Bbb C^{*}$.  To prove weak
convergence, we must prove that the Fourier transform for
the diagonal distribution,
$$(a_{*}\delta_{N,M}\hat {)}(\lambda )\equiv\int_{SU(N+M)}e^{-i\lambda
(log(a((\frac NM)^{1/2}k)))}dm(k),\tag 3.6.9$$
has a limit for each fixed $\lambda$ in
$$(\prod_{-\infty}^{+\infty}\Bbb R^{+}\hat {)}=\sum_{-\infty}^{+\infty}
\Bbb R.\tag 3.6.10$$
(a crucial point here is that $\lambda$ has only finitely many nonzero
components - this is one of the advantages of working with
formal completions; see chapter 1 of [Bill1]).  Now according to
Harish-Chandra (see (4.4.27) of the next section), (3.6.9) is equal
to
$$exp(-\frac i2log(\frac NM)\sum\lambda_j)\prod_{-N<j<k\le M}(1-\frac
i2\frac {\lambda_k-\lambda_j}{k-j})^{-1}.\tag 3.6.11$$
Since there are only finitely many nonzero $\lambda_j$, the exponential
factor regularizes the product, so that (3.6.11) has the following
limit as $N,M\to\infty$,
$$(\prod_{-\infty <j<k<+\infty}(1-\frac i2\frac {\lambda_k-\lambda_
j}{k-j})exp(\frac i2\frac {\lambda_j}{j-k}))^{-1}\tag 3.6.12$$

It follows now that the projection of $\{\delta_{N,M}\}$ to each of the
spaces $\Cal N^{\pm}$ and $\Cal H/\Bbb C^{*}$ has a unique weak limit as $
N,M\to\infty$.  Thus
there will be limit points for $\{\delta_{N,M}\}$ on the top stratum.  But
any such limit point will be $SU(2\infty )$ biinvariant, and the
Fourier transform of the diagonal distribution will be given by
the limit of (3.6.11).  In \S 4.4 we will prove that a unitarily
biinvariant probability is uniquely determined by its diagonal
distribution; see (4.4.9).  This implies that the limit point is
uniquely determined.  Thus $\{\delta_{N,M}\}$ has a unique limit as
$N,M\to\infty$.

We will now identify the limit.

Fix $N$.  Consider the map in the first paragraph of \S 3.5,
$\Cal N^{-}\cdot \Cal H/\Bbb C^{*}\cdot \Cal N^{+}\to \Bbb M_{\infty}$, with $
N$ in place of $n$.  By part (c) of (2.3.1)
the limit of $\delta_{N,M}$, as $M\to\infty$, is $(N^{-1/2})_{*}\nu_
1^{(N)}$, where $\nu_1^{(N)}$ is
the standard Gaussian, viewed as a measure on the embedded
subspace $\Cal N^{-}_{\infty}\cdot \Cal H_{\infty}\cdot \Cal N^{+}_{
\infty}\hookrightarrow \Cal N^{-}\cdot \Cal H/\Bbb C^{*}\cdot \Cal N^{
+}$, where the embedding
depends on N, as in the first paragraph of \S 3.5.

The following result is of independent interest, because it
shows that $\mu_0$ is a limit of Gaussians.

\proclaim{(3.6.13) Lemma} $(N^{-1/2})_{*}\nu_1^{(N)}$ converges to $
d\mu_0$ as
$N\to\infty$.  \endproclaim

Assuming the Lemma, we can conclude that
$$\lim_N\,\,\lim_M\,\,\delta_{N,M}=d\mu_0.\tag 3.6.14$$
This will conclude the proof of (3.6.1) (note that we can also
directly conclude that
$$\lim_M\,\,\lim_N\,\,\delta_{N,M}=d\mu_0;$$
the line of argument differs only in that one considers the
map in the second paragraph of \S 3.5, and one applies a
corresponding reflection principle to the proof of Lemma
(3.6.13) below).

\bigskip

\flushpar Proof of Lemma (3.6.13).  Project $\nu_1^{(N)}$ onto $G
L(2N,\Bbb C)$
via the map in (3.3.1).  The image is again the standard
Gaussian.  Now fix $n<N$.  Write $g\in GL(2N,\Bbb C)$ as a $3\times
3$
matrix, $g=(a_{ij}),$ precisely as in the proof of (3.3.2) (so $a_{
22}$ is
a $2n\times 2n$ matrix).  We must show that if g is distributed
according to the standard Gaussian, then the random matrix
$$N^{-1/2}(a_{22}-a_{21}a_{11}^{-1}a_{12})\tag 3.6.15$$
has limiting distribution $d\mu_0^{(n)}$.

View $a_{21}a_{11}^{-1},\quad a_{22}^{-1}a_{21}\in \Cal L(\Bbb C^
m,\Bbb C^{2n})\hookrightarrow Gr(m,\Bbb C^{m+2n})$.  These random
variables are distributed according to the invariant probability
on the Grassmannian (see section 1 of [Pi]).  Thus the limit of
(3.6.15) equals the limit of
$$N^{-1/2}(a_{22}^{-1}a_{21}a_{12}).\tag 3.6.16$$
The central limit theorem implies that $N^{-1/2}a_{21}a_{12}$ converges
to the standard Gaussian.  The matrix $a_{22}$ is independent of
this limit, and it also has the standard Gaussian distribution.
Thus the limit of (3.6.16) is $d\mu_0^{(n)}$.//

\flushpar(3.6.17)Remark. Just as in the case of $SL(\infty )$, it is
possible to check (3.6.12) directly for the individual
distributions of the $\sigma_j$.

\bigskip

\flushpar\S 3.7. Nonexistence of invariant measures on $\Cal G$.

\proclaim{(3.7.1)Proposition}There does not exist a left
$SU(2\infty )$-invariant probability on $\Cal G$, with the property that the
the matrix coefficients for the fundamental modules are
square-integrable.
\endproclaim

\flushpar Proof.  By uniqueness of invariant measures on flag
spaces, the existence of such a measure would imply that the
measures $\mu^{\vert \Cal L_{\Lambda}\vert^{2s}}$ on $\Cal F'$ are mutually
absolutely continuous, for
different $\Lambda$ and nonnegative integral $s$, and for any generalized
flag space $\Cal F'$.  But it is known that for the minimal flag
spaces, the Grassmannians, these different measures are
actually mutually singular ([Pi5]).//

\bigskip

\flushpar\S 3.8. Some Additional Comments.

Given the existence of the measures $\mu^{\vert \Cal G\vert^{2s}}$, we can
proceed to
define a number of unitary representations of $SU(2\infty ,\Bbb C
)$.   For
any $s>-1$ there is a natural unitary action of $SU(2\infty ,\Bbb C
)$ on the
(incomplete) space of square-integrable densities
$$\{\theta\in\Omega^0(\vert \Cal G^{*}\vert^s):\int\theta\bar{\theta }
d\mu^{\vert \Cal G\vert^{2s}}<\infty \}.\tag 3.8.1$$
By using the density $\vert\sigma_0\vert^s$ to trivialize $\vert
\Cal G^{*}\vert^s$ over $\{\sigma_0\ne 0\}$, this
action is identified with the restriction to a dense subspace of
the action
$$SU(2\infty ,\Bbb C)\times L^2(d\mu_s)\to L^2(d\mu_s):g,f\to\rho
(g,\cdot )^sf(g^{-1}\cdot (\cdot ))\tag 3.8.2$$
where
$$\rho (g,\cdot )=\vert\frac {g_{*}\sigma_0(\cdot )}{\sigma_0(\cdot
)}\vert .\tag 3.8.3$$
If $g\in SU(2n,\Bbb C)$, then $\rho (g,\cdot )$ is a function of finitely many
variables, i.e.  $\rho (g,\cdot )$ is the pullback of a function from
$GL(2n,\Bbb C)$ with respect to the horizontal arrows in (3.3.1).  Thus
the subspace of functions of finitely many variables is dense
and invariant.

If $s=k$ is a nonnegative integer, then we can consider the
unitary action of $SU(2\infty ,\Bbb C)$ on the (incomplete) space of
square-integrable sections of $(\Cal G^{*})^{\otimes k}$.  The subspace of
holomorphic sections is an invariant subspace.  By (2.1.6) of
Part I, every holomorphic section of $(\Cal G^{*})^{\otimes k}$ is regular,
meaning
that it is a function of finitely many variables, in the sense of
the preceding paragraph.  Thus all holomorphic sections of
$(\Cal G^{*})^{\otimes k}$ are square-integrable.  Any $L^2$ limit of
holomorphic
sections of $(\Cal G^{*})^{\otimes k}$ can be realized as a holomorphic section
of
the restriction of this bundle to $GL_{(2)}$; this can be proved by
modifying the argument for Grassmannians, see theorem (6.2) of
[Pi1].  It seems doubtful that there is a corresponding
attractive functional analytic realization of the orthogonal
complement of $H^0((\Cal G^{*})^{\otimes k})$.  Nonetheless we can assert the
following

\proclaim{(3.8.4)Proposition}All of the unitary representations
above extend continuously to strong operator continuous
unitary representations of $\tilde {U}_{(2)}$, the Kac-Moody completion of
$SU(2\infty ,\Bbb C)$ (see (0.1.1)).  \endproclaim

\flushpar Proof. This follows essentially from the fact that for
each of the representations above, the subspace of functions of
finitely many variables is dense. For any such function is fixed
by the subgroup
$$\{g=\left(\matrix a&0\\0&d\endmatrix \right)\in SU(2\infty ,\Bbb C
):g\epsilon_j=\epsilon_j,\vert j\vert\le N\},$$
for some $N$. It follows automatically from this condition that
the representations extend continuously as claimed; see
Proposition 6.11 of [Pi2].//

\bigskip

\centerline{ \S 4. The Cases $\frak g=o(2\infty ,\Bbb C),$ $o(2\infty
+1,\Bbb C),$ and $sp(\infty ,\Bbb C)$}

\bigskip

\flushpar\S 4.1. Realizations for $\frak g=o(2\infty ,\Bbb C)$.

Let $\{\epsilon_j:j\in \Bbb Z+1/2\}$ denote a basis for a complex vector space
$H_{alg}$, and identify $\Bbb C^n\oplus \Bbb C^m$ with
$$\sum_{-n-1<j<0}\Bbb C\epsilon_j\quad\oplus\quad\sum_{0<j<m+1}\Bbb C
\epsilon_j.\tag 4.1.1$$
Let $(\cdot ,\cdot )$ denote the symmetric complex bilinear form defined
by
$$(\epsilon_i,\epsilon_j)=\delta_0(i+j).\tag 4.1.2$$
The transpose with respect to this form is given by
$$(x^t)_{ij}=x_{-j,-i}.\tag 4.1.3$$
We can identify $SO(2n,\Bbb C)$ with those transformations which
preserve $(\cdot ,\cdot )$, and which fix the complement of $\Bbb C^
n\oplus \Bbb C^n$.  We
make the usual choice of diagonal Cartan subgroup and positive
roots (that is to say, the positive Borel subalgebra consists of
upper triangular matrices in $so(2n,\Bbb C))$.  We will index the
positive simple roots by the natural numbers; the
corresponding coroots $h_j$ are given by
$$h_1=diag(\cdot\cdot ,0,1,1;-1,-1,0,\cdot\cdot ),h_2=diag(..,0,1
,-1;1,-1,..),\tag 4.1.4$$
$$h_3=diag(..,0,1,-1,0;0,1,-1,0,..),h_4=diag(..,0,1,-1,0,0;0,0,1,
-1,0,..)$$
and so on, where the semicolon marks the transition from
negative to positive indices.

Let $\Bbb M$ denote the set of all matrices indexed by $\Bbb Z+1/
2$.  Write
$g\in \Bbb M$ as
$$g=\left(\matrix A&B\\C&D\endmatrix \right),\tag 4.1.5$$
where $D$ is indexed by positive $j$. The description of the
spaces $G/\Bbb C^{*}$ and $\Cal G^{\pm}/\Bbb C^{*}$ inside $\Bbb M$ is the same
as in the case of
$sl(2\infty ,\Bbb C)$, with the orthogonality condition thrown in. In
particular
$$\Cal H=\prod H_{(j)}\to \Bbb M:\prod a_j^{h_j}\to diag(d_j),\tag 4.1.6$$
where $d_{-1/2}=a_1/a_2$, $d_{-3/2}=a_1a_2/a_3$, $d_{-(2j-1/2)}=a_
j/a_{j+1}$,
$j>2$. The map of the center, $\Bbb C^{*}$, into $\Cal H$ is given by
$$\Bbb C^{*}\to \Cal H:a\to (a_j),\tag 4.1.7$$
where $a_1=a_2=a$, $a_j=a^2$, $j>2$.  From this one sees that
$\sigma_1=a_1$ and $\sigma_2=a_2$ (the fundamental matrix coefficients for the
spin representations) have level 1, while $\sigma_j=a_j$ has level 2 for
$j>2$.

\bigskip

\S 4.2. Realizations for $\frak g=o(2\infty +1,\Bbb C)$.

In this case we take a basis $\{\epsilon_j:j\in \Bbb Z\}$ for $H_{
alg}$. We define a
symmetric complex bilinear form $(\cdot ,\cdot )$ on $H_{alg}$ by the formula
(4.1.2). The formula for the transpose is given by (4.1.3). The
coroots are
$$h_1=diag(..,0,2;0;-2,0,..),h_2=diag(..,0,1,-1;0;1,-1,..),\tag 4.2.3$$
and so on.

In this case $\Bbb M$ will denote all matrices indexed by $\Bbb Z$, and we
will write a general $g\in \Bbb M$ as
$$g=\left(\matrix A&u&B\\l&m&r\\C&d&D\endmatrix \right),\tag 4.2.2$$
where $A$ ($D$, resp) is indexed by negative integers (positive
integers, resp). The realizations of $G/\Bbb C^{*}$, etc. are as before.
The map of $\Cal H$ to $\Bbb M$ is given by
$$\Cal H\to \Bbb M:\prod a_j^{h_j}\to diag(d_j),\tag 4.2.3$$
where $d_{-1}=a_1^2/a_2$, $d_{-j}=a_j/a_{j+1}$, $j>1$. The center, $
\Bbb C^{*}$, maps
into $\Cal H$ by
$$\Bbb C^{*}\to \Cal H:a\to (a_j),\tag 4.2.4$$
where $a_1=a$, and $a_j=a^2$, for $j>0$. In this case $\sigma_1$ is of level
1, and all of the other $\sigma_j$ are of level 2.

\bigskip

\flushpar\S 4.3. Realizations for $\frak g=sp(\infty ,\Bbb C)$.

In this case we take a basis $\{\epsilon_j:j\in \Bbb Z+1/2\}$ for $
H_{alg}$.  The
skew complex bilinear form is given by
$$(\epsilon_i,\epsilon_j)=sign(i)\delta_0(i+j).\tag 4.3.1$$
The transpose is given by
$$(x^t)_{ij}=sign(i)sign(j)x_{-j,-i}\tag 4.3.2$$
The coroots are given by
$$h_1=diag(..,0,1;1,0,..),h_2=diag(..,0,1,-1;1,-1,0,..)\tag 4.3.3$$
and so on.

We write a $g\in \Bbb M$ as a $2\times 2$ matrix as in (4.1.5). The map of $
\Cal H$ to
$\Bbb M$ is given by
$$\Cal H\to \Bbb M:\prod a_j^{h_j}\to diag(d_j),\tag 4.3.3$$
where $d_{-j}=a_j/a_{j+1}$, for all $j$. The center embedding is
$$\Bbb C^{*}\to \Cal H:a\to\prod a^{h_j},\tag 4.3.4$$
and $\sigma_j$ has level one for all $j$.

\bigskip

\flushpar\S 4.4. The projection of a biinvariant measure to the
diagonal.

\smallskip

In this subsection $G$ will denote a finite-dimensional semisimple
simply connected complex group.  We will show that a
$K$-biinvariant probability on $G$ is determined by its projection
to the diagonal.  This and various explicit formulas are simple
consequences of Harish-Chandra's theory of spherical analysis
for $G$. We will use [GV] and [Helg] as basic references.

As usual we fix a triangular decomposition
$$\frak g=\frak n^{-}\oplus \frak h\oplus \frak n^{+}.\tag 4.4.1$$
Let $\{h_j\}$ and $\{\sigma_j\}$ denote the corresponding sets of coroots and
fundamental matrix coefficients for $G$, respectively.  If
$\prod\sigma_j(g)\ne 0$, then we can write $g\in G$ uniquely as
$$g=l\cdot diag\cdot u,\tag 4.4.2$$
where $l$ (resp.  $u$) is lower (resp.  upper) triangular, and
$diag=\prod\sigma_j(g)^{h_j}\in H$.

We have
$$H=T\cdot A,\tag 4.4.3$$
where $T=K\cap H$, $\frak a=\frak h_{\Bbb R}$, and
$$exp:\frak a\to A\tag 4.4.4$$
is an isomorphism.  Let $\rho$ denote the sum of the positive roots
(for the complex linear adjoint action of $\frak h$ on $\frak g$).  We can
write Haar measure for $G$, with respect to the decomposition
(4.4.2), as
$$dm(g)=dm(l)\times dm(m)\times a^{2\rho}dm(a)\times dm(u),\tag 4.4.6$$
where $diag(g)=m\cdot a$ with respect to (4.4.3), and $dm(l)$ denotes
a Haar measure for $N^{-}$, and so on.

Consider the bounded linear map
$$L^1(G)\to L^1(A):f\to F,\tag 4.4.7$$
where
$$F(a)=a^{2\rho}\iint f(lmau)dm(m)dm(l)dm(u).\tag 4.4.8$$

\proclaim{(4.4.9)Proposition}The map (4.4.7) induces an injective
transformation
$$\Cal T:L^1(K\backslash G/K)\to L^1(A).\tag 4.4.10$$
\endproclaim

The map $\Cal T$ actually extends to an injective map of spaces of
finite Borel measures, and the inversion formula can be made
quite explicit; see Remark (4.4.32) below.

\flushpar Proof of (4.4.9).  By slight abuse of notation, we will
write $L^1(A,$ $a^{\rho}dm)^W$ for the space of functions on $A$ which
are $W$-invariant and $L^1$ with respect to $a^{\rho}dm(a)$ (it is not
implied that this $L^1$ space is $W$-invariant).  The map $\Cal T$ is a
composition $\Cal B\circ \Cal A$, where $\Cal A$ is the Abel transform
$$\Cal A:L^1(K\backslash G/K)\to L^1(A,a^{\rho}dm)^W:f\to a^{\rho}
\int_{N^{+}}f(au)dm(u)\tag 4.4.12$$
and
$$\Cal B:Image(\Cal A)\to L^1(A):f_1\to a^{\rho}\int_{N^{-}}f_1(\alpha
(la))dm(l),\tag 4.4.13$$
where we have written $g=\kappa\cdot\alpha\cdot n$ corresponding to the
Iwasawa decomposition $G=K\cdot A\cdot N^{+}$.  The injectivity of the
Abel transform is basically a general fact, because when it is
composed with the Fourier transform for $A$, it becomes
Harish-Chandra's realization of the Gelfand transform for the
semisimple commutative Banach $*$-algebra $L^1(K\backslash G/K)$; this also
explains why the image of $\Cal A$ consists of $W$-invariant functions
(see \S 3.3 of [GV]).  Thus it suffices to show that $\Cal B$ is
injective.

Note that $\alpha (la)=\alpha (l)a$.  We claim that for $\alpha :
N^{-}\to A$, $\alpha_{*}dm(l)$
is absolutely continuous with respect to $dm(a)$, hence
$$\aligned
\int_{N^{-}}f_1(\alpha (la))dm(l)=\int_Af_1(a\alpha )\Cal D(\alpha
)dm(\alpha )=(f_1*_A\Cal D^{*})(a)\endaligned
\tag 4.4.14$$
where $\Cal D$ is the Radon-Nikodym derivative, $\Cal D^{*}(\cdot
)=\Cal D((\cdot )^{-1})$, and
$*_A$ denotes convolution in $A$.  To justify this, consider the
Iwasawa factorization $l=\kappa\alpha n$.  We have the Bruhat
decomposition $\kappa =l\alpha^{-1}(\alpha n^{-1}\alpha^{-1})$, and this
implies that
$$\alpha (l)=\prod\sigma_j(\kappa (l))^{-h_j}.\tag 4.4.15$$
It follows that $\sigma_j(\kappa (l))$ is real, and it is always $
\le 1$ in
magnitude since $\kappa$ is unitary.  Now the range of the injective map
$$N^{-}\to \Cal F\equiv K/T:l\to\kappa (l)T\tag 4.4.16$$
is precisely the top stratum of the flag space.  Also the map
$$\Cal F\equiv K/T\to \{\vert z\vert\le 1\}^{rk(\frak g)}:k\to (\vert
\sigma_j(k)\vert )\tag 4.4.17$$
is surjective, and the inverse image of (0,..,0) is precisely the
complement of the top stratum.  Thus
$$\alpha :N^{-}\to A^{+}\equiv \{\alpha\in A:\sigma_j(\alpha )\ge
1,\forall j\}\tag 4.4.18$$
is a smooth surjective map, and this implies that $\alpha_{*}dm(l
)$ is
absolutely continuous. This justifies (4.4.14).

In the remainder of the proof, we will use the exponential map
to identify $\frak a$ with $A$.  We will denote the Fourier transform
by $\hat {f}$.  On a heuristic level, the proof is completed by
observing that if $f_1*\Cal D^{*}=0$, then $\hat {f}_1(\Cal D^{*}
\hat {)}=0$; $\Cal D^{*}$ is defined in a
wedge, hence $(\Cal D^{*}\hat {)}$ is nonzero a.e.,  because it is the boundary
values of a holomorphic function; hence $f_1=0$ a.e.  Below we
will observe that $(\Cal D^{*}\hat {)}$ is expressible in terms of
Harish-Chandra's {\bf c}-function, and this completes the proof.
However we will first present a less sophisticated argument.

We can write
$$e^{\rho (x)}=e^{2\sum x_j}\tag 4.4.19$$
where $x=\sum x_jh_j\in \frak a$. The positive Weyl chamber is
contained in $\{x\in \frak a:x_j\ge 0,\forall j\}$.  Because $f_1$ and $
f_1*\Cal D^{*}$ belong to
$L^1(\frak a,e^{\rho (x)}dx)^W$,
$$\hat {f}_1(\lambda )=\int_{\frak a}f_1(x)e^{-i\lambda (x)}dx\tag 4.4.20$$
is defined and holomorphic in the region $\Cal R=$
$\{\lambda\in \frak a_{\Bbb C}^{*}:Im(\lambda )^w(h_j)<2,\forall
j,\forall w\in W\}$, and similarly for $f_1*\Cal D^{*}$.
Because $Image(\Cal A)$ contains the Schwarz space, it follows that
the function
$$(\Cal D^{*}\hat {)}\equiv\frac {(e^{-\langle x,x\rangle /2}*\Cal D^{
*}\hat {)}}{(e^{-\langle x,x\rangle /2}\hat {)}}\tag 4.4.21$$
is holomorphic in the region $\Cal R$ (we will justify the use of the
notation $(\Cal D^{*}\hat {)}$ below; see (4.4.25)).   Because $I
mage(\Cal A)$ is a
convolution algebra, we have an equality of holomorphic
functions in $\Cal R$,
$$(f_1*\Cal D^{*}\hat {)}=(f_1*e^{-\langle x,x\rangle /2}*\Cal D^{
*}\hat {)}/(e^{-\langle x,x\rangle /2}\hat {)}=\hat {f}_1(\Cal D^{
*}\hat {)}.\tag 4.4.22$$
It follows that if $f_1*\Cal D^{*}=0$, for some $f_1\in Image(\Cal A
)$, then
$f_1=0$. This implies that $\Cal B$ is injective, completing the
proof.//

A far deeper analysis can be given as follows.  A basic
estimate of harmonic analysis on $G$ is that
$$\int_{N^{-}}\alpha (l)^{-(1+\epsilon )\rho}dm(l)<\infty\tag 4.4.24$$
for $\epsilon >0$, and this diverges logarithmically if $\epsilon
=0$ (see Thm
4.7.3 of [GV]). Thus we have
$$(\Cal D^{*}\hat {)}(\lambda )=\int_{\frak a}e^{-i\lambda (x)}\Cal D
(-x)dx=\int_{N^{-}}\alpha (l)^{i\lambda}dm(l)\equiv \bold c(-i\lambda
-\rho )\tag 4.4.25$$
for $Im(\lambda )(h_j)>2$ (see Thm 4.7.4 of [GV], and note that our
Fourier transform convention differs slightly from [GV]'s).  But
for complex $G$, the {\bf c}-function is rational,
$$\bold c(\lambda )=\frac {\pi (\rho )}{\pi (\lambda )},\quad whe
re\quad\pi (\cdot )=(\prod_{\alpha >0}\langle\cdot ,\alpha\rangle
)\tag 4.4.26$$
(see Thm 5.7 of Chapter IV of [Helg]). From our point of view
this means that $(\Cal D^{*}\hat {)}$, defined by (4.4.25), has a meromorphic
extension to all of $\frak a_{\Bbb C}^{*}$, and this is consistent with our
definition of
$(\Cal D^{*}\hat {)}$ in (4.4.21) above. The explicit formulas (4.4.25) and
(4.4.26)
for $(\Cal D^{*}\hat {)}$ also imply the following

\proclaim{(4.4.27)Proposition}$(\Cal D^{*}\hat {)}(\cdot )=\bold c
(-\rho -i(\cdot ))$ is uniformly
bounded away from zero on $\frak a^{*}$.  If $f\in C_{cpt}(K\backslash
G/K)$, then the
Harish transform $\Cal Hf$ is given by
$$\aligned
\Cal Hf(\lambda )&=\bold c(-\rho -i\lambda )^{-1}(a^{-\rho}\Cal T
f\hat {)}(\lambda )\\&=\bold c(-\rho -i\lambda )^{-1}(\Cal Tf\hat {
)}(\lambda -i\rho )\endaligned
.$$
\endproclaim\

This reduces the inversion of $\Cal T$ to the inversion of $\Cal H$, which is
the central issue in both [GV] and [Helg].

The {\bf c}-function also arises in the calculation of the projection
of the Haar measure for $K$ to the diagonal. Write $a(k)$ for the
$A$ component of $diag(k)$ in (4.4.2).

\proclaim{(4.4.28)Proposition}We have
$$(a_{*}dm(k)\hat {)}(\lambda )=\bold c(\rho -i\lambda ).$$
More generally, given a dominant integral functional $\Lambda$, if $
\sigma_{\Lambda}$
is a corresponding lowest weight, viewed as a matrix
coefficient, then
$$\int_Ka(k)^{-i\lambda}\vert\sigma_{\Lambda}(k)\vert^{2s}dm(k)=\bold c
(\rho +2s\Lambda -i\lambda )\int\vert\sigma_{\Lambda}(k)\vert^{2s}
dm(k).$$
\endproclaim

\flushpar Proof.  We can normalize the Haar measure for $N^{-}$ so
that with respect to the coordinate (4.4.17) for the top
stratum of $K/T$, the invariant probability is given by
$$\alpha (l)^{-2\rho}dm(l),\tag 4.4.29$$
where $l=\kappa\alpha n$ is the Iwasawa decomposition.
It follows that
$$\int_Ka(k)^{-i\lambda}dm(k)=\int_{N^{-}}a(\kappa (l))^{-i\lambda}
\alpha (l)^{-2\rho}dm(l)$$
$$=\int_{N^{-}}\alpha (l)^{-(-i\lambda +2\rho )}dm(l)=\bold c(\rho
-i\lambda ).\tag 4.4.30$$

The second part follows along the same lines, using
$$\vert\sigma_{\Lambda}(g)\vert^{2s}=a(g)^{2s\Lambda}.//\tag 4.4.31$$

\smallskip

\flushpar(4.4.32)Remark. It is possible to formulate a more
general result than (4.4.9) as follows. Harish-Chandra's theory
implies that the Abel transform actually induces an
isomorphism
$$\Cal M(K\backslash G/K)\to \Cal M(A)^W,$$
where $\Cal M$ denotes the space of finite Borel measures (in fact,
more generally still, the Abel transform induces an
isomorphism of spaces of tempered distributions, where the
Schwarz space of $G$ must be suitably defined; see [GV]).  The
map $\Cal B$ extends to a map
$$\Cal M(A)^W\to \Cal M(A):\nu\to \Cal D^{*}*_A\nu ;$$
this map is injective by the first statement in (4.4.27), and in
practice we can use the second part of (4.4.27) and a limiting
procedure to invert the map.

\bigskip

\flushpar\S 4.5. Existence of the measure $\mu_0$.

In this subsection we will give a uniform proof of the
existence of a $K(\infty )$-biinvariant measure $\mu_0$ on $\Cal G
/\Bbb C^{*}$.

Let $\bar {G}(l)$ denote $SO(2l,\Bbb C),$ $SO(2l+1,\Bbb C)$, or $
Sp(l,\Bbb C)$, let $\bar {H}(l)$
denote the diagonal Cartan subgroup, and let $\bar {K}(l)$ denote the
unitary form
$$\bar {K}(l)=\bar {G}(l)\cap U(H),\tag 4.5.1$$
where $H$ is the Hilbert space completion of $H_{alg}$ gotten by
declaring $\{\epsilon_j\}$ to be an orthonormal basis.  In each case, there
is a natural inclusion
$$U(l,\Bbb C)\hookrightarrow\bar {K}(l):k\to g,\tag 4.5.2$$
where $A(g)=k$, in the notation of the preceding subsections.

Let $\delta_{\bar {K}(l)}$ denote the invariant probability on $\bar {
K}(l)$.

\proclaim{(4.5.3)Theorem}The measures $\delta_{\bar {K}(l)}$ converge weakly as
$l\to\infty$ to a $K(\infty )$ biinvariant probability measure on $
\Cal G/\Bbb C^{*}$.
\endproclaim

\flushpar Proof.  It suffices to check weak convergence for the
top stratum $\Cal N^{-}\cdot \Cal H/\Bbb C^{*}\cdot \Cal N^{+}\subset
\Cal G/\Bbb C^{*}$ (the argument here is
identical to the one in the proof of (3.6.1), see the paragraph
following (3.6.7)).  The key, as in the proof of (3.6.1), is to
check weak convergence for the diagonal distributions.

Our approach to proving this is the same as in the proof of
(3.6.1), modulo notational details.  There is an identification
$$\prod_1^{\infty}\Bbb C^{*}\to \Cal H/\Bbb C^{*},\tag 4.5.4$$
where $(z_j)$ maps to $diag(d_k)$, and for types $C$ and $D$,
$d_{-1/2}=z_1$, $d_{-3/2}=z_2$,..,  and for type $B$, $d_0=1$, $d_{
-1}=z_1$,..
The diagonal distribution will be invariant with respect to the
unitary part of this group.  Write $a=(a_j)$, where $a_j=\vert z_
j\vert$,
and regard $a$ as a function on the top stratum, by composing
with the projection to $\Cal H/\Bbb C^{*}$.  To prove weak convergence of
the diagonal distributions is equivalent to proving pointwise
convergence of the Fourier transforms,
$$(a_{*}\delta_{\bar {K}(l)}\hat {)}(\lambda )=\int_{\bar {K}(l)}
a(k)^{-i\lambda}dm(k),\tag 4.5.5$$
for $\lambda$ in the space
$$(\prod_1^{\infty}\Bbb R^{+}\hat {)}=(\sum_1^{\infty}\Bbb R).\tag 4.5.6$$
To apply (4.4.28), we calculate that for type $D$, we have
$\rho_j=2(j-1)$; for type $C$, we have $\rho_j=2j$; for type $B$, we have
$\rho_j=2j-1$.  We then obtain the following formulas for (4.5.5),
$$\prod_{0<j<k\le l}(1-\frac i2\frac {\lambda_k-\lambda_j}{k-j})^{
-1}\prod_{0<p<q\le l}(1-\frac i2\frac {\lambda_p+\lambda_q}{p+q-1}
)^{-1},$$
$$\prod_{0<j<k\le l}(1-\frac i2\frac {\lambda_k-\lambda_j}{k-j})^{
-1}\prod_{0<p\le q\le l}(1-\frac i2\frac {\lambda_p+\lambda_q}{p+
q})^{-1},\tag 4.5.7$$
$$\prod_{0<j<k\le l}(1-\frac i2\frac {\lambda_k-\lambda_j}{k-j})^{
-1}\prod_{0<p<q\le l}(1-\frac i2\frac {\lambda_p+\lambda_q}{p+q-1}
)^{-1}\prod_{0<r\le l}(1-i\frac {\lambda_r}{2r-1})^{-1},$$
for types $D$, $C$, and $B$, respectively.  Each of these does have a
pointwise limit as $l\to\infty$.  To understand this, observe that if
we fix $j$, then for large $k$, because $\lambda_k=0$, from the first
factor we obtain a term of the form
$$(1-\frac i2\frac {-\lambda_j}{k-j})^{-1},\tag 4.5.8$$
and from the second factor we obtain the conjugate of this
term by taking $p=j$, and solving for $q$.  The diagonal in the
second factor for type $C$ and the third factor for type $B$ are
independent of large $l$, because $\lambda_r=0$ for large $r$.  Thus the
infinite products converge.

We can now complete the proof as follows.  We know that the
$\Cal N^{\pm}$ and diagonal distributions of $\{\delta_{\bar {K}(
l)}\}$ have weak limits as
$l\to\infty$.  This implies the existence of weak limit points for this
sequence on the top stratum.  But these limit points, regarded
as measures on $\Cal G/\Bbb C^{*}$, will be $K(\infty )$-biinvariant, and by
(4.4.4)
any such measure is determined by its diagonal distribution.
This implies that there is a unique limit point.//

\smallskip

It would clearly be desirable to give a direct computational
proof of this result, as we did for type $A$ in the previous
section. This would involve finding explicit expressions
for the projection of $\mu_0$ to $G(l)$ in these cases, and I have
succeeded at this.

\bigskip

\flushpar\S 4.6. Existence of biinvariant measures with values in
$\vert \Cal G\vert^{2s}$.

The inclusion of $G(\infty )\subset \Cal G$ defines a natural trivialization of
the
restriction of the $\Bbb C^{*}$-bundle $\Cal G\to \Cal G/\Bbb C^{
*}$; explicitly, the
trivialization is given by
$$G(\infty )\times \Bbb C^{*}\to \Cal G:g,z\to g\cdot\iota (z)=\iota
(z)\cdot g\tag 4.6.1$$
where $\iota :\Bbb C^{*}\to \Cal H$ is the isomorphism of $\Bbb C^{
*}$ with the center of $G$
given by (4.1.7), (4.2.4) or (4.3.4) (for the definition of $G$, see
(1.3.7)).  This induces a natural trivialization of all associated
line bundles over $G(\infty )$.

Now fix $s>0$.  Let $\mu^{(l)}=\mu^{\vert \Cal G\vert^{2s},(l)}$ denote the
unique
$K(l)$-biinvariant measure on $K(l)$ having values in the line
bundle $\vert \Cal G\vert^{2s}$ determined by the normalization condition that
$$\vert\sigma_1\vert^{2s}d\mu^{(l)}=\frac 1Edet\vert A(k)\vert^{2
s}dm(k)\tag 4.6.3$$
is a probability measure, where the notation is that established
in sections 4.1-4.3.  In terms of the trivialization (4.6.1), $\mu^{
(l)}$ is
a multiple of the standard metric on $\Bbb C$ coupled to the Haar
measure for $K(l)$.

It is presently unclear to me how to formulate a general (and
useful) notion of weak convergence for measures having values
in line bundles.  The following formulation is somewhat ad hoc.

\proclaim{(4.6.3)Theorem}The sequence $\{\mu^{\vert \Cal G\vert^{
2s},(l)}\}$ has a unique
weak limit with respect to sections of the form $\phi\theta$, where $
\phi$
is a bounded continuous function of $\Cal G/\Bbb C^{*}$ and $\theta$ is of the
form
$\vert\sigma\vert^{2p}$, where $\sigma$ is one of the fundamental sections and
$p=s-level(\sigma )$.  \endproclaim

\flushpar Proof.  The argument is essentially the same as for
(4.5.3).  We will carry out the proof for $\sigma =\sigma_1$.  It will be
clear from the calculation below that the other cases can be
handled similarly.

We must show that $\mu_s^{(l)}$ has a unique weak limit with respect
to bounded continuous functions on the top stratum
$\Cal N^{-}\cdot \Cal H/\Bbb C^{*}\cdot \Cal N^{+}$.  Weak convergence for the
$
\Cal N^{\pm}$ distributions
follows from (1.1.1).

In considering the diagonal distributions, we use the same
notation as in the proof of (4.5.3). We have
$$(a_{*}\mu_s^{(l)}\hat {)}(\lambda )=\int_{K(l)}a(k)^{-i\lambda}
det\vert A(k)\vert^{2s}dm(k)$$
$$=\int_{K(l)}a(k)^{-i\lambda}\prod\vert\sigma_j(k)\vert^{2s}dm(k
)=\bold c_l(\rho +2s\Lambda_1-i\lambda )\tag 4.6.4$$
(see the proof of (4.4.28)), where $\bold c_l$ denotes the {\bf c}-function for
$G(l)$.  This is given explicitly by a slight modification of the
formulas (4.5.7), where we replace $p+q$ by $p+q+s$ and $r$ by
$r+s$.  The same argument as in the proof of (4.5.3) clearly
applies to show that the $l\to\infty$ limit exists.

We now know that weak limits exist on the top stratum. The
limit is uniquely determined by its diagonal distribution. This
follows directly from (4.4.4), using the fact that line bundles
over groups are trivial.//

\bigskip

\flushpar\S 4.7. Additional Comments.

The argument we gave in section 3.7, to show that there does
not exist an appropriate biinvariant measure on
$\Cal G=SL(2\infty )_{formal}$, depended on special features of the type $
A$
Grassmannian homogeneous spaces. Our comments in section 3.8
carry over directly.

There is a very general method, due to Vershik, which should
be applicable to determining all of the $K(\infty )$-biinvariant
probabilities on $\Cal G$ (of course we contend that there is only one
such probability); see section 3 of [OV].

\bigskip
\centerline{Part III. Loop Groups}

\bigskip

\centerline{\S 0. Introduction for Part III.}

\bigskip

\flushpar\S 0.1. The Kac-Moody group.

Throughout this part of the paper, $G=\dot {G}$ will denote a
connected, simply connected complex Lie group with simple Lie
algebra $\frak g=\dot {\frak g}$.  The Kac-Moody Lie algebra that we will
primarily consider in this part of the paper is the untwisted
affine Lie algebra associated to $\dot {\frak g}$.  We will always use the
standard realization of this algebra as the universal central
extension of the regular loops in $\dot {\frak g}$.  On the occasions when
there is a potential clash of notation involving symbols
attached to $\dot {\frak g}$, and symbols attached to the associated affine
algebra, we will follow the convention in [Kac1] and place a dot
over the symbols attached to $\dot {\frak g}$.  Otherwise we will simply
write $\frak g$ in place of $\dot {\frak g}$, and so on.

The affine Kac-Moody group associated to $\frak g$ is the universal
central extension of the polynomial loop group,
$$0\to \Bbb C^{*}\to\tilde {L}_{pol}G\to L_{pol}G\to 0.$$
(Sometimes the semidirect product of $Rot(S^1)$ with $\tilde {L}_{
pol}G$ is
also referred to as the affine Kac-Moody group associated to
$\frak g$).  We will be interested in the analytical realization of this
extension due to Segal (see chapter 6 of [PS]).  In this approach
one fixes a representation
$$\pi :G\to U(N),$$
and the extension is realized by the determinant line over the
family of Toeplitz operators parameterized by $L_{pol}G$,
$\{A(g)\vert g\in L_{pol}G\}$.  We will recall this realization in \S 1 below.

\bigskip

\flushpar\S 0.2. Constructing measures on the formal completion.

We have already noted that an observation of the Malliavins,
the asymptotic invariance of Wiener measures, implies the
existence of an $L_{pol}K$-invariant probability measure on the
formal completion of the flag space (see \S 3.3 of Part I).  We
will prove a precise version of asymptotic invariance in \S 4.1
below.  In this introduction we will mainly address the
problem of coupling in hermitian structures on line bundles.

In Segal's realization one of the elemental matrix coefficients
for the Kac-Moody group $\tilde {L}_{pol}G$ can formally be written as
$$\tilde {g}\to detA(\tilde {g}).$$
Hence one of the measures that we wish to construct can be
written heuristically as
$$det\vert A(g)\vert^{2s}\Cal D(g),\tag 0.2.1$$
where $\Cal D(g)$ denotes the nonexistent Haar measure for the
unitary real form $LK$, and $s=1$.  Rather than directly attacking
this, we will first consider a two parameter family of
measures which can be written heuristically as
$$d\nu_{\beta ,s}(g)=\frac 1Edet\vert A(g)\vert^{2s}d\nu_{\beta}(
g),\tag 0.2.2$$
where $\nu_{\beta}=\nu^{*=*}_{\beta}$ denotes the Wiener probability measure on
the
continuous loop space, $L_CK$, $\beta >0$ denotes inverse temperature,
and $E$ is a normalization constant.  Analogously the Wiener
measure can be written heuristically as
$$\frac 1Eexp(-\beta \Cal E(g))\Cal D(g),\tag 0.2.3$$
where $\Cal E(g)$ is the energy of a loop,
$$\Cal E(g)=\frac 12\int\langle g^{-1}dg\wedge *g^{-1}dg\rangle ,$$
and $E$ is a normalization constant.

The basic analytical problem is that the Toeplitz determinant
in (0.2.2) is zero on the complement of the Kac-Moody
completion of the loop group, $L_{W^{1/2}}K$, the group of (measure
classes of) loops which belong to the Sobolev space $W^{1/2}$, and
the Wiener measure is supported on this complement.  This
means that the normalization constant and the determinant in
(0.2.2) do not have independent meanings, only the combination
has a meaning, as in (0.2.3).  In \S 3 below we will show that
(0.2.2) can properly be written as
$$\frac 1Zexp(-s\sum_{n>0}n(\vert\hat {g}(n)\vert^2-E\vert\hat {g}
(n)\vert^2))det_2\vert A(g)\vert^{2s}d\nu_{\beta}(g),\tag 0.2.4$$
where $\hat {g}$ denotes the Fourier transform, $E(\cdot )$ denotes
expectation with respect to $\nu_{\beta}$, the exponential and regularized
determinant are well-defined functions relative to Wiener
measure, and $Z$ is now a finite normalization constant.

The next step is to prove existence of the limit
$$d\mu_s(g)=\lim_{\beta\to 0^{+}}d\nu_{\beta ,s}(g),\tag 0.2.5$$
where this is to be understood as a weak limit of probability
measures on the formal completion of the loop space,
$$L_{formal}G=(\tilde {L}_{pol}G)_{formal}/\Bbb C^{*}.\tag 0.2.6$$
We have completely succeeded at this only in the case $s=0$.
The obstacle in the case $s>0$ is the lack of a certain technical
$L^q$ estimate that we explain in \S 4.2.

Thus our main result is the existence of a $L_{pol}K$ biinvariant
probability measure on $L_{formal}G$.  It is very likely that this
measure can be calculated in an explicit manner in several
different ways. Here we have indulged in just one bit of
speculation about the diagonal distribution and Harish-Chandra's
{\bf c-}function in \S 4.4.

\bigskip

\flushpar\S 0.3. Contents.

In \S 2, in addition to describing the formal completion of the
loop group in more concrete terms, we have described another
completion, the hyperfunction loop space, $L_{hyp}G$.  The
advantages of this completion are that (1) the group of analytic
diffeomorphisms of the circle acts naturally, and (2) given a
closed curve on a Riemann surface, there is a natural
projection from $L_{hyp}G$ to the moduli space of holomorphic
$G$-bundles on the surface.  In this paper this completion will
not play much of a role.  However in a projected sequel we
will show that the measures we construct are supported on
$L_{hyp}G$ and are reparameterization invariant, and we hope to
show that they project to canonical measures on the moduli
space.

In \S 3 we construct the measures $\nu_{\beta ,k}$. In \S 4.1 we prove a
precise form of approximate invariance of Wiener measures.
In \S 4.2  and \S 4.3  we prove existence of invariant measures on
the formal flag and loop spaces, respectively. Finally in \S 4.4
we discuss Harish-Chandra's {\bf c-}function.

Although we have not done so explicitly, it should be noted that
our methods apply to twisted affine algebras directly.

\bigskip

\centerline{\S 1. Extensions of Loop Groups.}

\bigskip

\flushpar\S 1.1. The Universal Central Extension

Suppose that $K$ is a simply connected compact Lie group with
simple Lie algebra $\frak k$.  The group of real analytic loops, $
L_{an}K$,
has a distinguished central extension,
$$0\to \Bbb T\to\tilde {L}_{an}K\to L_{an}K\to 0\,\,,\tag 1.1.1$$
which can be characterized in several different ways:  it is
the universal central extension of $L_{an}K$, in the category of Lie
groups; it is the unique simply connected central extension of
$L_{an}K$ by $\Bbb T$; it is the real analytic completion of the unitary
real form of the affine Kac-Moody group associated to $\frak g$.

Segal discovered a relatively explicit model for $\tilde {L}_{an}
K$,
which we will briefly recall (see [PS]).  Fix a nontrivial
representation $\pi :K\to U(N)$. Consider the polarized Hilbert space
$$H=L^2(S^1,\Bbb C^N)=H_{+}\oplus H_{-}\,\,,$$
where $H_{+}$ is the Hardy space of boundary values of
holomorphic functions on the disk, and the induced map from
loops to multiplication operators
$$L_{an}K\to U(H):g\to M_g=\pmatrix A(g)&B(g)\\C(g)&D(g)\endpmatrix \,\,
.\tag 1.1.2$$
Then the operator $A(g):H_{+}\to H_{+}$, the Toeplitz operator
associated to $g$ and $\rho$, is an index zero Fredholm operator,
while $C(g):H_{+}\to H_{-}$, the Hankel operator associated to $g$, is in
$\Cal L_p$ for all $p$.  Define
$$\Cal E=\{(g,q)\in L_{an}K\times U(H_{+})|A(g)q=1\mod \Cal L_1\}\,\,
.$$
Under pairwise multiplication $\Cal E$ is a group, and the
projection $(g,q)\to g$ defines an extension
$$0\to U(H_{+})_1\to \Cal E\to L_{an}K\to 0\,\,.$$
The quotient of this extension by $SU(H_{+})_1$ defines a central
extension by $\Bbb T$ which depends upon $\pi$:
$$0\to \Bbb T\to\tilde {L}_{an}^{\pi}K\to L_{an}K\to 0\,\,.\tag 1.1.3$$
The group $\tilde {L}_{an}K$ is the universal covering of this extension; to
be precise
$$0\to\mu_m\to\tilde {L}_{an}K\to\tilde {L}^{\pi}_{an}K\to 0,\tag 1.1.4$$
where $\mu_m\subseteq \Bbb T$ is the group of $m^{th}$ roots of 1, $
h_{\theta}\in \frak g$ is the
coroot corresponding to the highest root $\theta$, and
$m=\frac 12tr(d\pi (h_{\theta})^{*}d\pi (h_{\theta}))$.

The extensions of the preceding paragraph can be complexified,
using the same method.  One simply replaces $U(H_{+})$ by $GL(H_{
+})$
in the definition of $\Cal E$.  The subgroup of $\tilde {L}_{an}G$ covering the
group of polynomial loops is the affine Kac-Moody group
associated to $\frak g$, in the sense introduced by Kac and Peterson
(see [Kac]).  The complex line bundle associated to the $\Bbb C^{
*}$
bundle
$$\tilde {L}^{\pi}_{an}G\to L_{an}G$$
is $A^{*}Det^{*}$, the pullback of the dual determinant line bundle on
Fredholm operators, relative to the map $g\to A(g)$.  In fact
Segal's construction entails a construction of $Det$ over
Fredholm operators of index zero.

There is a canonical holomorphic function, $\sigma_0^{\pi},$ which is defined
on $\tilde {L}^{\pi}_{an}G$ by
$$\sigma^{\pi}_0([g,q])=\det\,(A(g)q)\,\,.\tag 1.1.5$$
Since $\sigma_0^{\pi}(\tilde {g}\cdot\lambda )=\lambda\sigma_0^{\pi}
(\tilde {g}$) whenever $\lambda\in \Bbb C^{*},$ $\sigma^{\pi}_0$ defines a
holomorphic section of $A^{*}Det;$ this section is the pullback
$A^{*}det$, where $det$, the determinant of a Fredholm operator, is
the canonical section of $Det$.

It is a nonobvious fact (which follows from Borel-Weil theory)
that the pullback of $\sigma^{\pi}_0$ to $\tilde {L}_{an}G$ has a holomorphic $
m^{th}$ root;
we will denote the root with value 1 at the identity by $\sigma_0$.
Restricted to the unitary real form $\tilde {L}_{an}K$, we have
$$|\sigma_0(\tilde {g})|^{2k}=\det|A(g)|^{2s}\,\,,\,\,s=k/m\,\,.\tag 1.1.6$$

\smallskip

\flushpar(1.1.7)Remarks (1) There is a unique automorphic lift
to $\tilde {L}_{an}G$ of the action by $\Cal D^{+}_{an}$ on $L_{a
n}G$ by reparameterization.

(2) For any compact connected Lie group $K$ and representation
$\pi :K\to U(N)$, one can define a central extension $\tilde {L}^{
\pi}_{an}K$ using the
method above, together with some additional arguments to
handle the nonidentity components.  We will make occasional
reference to this extension for the identity component of $L_{an}
\Bbb T$.

(3) The Kac-Moody completion of $\tilde {L}_{pol}K$, as a group, is the
unitary central extension
$$0\to \Bbb T\to\tilde {L}_{W^{1/2}}K\to L_{W^{1/2}}K\to 0,\tag 1.1.8$$
which is constructed in [PS].

\bigskip

\flushpar\S 1.2. Segal Reciprocity.

There is a canonical local cross-section of the extension
$\tilde {L}^{\pi}_{an}G\to L_{an}G$ which is determined by $\sigma^{
\pi}_0$. It is given by
$$g\to [g,A(g)^{-1}],$$
which is well-defined on the set where $A(g)$ is invertible.  This
section induces a vector space splitting
$$\tilde {L}_{an}\frak g\cong L_{an}\frak g\oplus \Bbb Cc,\tag 1.2.1$$
relative to which the bracket is given by
$$[(\xi ,s),(\eta ,t)]=[\xi ,\eta ]+\omega (\xi ,\eta )c\,\,,$$
where
$$\omega (\xi ,\eta )=\frac 1{2\pi i}\,\,\int_{S^1}\,\,\text{Trace}\,\,
(\xi\wedge d\eta )\,\,.$$

Given a compact Riemann surface $\Sigma$ together with an analytic
identification of its boundary $\partial\Sigma\cong S^1$, there are
inclusions
$$L^{\Sigma}_{an}\frak g\to L_{an}\frak g\,\,,\quad and\quad L^{\Sigma}_{
an}G\to L_{an}G\,\,,\tag 1.2.2$$
where $L^{\Sigma}_{an}$ denotes boundary values of holomorphic functions
on $\Sigma$.  The induced extension,
$$\tilde {L}^{\Sigma}_{an}\frak g\to L_{an}^{\Sigma}\frak g,\tag 1.2.3$$
has a natural Lie algebra cross-section given by $\xi\to (\xi ,0)$, by
Cauchy's theorem applied to the cocycle formula above.  Since
$L_{an}^{\Sigma}G$ is simply connected, this induces a splitting of the
corresponding group extension as well.  These splittings are
compatible with sewing of Riemann surfaces, as we will spell
out in Part IV.

\bigskip

\flushpar Remark.  Extensions having this reciprocity property
exist for all connected compact $K$, and they are parameterized
by $H^4(BK;\Bbb Z)$ (see [Pi4]).

\bigskip

\centerline{\S 2. Hyperfunction, Formal and Measurable Loop Spaces.}

\bigskip

\flushpar\S 2.1. The Hyperfunction Loop Space.

Given an open set $\Omega\subset \Bbb C$, let $\Cal O(\Omega )$ denote the
Frechet space of
holomorphic functions on $\Omega$ with the topology of uniform
convergence on compact sets. Given a compact set $S\subseteq \Bbb C$, the
inductive limit
$$\Cal O(S)=\lim_{S\subset\Omega}\Cal O(\Omega )\,\,,\tag 2.1.1$$
with the inductive limit topology, is a reflexive Frechet space.
By definition the space of hyperfunctions on $S$ is the dual
space $\Cal O(S)^{*}$.  We will denote the space of hyperfunctions on $
S^1$
by
$$L_{hyp}\Bbb C=\Cal O(S^1)^{*}\,\,.\tag 2.1.2$$

Consider the one-sided inductive limits
$$\Cal O(S^1_{\pm})=\lim_{r\nearrow 1}\Cal O(\{r<|z|^{\pm 1}<1\})\,\,
.\tag 2.1.3$$
For our purposes the structure of $L_{hyp}\Bbb C$ is most usefully
described by the fact that there is an exact sequence of
Frechet spaces
$$0\to L_{an}\Bbb C\overset {\Delta}\to {\longrightarrow}\Cal O(S^
1_{-})\oplus \Cal O(S^1_{+})\overset {\Phi}\to {\longrightarrow}L_{
hyp}\Bbb C\to 0\tag 2.1.4$$
where $\Delta :f\to (f,f)$ and $\Phi :f,g\to\Phi_{f,g}$,
$$\Phi_{f,g}(u\in \Cal O(S^1))=\underset|z|=R\to\int ufdz-\underset
|z|=r\to\int\,\,ugdz$$
for sufficiently small $R-r$, where $r<1<R$.  This is the
original realization of hyperfunctions due to Sato.

Given an open annulus $A^{\circ}$, there is a group extension
$$0\to \Bbb C^{*}\to\tilde {G}(\Cal O(A^{\circ}))\to G(\Cal O(A^{
\circ}))\to 0,\tag 2.1.5$$
which by restriction injects into the extension $\tilde {L}_{an}G
\to L_{an}G$,
for each positively oriented and analytically parameterized
curve generating the homotopy of $A^{\circ}$.  By taking direct limits
we can form the extensions
$$0\to \Bbb C^{*}\to\tilde {G}(\Cal O(S^1_{\pm}))\to G(\Cal O(S^1_{
\pm}))\to 0.\tag 2.1.6$$
There are natural injections
$$\tilde {L}_{an}G\to\tilde {G}(\Cal O(S^1_{\pm})).$$

The hyperfunction loop space of $G$ is defined by
$$L_{hyp}G=G(\Cal O(S^1_{-}))\times_{L_{an}G}G(\Cal O(S^1_{+}))\tag 2.1.7$$
where $(g_1g,g_2)\sim (g_1,gg_2)$, whenever $g\in L_{an}G$. Similarly the
hyperfunction completion of the affine Kac-Moody group
associated to $G$ is defined by
$$\tilde {L}_{hyp}G=\tilde {G}(\Cal O(S^1_{-}))\times_{\tilde {L}_{
an}G}\tilde {G}(\Cal O(S^1_{+})).\tag 2.1.8$$
Note that there is an inclusion  $\tilde {L}_{an}G\hookrightarrow
\tilde {L}_{hyp}G:\tilde {g}\to [\tilde {g},1]=[1,\tilde {g}]$.

Unless $G$ is abelian neither $L_{hyp}G$ nor $\tilde {L}_{hyp}G$ is a group,
and
these spaces do not have ``unitary real forms'' (i.e.  the
involution that defines $LK\hookrightarrow LG$ does not extend).  I do not know
whether there is some alternate construction which suggests a
formulation of nonabelian hyperfunctions for spaces other than
the circle.

\bigskip

\flushpar\S 2.2. The Birkhoff Decomposition and Complex Structure.

The Birkhoff decomposition for $L_{an}G$ induces a decomposition for
$L_{hyp}G$, by taking direct limits:
$$L_{an}G=\bigsqcup G(\Cal O(\bar {D}^{-}))\cdot\lambda\cdot G(\Cal O
(\bar {D}^{+}),\tag 2.2.1$$
implying
$$G(\Cal O(S^1_{-}))=\bigsqcup G(\Cal O(D^{-}))\cdot\lambda\cdot
G(\Cal O(\bar {D}^{+})),\tag 2.2.2$$
$$G(\Cal O(S_{+}^1))=\bigsqcup G(\Cal O(\bar {D}^{-}))\cdot\lambda
\cdot G(\Cal O(D^{+})),\tag 2.2.3$$
and finally
$$L_{hyp}G=\bigsqcup G(\Cal O(D^{-}))\cdot\lambda\cdot G(\Cal O(D^{
+})),\tag 2.2.4$$
where the disjoint unions are over $\check {T}$, the group of
homomorphisms from $S^1$ into a maximal torus $T$ of $K$, and $D^{
\pm}$
are the open disks with boundary $S^1$ centered at 0 and $\infty$,
respectively.

The stratum corresponding to $\lambda =1$ is open, dense and
holomorphically equivalent to
$$G(\Cal O(D^{-}))_1\times G\times G(\Cal O(D^{+}))_1.\tag 2.2.5$$
Here $G(\Cal O(D^{+}))_1=\{g\in G(\Cal O(D^{+}))|g(0)=1\}$; this is a complex
profinite nilpotent Lie group which, as a complex manifold, is
holomorphically equivalent to its Lie algebra:
$$G(\Cal O(D^{+}))_1\cong \Cal O(\Omega^{(1,0)}(D^{+})\otimes \frak g
)\tilde {=}\frak g(\Cal O(D^{+})_0),\tag 2.2.6$$
where $g\leftrightarrow\omega =(\partial g)g^{-1}=\partial x\leftrightarrow
x$.  For $D^{-}$ a similar
description applies with $\infty$ in place of 0. $_{}$Note that the
exponential map
$$exp:\frak g(\Cal O(D^{+})_0)\to G(\Cal O(D^{+}))_1$$
is injective, but it is not surjective, as observed by Goodman
and Wallach (see Remark (ii) following (8.4.5) of [PS]).

By taking direct limits
$$L_{hyp}G=\bigcup\lambda\cdot G(\Cal O(D^{-}))_1\cdot G\cdot G(\Cal O
(D^{+}))_1\,\,,\tag 2.2.7$$
a union of dense open sets.  Similarly, by applying Segal reciprocity
to $G(\Cal O(D^{\pm}))$,
$$\tilde {L}_{hyp}G=\bigcup\tilde{\lambda}\cdot G(\Cal O(D^{-}))_
1\cdot G\cdot \Bbb C^{*}\cdot G(\Cal O(D^{+}))_1\,\,,\tag 2.2.8$$
where $\tilde{\lambda}$ denotes a choice of element in $\tilde {L}_{
pol}G$ which covers
$\lambda\in\,\,\text{Hom}\,\,(S^1,T)$.

\smallskip

\flushpar(2.2.9) Remarks.  (a) The group of
orientation-preserving real analytic homeomorphisms of the
circle, $\Cal D^{+}_{an}$, acts holomorphically on $L_{hyp}G$ and $
\tilde {L}_{hyp}G$, naturally
extending its actions on the analytic loop space and its
extension.  Note that there does not exist an induced action of
$\Cal D^{+}_{an}$ on the corresponding flag space, because $G(\Cal O
(D^{+}))$ is not
invariant with respect to $\Cal D^{+}_{an}$.

(b) Given a Riemann surface $\hat{\Sigma }=\Sigma^{-}\circ\Sigma^{
+}$ as in \S 0.4, then
$$G(\Cal O(\Sigma^{-0}))\backslash L_{hyp}G/G(\Cal O(\Sigma^{+0})
)\cong H^1(\hat{\Sigma };\Cal O_G)\cong G(\Cal O(\Sigma^{-}))\backslash
L_{an}G/G(\Cal O(\Sigma^{+}))\,\,,$$
where $H^1(\hat{\Sigma },\Cal O_G)$ denotes the set of equivalence classes of
holomorphic $G$-bundles. The point of this, from the viewpoint
of this paper, is that this space is canonically symplectic; we
believe that the biinvariant probability measure that we
construct in \S 4 projects to the normalized symplectic volume
element.

\smallskip

\bigskip

\flushpar\S 2.3. The formal loop space.

In part I we constructed a formal completion of the Kac-Moody
group $\tilde {L}_{pol}G$, which we denoted by $\Cal G$.  We will refer to the
quotient $\Cal G/C_0$ as the formal loop space of $G$, and we will
denote it by $L_{formal}G$.  Our purpose here is to briefly note
that there is a more direct way to define the formal loop
space.

Let $\Bbb C[[z]]$ denote the ring of formal series in the indeterminant
$z$, let $\Bbb C((z))$ denote the algebra of meromorphic formal series
(that is, $\Bbb C((z))$ consists of formal expressions of the form
$$\sum_{-N}^{\infty}c_nz^n,\tag 2.3.1$$
where $c_n\in \Bbb C$, and one inverts the series by factoring out the
lowest power of $z$ and using the geometric series), and let
$\Bbb C(z)$ denote the subring of finite meromorphic series (that is,
finite Laurent series).  Note that
$$L_{pol}G=G(\Bbb C(z)).$$

Since $\Bbb C((z))$ is an algebra, it is possible to define $G(\Bbb C
((z)))$ for
any algebraic group $G$.  For simply connected $G$, making use of
the fact that inversion in $G$ is a polynomial, it is possible to
realize $G(\Bbb C((z)))$ very concretely in each case (for example
$SL(n,\Bbb C((z)))$ consists of $n\times n$ matrices $g$ with coefficients in
$\Bbb C((z))$ which satisfy the equation
$$g\cdot adj(g)=1,$$
where $adj(\cdot )$ denotes the adjugate; for $SO(n,\Bbb C((z)))$ we replace
$adj(\cdot )$ by transpose, and so on; for $GL_n$, which has a more
subtle realization, see \S 1 of [BL]).

We claim that the formal loop space of $G$ is given by
$$L_{formal}G=G(\Bbb C((z^{-1})))\times_{G(\Bbb C(z))}G(\Bbb C((z
))).\tag 2.3.2$$
Comparing (2.3.2) with the definition of the formal completion
in Part I, we see that the essence of our claim is that the
natural inclusion
$$\Cal G^{+}/C_0\equiv G(\Bbb C(z))G(\Bbb C[[z]])\hookrightarrow
G(\Bbb C((z)))\tag 2.3.3$$
is surjective, or equivalently, that the induced inclusion of flag
spaces
$$G(\Bbb C(z))/G(\Bbb C[z])\hookrightarrow G(\Bbb C((z)))/G(\Bbb C
[[z]])\tag 2.3.4$$
is surjective.  This follows from the existence of a Birkhoff
decomposition for the right hand side of (2.3.3) (and hence also
(2.3.4)).  From one perspective this is a consequence of the
general ``gaga'' principle of algebraic geometry (see \S 1 of [BL]).
One can also mimick the method in chapter 8 of [PS] (one
considers the polarized vector space
$$V=V_{-}\oplus V_{+}=\frak g(\Bbb C[z^{-1}])_0\oplus \frak g(\Bbb C
[[z]]),$$
the Grassmannian model, and so on). This has the advantage of
yielding a concrete method for producing the unique factorization
$$G(\Bbb C((z)))=K(\Bbb C(z))\times A\times \Cal N^{+}.$$

The Kac-Moody extensions
$$0\to \Bbb C^{*}\to \Cal G^{\pm}\to G(\Bbb C((z^{\pm})))\to 0\tag 2.3.5$$
were constructed from the Kac-Peterson point of view in \S 1 of
Part I.  One can also construct these extensions by modifying
Segal's construction, outlined in \S 1.1.  In the definition of the
extension $\Cal E$ one considers $q$'s which are asymptotically block
upper (respectively, lower) triangular, where the blocks are
determined by the representation $\pi$.

\bigskip

\flushpar\S 2.4. Relations with functional loop spaces.

In subsequent sections we will construct unitarily invariant
measures on the space $L_{hyp}G$. The question naturally arises
whether these measures are supported on a smaller space having
a unitary form (see \S 5). We do not have significant results in this
direction, but we would like to comment on some issues of
intrinsic interest that arise in trying to construct a support
with reasonable geometric properties.

Relative to the rotation invariant measure class on $S^1$, let
$L_{meas}G$ denote the Polish (i.e.  complete separable metric)
topological group obtained by equipping the group of
equivalence classes of measurable maps $S^1\to G$ with the
topology of convergence in measure.

Define the set of all possible Riemann-Hilbert factorizations,
$L_{RH}G$, as the subset of $L_{hyp}G$ satisfying $[g,h]\in L_{RH}
G$ if and
only if $g(r^{-1}e^{i\theta})$ and $h(re^{i\theta})$ converge nontangentially
as $
r\uparrow 1$, for
a.e.  $\theta$.  We define $L_{RH}K$ in the same way, but require the
boundary values to be in $K$, for a.e.  $\theta$.

This can be expressed more directly as follows. Let
$$\aligned
\Cal O(\{r_0<\vert z\vert <1\})_{meas}&=\{f\in \Cal O(\{r_0<\vert
z\vert <1\}:f(re^{i\theta})\\&converges\quad nontangentially\quad
as\quad r\uparrow 1,fora.e.\theta \}.\endaligned
$$
The natural boundary map,
$$\partial :\Cal O(\{r<|z|<1\})_{meas}\to L_{meas}\Bbb C,\tag 2.4.1$$
is an inclusion (see 1.9, chapter XIV of [Zygmund]). Define
$$\Cal O(S^1_{\pm})_{meas}=\lim_{r\uparrow 1}\Cal O(\{r<\vert z\vert^{
\pm 1}<1\})_{meas}.\tag 2.4.2$$
There are inclusions
$$G(\Cal O(S^1_{\pm}))\hookrightarrow L_{meas}G.$$
Then essentially by definition
$$L_{RH}G=G(\Cal O(S^1_{-})_{meas})\times_{L_{an}G}G(\Cal O(S^1_{
+})_{meas}).\tag 2.4.3$$

The stratification of $L_{hyp}G$ restricts to a stratification of
$L_{RH}G$, where $G(\Cal O(D^{\pm}))$ is replaced by $G(\Cal O(D^{
\pm})_{meas})$.
Unfortunately, whereas we know that $G(\Cal O(D))_1$ is
holomorphically equivalent to a linear Frechet space, at this
point we can only say that $G(\Cal O(D)_{meas})_1$ is topologically
contractible. It seems unlikely that $L_{RH}G$ is a complex
manifold.  Also, the natural boundary map,
$$L_{RH}G\to L_{meas}G,\tag 2.4.4$$
is not injective, and it is not clear at all how to characterize
the image.  The extent to which it is not injective measures
lack of uniqueness for the Riemann-Hilbert factorization of
loops (It would clearly be desirable to find a ``maximal domain''
for Riemann-Hilbert factorization; see (2.4.7) below for what I
believe is the strongest known result).

These considerations explain part of the following diagram:
$$\matrix \quad&&L_CG&&\\\quad&\swarrow&\downarrow&\quad\quad\\&&\\
\prod_{S^1}G&&L_{QC}G&\quad&\quad\\&&\downarrow&\searrow&\quad\\&
\quad&L_{meas}G&\leftarrow&L_{RH}G&\rightarrow&L_{hyp}G\endmatrix \tag 2.4.5$$
Only the map $L_{RH}G\to L_{meas}G$ is noninjective.

The inclusion of continuous loops into the hyperfunction loop
space follows from the basic theorem about Riemann-Hilbert
factorization for continuous loops (for one approach see [CG]).
More generally, one can consider Douglas's algebra of
quasi-continuous functions,
$$QC=(C+H^{\infty})\cap (C+H^{\infty})^{conj}=VMO\cap L^{\infty}.\tag 2.4.6$$
(see Chapter 6 of [Sarason]). The basic existence and uniqueness
result is the following

\proclaim{(2.4.7) Proposition } Given $g\,\,\in\,\,L_{QC}G$, there exist
$g_{\pm}\,\,\in\,\,G(\Cal O(D^{\pm}))$ and $\lambda\,\,\in\,\,\check {
T}$, such that

(i)  $g_{\pm}(re^{i\theta})\to g_{\pm}(e^{i\theta})$ in $L^p$ as $
r\uparrow 1$, for all $1<p<\infty$.

(ii)  $g=g_{-}\cdot\lambda\cdot g_{+}\,\,\text{a.e. $\theta$.}$

(iii) $\lambda$ is determined up to conjugation by $W$, the Weyl group,

and

(iv)  if $\lambda =1$ and $g_{-}(\infty )=1$, then $g_{\pm}$ are uniquely
determined.
\endproclaim

The Grassmannian approach of [PS] can be modified to give a proof of
this more general result.  This is so because of the following
observations:

(1) the ideal of Hilbert-Schmidt operators can be replaced by compact
operators;

(2) the quasi-continuous loops are precisely the loops which belong to
the restricted general linear group corresponding to compact
operators;

(3) $L^2(S^1)$ can be replaced by $L^p(S^1)$ for any $1<p<\infty$.  (the
point is that the projection onto the Hardy subspace is
continuous for this range of $p$).

\bigskip

\flushpar Remark (2.4.8).  Segal's method gives a model for the
$C^{*}$ bundle $\tilde {L}_{QC}G$ induced by the inclusion $L_{QC}
\to L_{hyp}G$, namely
$\tilde {L}_{QC}G$ is the simply connected $m$-fold covering of $
\tilde {L}^{\rho}_{QC}G$, where
this bundle is the $\Bbb C^{*}$ bundle associated to the complex line
bundle $A^{*}Det^{*}$.  This bundle is not a group, and the $C^{*}$-bundle
$\tilde {L}_CK\to L_CK$ does not have a $\Bbb T$-reduction.

\bigskip
\centerline{\S 3. Existence of the Measures $\nu_{\beta ,k}$, $\beta
>0.$}

\bigskip

Fix a representation $\pi :K\to U(N).$ Given $\beta >0,$ $\beta\cdot
trace$ induces
an Ad-invariant inner product on $\frak k\subseteq u(N,\Bbb C)$, hence a
bi-invariant Riemannian metric on $K$.  Let $\nu_{\beta}^{*=*}$ denote the
corresponding Wiener probability measure on the continuous
loop space $L_CK=Path^{*=*}K$.  In this section we will establish
the existence of a probability measure on $L_CK$, which on a
formal level can be written as
$$d\nu_{\beta ,k}=\frac 1E\vert\sigma_0\vert^{2k}d\nu_{\beta}=\frac
1E\det|A(g)|^{2s}\,\,d\nu^{*=*}_{\beta}(g)\,\,.\tag 3.0.1$$
These formal expressions indicate the quasi-invariance
properties which characterize the measure.  The number $s\ge 0$
is related to $k$ as in (1.1.6).

\bigskip

\flushpar\S 3.1. Abelian vs. Nonabelian Cases.

To motivate the regularization which arises in constructing the
measure, first consider the case $K=\Bbb T$, together with the
defining representation.  In this case
$$(L_C\Bbb T)_0\cong \Bbb T\times\left\{x\,\,\in\,\,L_C\Bbb R|\int
x=0\right\},\tag 3.1.1$$
where $g=\lambda\cdot\exp(ix)$.  In these coordinates
$$d\nu_{\beta}(g)=(2\pi )^{-1}d\theta\times\overset {}\to {\prod_{
n>0}}(2\pi )^{-1}\beta n^2\exp(-\beta n^2|x_n|^2)dLeb(x_n)\,\,,\tag 3.1.2$$
where $\{x_n\}$ are the Fourier coefficients of $x$.  Szego's formula
for Toeplitz determinants asserts that in this abelian context
$$\det|A(g)|^{2k}=\exp(-k\sum^{\infty}_1n|x_n|^2)\,\,.\tag 3.1.3$$
Thus the desired measure is
$$d\nu_{\beta ,k}(g)=(2\pi )^{-1}d\theta\times\prod_{n>0}(2\pi )^{
-1}(\beta n^2+kn)\exp(-(\beta n^2+kn)|x_n|^2)dLeb(x_n)$$
$$=(2\pi )^{-1}d\theta\times\frac 1{\Cal Z}exp(-k\sum_{n>0}n(\vert
x_n\vert^2-E\vert x_n\vert^2))d\nu_{\beta}(x),\tag 3.1.4$$
where $E|x_n|^2=(\beta n^2)^{-1}$ is the expected value relative to $
\nu_{\beta}$, and
$\Cal Z$ is the finite number
$$\Cal Z=\prod_{n>0}\{(1+\frac k{\beta n})^{-1}exp(\frac k{\beta
n})\}=\frac {\Gamma (1+k/\beta )}{e^{\gamma k/\beta}}\tag 3.1.5$$
where $\gamma$ is Euler's constant.

The important point here is that the series $\sum n(|x_n|^2-E|x_n
|^2)$
represents a well-defined random variable with respect to $\nu_{\beta}$.
It is trivial to check that the series converges absolutely in
the $L^2$ sense; more subtle is the fact that the series converges
conditionally for a.e.$ $ $x$ $[\nu_{\beta}]$.  This follows from the following
general result (see theorem 22.6 of [Bill2]),

\proclaim{Kolmogorov's Theorem}If the random variables $X_n$ are
independent, each has mean zero, and $\sum E|X_n|^2$ is finite, then
$\sum X_n$ converges conditionally with probability one.  \endproclaim

In our case $X_n=n(\vert x_n\vert^2-E\vert x_n\vert^2)$ and $\sum
E\vert X_n\vert^2=const\sum n^{-2}$.  It
can be checked that in our case, the series $\sum X_n$ converges
absolutely with probability zero.

In the nonabelian case there is not an explicit Szego type
formula for $\det|A|^2$.  Instead we will consider the expression

$$\aligned
det\vert A(g)\vert^{2s}&=det(1-\vert C(g)\vert^2)^s\\&=exp(-s\cdot
tr\vert C(g)\vert^2)det((1-\vert C(g)\vert^2)e^{\vert C(g)\vert^2}
)^s\\&=exp(-s\sum n\vert\hat {g}(n)\vert^2)det_2\vert A(g)\vert^{
2s}\endaligned
\tag 3.1.6$$
The measure $\nu_{\beta ,k}$ will be given by
$$\frac 1{\Cal E}\exp(-s\sum^{\infty}_1n(|\hat {g}(n)|^2-E|\hat {
g}(n)|^2))det_2|A(g)|^{2s}d\nu^{*=*}_{\beta}(g)\,\,.\tag 3.1.7$$
where the exponential and regularized determinant are well-defined
random variables relative to $\nu^{*=*}_{\beta}$, and $\Cal E$ is a finite
number.

\bigskip

\flushpar\S 3.2. Statement of the Main Result.

On a formal level the Wiener measure $\nu^{*=*}_{\beta}$ can be written as
$$d\nu^{*=*}_{\beta}(g)=\frac 1E\exp(-\beta \Cal E(g))\prod_{\theta}
dg(\theta )\,.\tag 3.2.1$$
where $\Cal E(g)$ is the energy of the loop $g$.  As this expression
suggests, $\nu^{*=*}_{\beta}$ is quasi-invariant under both left and right
translations by group elements from $L_{W^1}K$, the group of loops
of finite energy.  This is proven by first establishing the
analogous fact for the Wiener measure $\nu_{\beta}^{1,*}$ on the path space
$Path^{1,*}K$, by using Ito's map to reduce to a linear problem, then
localizing to loops (see [MM1] and the references there; we will
give a proof of the localization step in (4.1.12) below).

It would be highly desirable to formulate a characterization of
Wiener measure in terms of its invariance properties.  But
because it is not known how to directly write down the
Radon-Nikodym cocycle (this involves Ito's map, and the
definition of the map depends upon the background measure), it
is not so clear how to do this.  We will content ourselves
with characterizing the measures $\nu_{\beta ,k}$ as measures which are
absolutely continuous with respect to $\nu^{*=*}_{\beta}$.  This at least has
the virtue of showing that all means of regularizing the formal
expression (3.0.1) give the same result.

The main result is the following

\proclaim{(3.2.2) Theorem } (a) There exists, up to a multiple, a
unique measure $\mu <<\nu^{*=*}_{\beta}$ satisfying
$$\frac {d\mu (g^{-1}_Lgg_R)}{d\mu (g)}=\left|\frac {[(\tilde {g}_
L,\tilde {g}_R)\cdot\sigma_0](g)}{\sigma_0(g)}\right|^{2k}\frac {
d\nu^{*=*}_{\beta}(g^{-1}_Lgg_R)}{d\nu^{*=*}_{\beta}(g)}$$
for all $\tilde {g}_L,\tilde {g}_R\in\tilde {L}_{W^1}K$.

(b) The measure $\mu$ is finite.\endproclaim

\proclaim{(3.2.3)Remark} To make sense of this statement, one must
know that $\{\sigma_0=0\}$ has $\nu^{*=*}_{\beta}$ measure zero.  This follows
from
(3.1.1) of Part I.  \endproclaim

The bulk of the proof will be given in subsequent subsections.
The expression (3.1.7) for the measure will be justified in \S 3.3,
invariance will be established in \S 3.4, and finiteness will be
proven in \S 3.5.  Here we will give the

\bigskip

\flushpar Proof of uniqueness.  The measure $\mu$ is in fact
characterized by its quasi-invariance properties with respect
to either left or right translations (together with the fact that
we always assume it is absolutely continuous with respect to
$\nu_{\beta}^{*=*}$).  For definiteness suppose that $\mu$ is quasi-invariant
with
respect to the left translation action
$$L_{an}K\times L_CK\to L_CK:g_L,g\to g_L\cdot g,$$
and suppose that the Radon-Nikodym cocycle is as in the
proposition.  Gross has proven that $\nu_{\beta}^{*=*}$ is ergodic with
respect to left translation (and similarly for the right action); see
[Gr].  This implies that the measure $\mu$ must actually be
equivalent $ $to $\nu^{*=*}_{\beta}$.  If $\rho d\mu$ is another such measure,
then
ergodicity implies that $\rho$ = constant a.e.  $//$

\bigskip

\flushpar\S 3.3. First part of the proof of (3.2.2): justifying the
expression for the measure.

To justify the expression (3.1.7) for $\mu$, we need to recall several
facts:

\smallskip

(3.3.1) The function $\det_2:\,1+\Cal L_2\to \Bbb C:\,\,\,1+A\to\det
((1+A)e^{-A})$, is well
defined and continuous.

(3.3.2) Let $\rho$ denote the function
$$\rho (x)=(|x|\log|x|^{-1})^{\frac 12}\,.$$
The space of $\rho -$Lipschitz continuous functions on the interval
$[0,2\pi ]$, denoted $Lip(\rho )$, is defined by the norm
$$|f|_{Lip(\rho )}=|f|_{\infty}+\sup_{s<t}\frac {|f(s)-f(t)\vert}{
\rho (s-t)}.$$
The space $Lip(\rho )$ is a Banach algebra. The Banach Lie
group $L_{Lip(\rho )}K$ has full measure relative to $\nu^{*=*}_{
\beta}$, by
Levy's modulus of continuity for Brownian motion (see [McKean]).

(3.3.3) Write $M_g=\pmatrix A(g)&B(g)\\C(g)&D(g)\endpmatrix $ as in (1.1.2).
Then $\pmatrix
\quad&B\\C&\quad\endpmatrix $ is in
the Schatten class $\Cal L_p$ if and only if $g\in B_p^{1/p}$, the $
p$-Besov
space defined by the norm
$$|g|^p_{B_p^{1/p}}=\int\Delta^{-2}\int |g(t+\Delta )-g(t)|^pdt\,
d\Delta .$$
(see [Peller]).

(3.3.4) $Lip(\rho )\subset B_p^{1/p}$, for all $p>2$. This is easy to check.

(3.3.5) For $g$ as in (3.3.3), $\left(\matrix \quad&B\\C&\endmatrix \right
)$ has kernel
$$K_g=\frac 1{2\pi i}\frac {g(\zeta )-g(z)}{\zeta -z}d\zeta .$$
Thus we can write the Hilbert-Schmidt norm of C as
$$tr|C|^2=\sum_{n>0}n|\hat {g}(n)|^2=\frac 1{8\pi^2}\iint\frac {|
g(t+\Delta )-g(t)|^2}{|e^{i\Delta}-1|^2}dt\,d\Delta$$

Let $E$ denote the expectation functional relative to the measure
$\nu_{\beta}$.  For $g\in L_CG$ consider the integral
$$\aligned
\frac 1{8\pi^2}\iint_{\vert e^{i\Delta}-1\vert >\delta}&\frac {|g
(t+\Delta )-g(t)|^2-E|g(t+\Delta )-g(t)|^2}{|e^{i\Delta}-1|^2}dt\,
d\Delta\\&=\sum_n\left(\frac 1{2\pi}\int_{|e^{i\Delta}-1|>\delta}\frac {
|e^{in\Delta}-1|^2}{|e^{i\Delta}-1|^2}d\Delta\right)(|\hat {g}(n)
|^2-E|\hat {g}(n)|^2)\endaligned
\tag 3.3.6$$
where $t$ and $\Delta$ vary from $0$ to $2\pi$, modulo the restriction noted.
The expression (3.1.6) for the measure $\mu$ is justified by the
following

\proclaim{(3.3.7) Proposition }  The integral in (3.3.6) has a limit
as  $\delta\to 0$,  the sum
$$\sum_{n>0}n(\vert\hat {g}(n)\vert^2-E\vert\hat {g}(n)\vert^2)$$
converges conditionally, and these two limits are equal, for a.e.
g $[\nu_{\beta}^{*=*}]$.  \endproclaim

\proclaim{(3.3.8) Lemma} Let $I_n(\delta )=(2\pi n)^{-1}\int_{\delta}^{
2\pi -\delta}\frac {|e^{in\theta}-1|^2}{|e^{i\theta}-1|^2}d\theta$, and
suppose $\{a_n\}$, is a sequence of complex numbers.  If $\sum a_
n$
converges conditionally, then $\lim_{\delta\to 0}\sum a_nI_n(\delta
)$ exists and equals
$\sum a_n.$
\endproclaim

\flushpar Proof of (3.3.8).

A straightforward calculation shows that
$$I_n-I_{n+1}=\frac 2{n(n+1)}\frac {cos(\delta /2)-cos(\delta /2+
n\delta )}{sin(\delta /2)}\tag 3.3.9$$
Note that this quantity is infrequently negative, for small $\delta$.
In fact if $n_l$ denotes the $l$th integer for which (3.3.9) is
negative, then
$$n_l\ge const\cdot\frac l{\delta}\tag 3.3.10$$

Using this it is easy to show that there is a uniform bound for
$\sum |I_n(\delta )-I_{n+1}(\delta )|\,,$ for $0\le\delta\le\pi /
2$. For
$$\aligned
\sum_{n>0}\biggl|I_n(\delta )-I_{n+1}(\delta )\biggr|&=\sum_{n>0}
(I_n(\delta )-I_{n+1}(\delta ))+2\sum_{*}\,(I_{n+1}(\delta )-I_n(
\delta ))\\&=I_1(\delta )+2\sum_{*}\,(I_{n+1}(\delta )-I_n(\delta
)),\endaligned
$$
where in the second sum the $``*$'' indicates that we sum only
over those $n$ for which $\cos(\delta /2+n\delta )\ge\cos(\delta
/2)$.  We always have
$$\aligned
\biggl|I_n(\delta )-I_{n+1}(\delta )\biggr|\le const\cdot\frac 1{
n(n+1)}\cdot\frac 1{\delta}\,.\endaligned
\tag 3.3.11$$
Together with (3.3.10) this implies that
$$\sum_{*}(I_{n+1}(\delta )-I_n(\delta ))\le const\sum_{l>0}\frac
1{(l/\delta )(l/\delta +1)}\cdot\frac 1{\delta}.$$
Thus this sum tends to zero as $\delta\to 0$, which establishes the uniform
bound.

For the proof of (i), we can suppose $s_N=\sum_1^Na_n\to 0$ as $N
\to\infty\,.$ Then
$$\sum I_n(\delta )a_n=\sum (I_n(\delta )-I_{n+1}(\delta ))s_n$$
$$\Rightarrow\biggl|\sum I_n(\delta )a_n\biggl|\le\biggr|\sum^M(I_
n(\delta )-I_{n+1}(\delta ))s_n\biggr|+(\sup_{n\ge M}|s_n|)\sum^{
\infty}_M|I_n(\delta )-I_{n+1}(\delta )|$$
The uniform bound of the second paragraph, the fact $s_n\to 0\,,$
and the fact $I_n(\delta )\to 1$ as $\delta\to 0$ now imply that
$\lim_{\delta\to 0}\sum\,I_n(\delta )a_n=0$.  //

\bigskip

\flushpar Proof of Proposition (3.3.7).

The strategy of the proof is standard, and it involves three
steps:  (1) verify a similar result for the linear path space
$Path^{0,*}\frak k$; (2) use Ito's map to verify the result for $
Path^{1,*}K$;
and (3) localize the result to based loops.  Throughout this
proof $t,\Delta$ and $\tau$ will be constrained by $0\le t,\tau ,
\Delta\le 2\pi\,$; additional
constraints will be explicitly noted.

\bigskip

\flushpar Step (1).

Let $\nu^{0,*}_{\beta}$ denote the Wiener measure with inverse temperature
$\beta$ on the continuous pathspace $Path^{0,*}\frak k$.  In this first step we
will establish

\proclaim{(3.3.12) Lemma}
$$\lim_{\delta\to 0}\underset{|\tau -t|>\delta}\to\iint\left(\frac {
|w(\tau )-w(t)|^2-E|w(\tau )-w(t)|^2}{|\tau -t|^2}\right)d\tau\,d
t$$
exists for a.e. $w\quad [\nu^{0,*}_{\beta}]$.\endproclaim

\flushpar Proof of (3.3.12).

Recall that for the probability measure
$$\prod_{n>0}(2\pi )^{-1}n^2\beta exp(-n^2\beta\vert x_n\vert^2)d
Leb(x_n)\tag 3.3.13$$
on the space $\{x\in L_C\frak k\vert\,\,\int x=0\}$, Kolmogorov's theorem on
random
series implies that $\sum n(|x_n|^2-E|x_n|^2)$ converges conditionally
with probability one.  By (i) of (3.3.8),
$$\lim_{\delta\to 0}\frac 1{8\pi^2}\iint_{\vert e^{i\Delta}-1\vert
>\delta}\frac {|x(t+\Delta )-x(t)|^2-E|x(t+\Delta )-x(t)|^2}{|e^{
i\Delta}-1|^2}dt\,d\Delta\,\tag 3.3.14$$
$$=\lim_{\delta\to 0}\sum_{n>0}I_n(\delta )n(|x_n|^2-E|x_n|^2)=\sum_{
n>0}n(|x_n|^2-E|x_n|^2)\,\tag 3.3.15$$
with probability one.

We claim that this in turn implies that
$$\lim_{\delta\to 0}\frac 1{8\pi^2}\iint_{\vert\tau -t\vert >\delta}\frac {
|x(\tau )-x(t)|^2-E|x(t)-x(t)|^2}{\vert\tau -t\vert^2}dt\,d\tau\tag 3.3.16$$
exists for a.e.  $x$.  The essential difference between this claim
and what we established in the preceding paragraph is that
because of the exponential in (3.3.14), we also cut out the
region where the difference $\Delta =\tau -t$ was close to $2\pi\,
.$
However, it is easy to show that this region can be ignored.
In fact for $\delta$ small
$$\int^{\delta}_0\left(\int^{2\pi}_{2\pi -\delta +\tau}\frac {|x(
\tau )-x(t)|^2}{|e^{i\tau}-e^{it}|^2}dt\right)d\tau <\infty\tag 3.3.17$$
for a.e.  $x$.  For by writing $t=2\pi +\tau -\Delta$ and using the
periodicity of the integrand, this integral equals
$$\int^{\delta}_0\int^{\delta}_{\tau}\frac {|x(\tau +\Delta )-x(\tau
)|^2}{|e^{i\Delta}-1|^2}d\Delta\,d\tau .\tag 3.3.18$$
By Levy's modulus of continuity this is bounded by
$$\text{(constant)}\ \int^{\delta}_0\left(\int^{\delta}_{\tau}\frac {
\Delta\log\Delta^{-1}}{\Delta^2}d\Delta\right)d\tau ,\tag 3.3.19$$
which is finite.  The same argument shows that the expected
value of this double integral is finite.  This establishes our
modified claim, i.e.  the existence of the limit (3.3.16).

Let $Path^{0,0}\frak k$ denote the based path space, and $\nu^{0,
0}_{\beta}$ the Brownian
Bridge.  The map
$$Path^{0,0}\frak k\to \{x\in L_C\frak k\vert\quad\int x=0\}:w\to
x\tag 3.3.20$$
where $x(t)=w(t)-\frac 1{2\pi}\int^{2\pi}_0w(\tau )d\tau$, maps $
\nu^{0,0}_{\beta}$ to the measure
in (3.3.13).  It follows that the limit in the statement of the
Lemma, (3.3.12), exists for a.e.  $w$ relative to the Brownian
Bridge.

Let $Path^{0,\chi}\frak k=\{w_1\in Path^{0,*}\frak k|w_1(2\pi )=\chi
\in \frak k\}\,$, and $\nu^{0,\chi}_{\beta}$ the
corresponding Brownian Bridge.  The map
$$Path^{0,0}\frak k\to Path^{0,\chi}\frak k:w\to w_1\,,\tag 3.3.21$$
where $w_1(t)=w(t)+\frac t{2\pi}\chi$, sends $\nu^{0,0}_{\beta}$ to a measure
equivalent to
$\nu^{0,\chi}_{\beta}\,.$ Also
$$\aligned
&\iint\biggl|\frac {|w(\tau )-w(t)|^2-|w(t)-w(\tau )-\frac {\tau
-t}{2\pi}\chi |^2}{|\tau -t|^2}\biggr|d\tau dt\\\le\frac 1{2\pi}&
\iint\frac {|(\tau -t)\chi |^2}{|\tau -t|^2}+\frac 2{(2\pi )^{\frac
12}}\iint\frac {|w(\tau )-w(t)|\,\,|(\tau -t)\chi |}{|\tau -t|^2}
<\infty\endaligned
\tag 3.3.22$$

for a.e. $w$ (Levy's modulus of continuity implies that the second
integral is finite a.e. $w$), and similarly for the expected values.
This shows that the limit in (3.3.12) exists for a.e. $w$ relative to each
Brownian Bridge.  This completes the proof of (3.3.12) and step
(1).//

\smallskip

\flushpar Step (2).

Ito's map
$$I:Path^{0,*}\frak k\to Path^{1,*}K:w\to g\tag 3.3.23$$
is defined by the stochastic integral equation
$$g(t)=1+\int^t_0g(\tau )\circ dw(\tau )\,,\tag 3.3.24$$
interpreted in the sense of Fisk-Stratonowich.  The map $I$ is a
measureable map (defined for a.e.  $w\quad [\nu_{\beta}^{o,*}])$ which pushes $
\nu_{\beta}^{0,*}$
forward to the Wiener measure $\nu_{\beta}^{1,*}$ on $Path^{1,*}K$ (see Chapter
III
of [IW], [MM1] or [McKean]).

In this step we will prove

\proclaim{(3.3.25) Lemma }
$$E\iint\biggl|\frac {|g(\tau )-g(t)|^2-|w(\tau )-w(t)|^2}{|\tau
-t|^2}\biggr|dt\,d\tau <\infty\,.$$
\endproclaim

This Lemma implies that a.e. $w\quad [\nu^{0,*}_{\beta}]\,,$
$$\aligned
&\iint\frac {\biggl||g(\tau )-g(t)|^2-|w(\tau )-w(t)|^2\biggr|}{|
\tau -t|^2}<\infty\,,\\\text{and}&\iint\frac {\biggl|E|g(\tau )-g
(t)|^2-E|w(\tau )-w(t)|^2\biggr|}{|\tau -t|^2}<\infty\,.\endaligned
\tag 3.3.26$$

Hence together with Lemma (3.3.12) of step (1), this implies

\proclaim{(3.3.27) Corollary}  For
a.e. $g\quad [\nu_{\beta}^{1,*}],$ the limit
$$\lim_{\delta\to 0}\underset{|\tau -t|>\delta}\to\iint\frac {|g(
\tau )-g(t)|^2-E|g(\tau )-g(t)|^2}{|\tau -t|^2}$$
exists.
\endproclaim

\flushpar Proof of (3.3.25).

For $t$ fixed, let $h(\Delta )=g(t)^{-1}g(t+\Delta )$ and $\omega
(\Delta )=w(t+\Delta )-w(t)\,.$
Then
$$h(\Delta )=1+\int^{\Delta}_0h(\tau )\circ d\omega (\tau )\,,\tag 3.3.26$$
and it is shown in section 4.8 of [McKean] that this can be solved
by a Neumann series:
$$\aligned
h(\Delta )&=\sum^{\infty}_0\eta_n(\omega ,t;\Delta ),\quad\text{where}\,\,
\eta_0=1\,,\quad\text{and}\\\eta_n(\Delta )&\equiv\eta_n(\omega ,
t;\Delta )=\int^{\Delta}_0\eta_{n-1}(\tau )\circ d\omega (\tau )\,
.\endaligned
\tag 3.3.27$$
In particular
$$\aligned
\eta_1(\Delta )&=\omega (\Delta ),\quad\text{and}\\\eta_2(\Delta
)&=\int^{\Delta}_0\omega (\tau )\circ d\omega (\tau )\\&=\text{l.i.p.}
\sum\frac {\omega (t_i)+\omega (t_{i-1})}2\left((\omega (t_i)-\omega
(t_{i-1})\right)\endaligned
\tag 3.3.28$$
where the limit in probability is taken over all partitions
$0=t_0<...<t_{\ell}=\Delta ,$ as the mesh tends to zero.  Hence $
\Bbb Re(\eta_1)=0$,
and
$$\Bbb Re(\,tr\,\eta_2(\Delta ))=\frac 12\Bbb Re(\,tr\,\omega (\Delta
)^2)\,.\tag 3.3.29$$
Now suppose that $\tau >t$ and write $\Delta =\tau -t\,.$  Then
$$\aligned
|g(\tau )-g(t)|^2-|w(\tau )-w(t)|^2&=|h(\Delta )-1|^2-|\omega (\Delta
)|^2\\&=-2\Bbb Re\,tr\,\left\{h(\Delta )-(1+\omega (\Delta )+\frac
12\omega (\Delta )^2)\right\}\\&=-2\Bbb Re\,tr\,\sum^{\infty}_3\eta_
n(\Delta )\,.\endaligned
\tag 3.3.30$$

In section 4.8 of [McKean] it is shown that
$$E[tr\,\eta_n(\Delta )^{*}\eta_n(\Delta )]\le c\,\int^t_0E[tr\,\eta_{
n-1}(\tau )^{*}\eta_{n-1}(\tau )]d\tau\,,\tag 3.3.31$$
where $c$ depends only upon $N$ (where $\pi :K\to U(N)).$
In particular
$$E\,tr\,\eta_n(\Delta )^{*}\eta_n(\Delta )\le (c\Delta )^n/_{n!}\,
.\tag 3.3.32$$
It follows that
$$\aligned
E\,tr|\eta_n(\Delta )|&\le NE\{(tr|\eta_n(\Delta )|^2)^{\frac 12}
\}\le N(E\,tr|\eta_n(\Delta )|^2)^{\frac 12}\\&\le N\left(\frac {
(c\Delta )^n}{n!}\right)^{\frac 12}\,.\endaligned
\tag 3.3.33$$

This together with (3.3.32) implies that
$$E\iint\frac {\bigl||g(\tau )-g(t)|^2-|w(\tau )-w(t)|^2\bigr|}{|
\tau -t|^2}$$
$$\le\text{const}\sum^{\infty}_{n=3}(n!)^{-1/2}\iint\frac {|\tau
-t|^{n/2}}{|\tau -t|^2}<\infty\,.\tag 3.3.34$$
This proves (3.3.25) $//$

\bigskip

\flushpar Step 3.

{}From (3.3.25) we can infer that
$$\lim_{\delta\to 0}\underset|\tau -t|>\delta\to\iint\frac {|g(\tau
)-g(t)|^2-E|g(\tau )-g(t)|^2}{|\tau -t|^2}\tag 3.3.35$$
exists for a.e.  $g\in Path^{1,*}K\,,$ relative to almost every Brownian
Bridge.  Given $h\in Path_{C^{\infty}}K$, the convolution of measures
$\delta_h*\nu^{1,1}_{\beta}$ is equivalent to the measure $\nu_{\beta}^{
1,h(2\pi )}$ on $Path^{*,*}K$.
Since
$$\iint\frac {\bigl||h(\tau )g(\tau )-h(t)g(t)|^2-|g(\tau )-g(t)|^
2\bigr|}{|\tau -t|^2}\le\iint\frac {2|h(\tau )-h(t)|^2}{|\tau -t|^
2}<\infty\tag 3.3.36$$
(by the obvious pointwise estimate of the numerator, using the
unitarity of $h$ and $g$), it follows that (3.3.35) exists relative to
all Brownian Bridges, in particular for based loops.

Arguing as in the proof of (3.3.12), using Levy's modulus of
continuity, we can replace $\vert\tau -t\vert$ by $\vert e^{i\tau}
-e^{it}\vert$ in (3.3.35).  This
completes the proof of the existence of the limit of (3.3.6) as
$\delta\to 0$.

We now know that
$$\lim_{\delta\to 0}\sum_{n>0}I_n(\delta )(|\hat {g}(n)|^2-E|\hat {
g}(n)|^2)$$
exists a.e.  $g\in L_CG$.  The fact that $\sum n(\vert\hat {g}(n)
\vert^2-E\vert\hat {g}(n)\vert^2)$
converges conditionally for a.e. $g$ now follows from (ii) of
Lemma (3.3.8).  This completes the proof of (3.3.7)//

\bigskip

\flushpar\S 3.4. Proof of Invariance.

Consider the measure
$$d\mu (g)=det_2\vert A(g)\vert^{2s}exp\left(-s\sum n(|\hat g(n)|^
2-E\vert\hat g(n)\vert^2)\right)d\nu^{*=*}_{\beta}(g)\tag 3.4.1$$
where $M_g=\pmatrix A&B\\C&D\endpmatrix $.  This is well-defined by Proposition
(3.3.7).

If $g_L\in L_{W^1}G,$ then $C_L$ is in $\Cal L_1$, where $M_{g_L}
=\pmatrix A_L&B_L\\C_L&D_L\endpmatrix \,.$ Hence
$$|C|^2-|C(g_Lg)|^2=C^{*}(1-D_L^{*}D_L)C-A^{*}C_L^{*}C_LA-A^{*}C_
L^{*}D_LC+C^{*}D_L^{*}C_LA$$
is also trace class for a.e. $g$, and
$$\sum n(|\hat {g}(n)|^2-E\vert\hat {g}(n)\vert^2)-\sum n(|(g_Lg\hat {
)}(n)|^2-E\vert (g_Lg\hat {)}(n)\vert^2)=tr\{|C|^2-|C(g_Lg)|^2\}\,
.$$
{}From this it is easy to see that $\mu$ has the desired
Radon-Nikodym derivative relative to left translations; a
similar calculation applies for right translations. This completes
the proof of part (a) of the main result, (3.2.1).//

\bigskip

\flushpar\S 3.5. Proof of Finiteness.

Suppose that $g\in LG.$  Because $|A|^2+|C|^2=1\,,$ we have $0\le
|C|^2\le 1\,.$
It follows from this that
$$0\le\det\left((1-|C|^2)e^{|C|^2}\right)\le 1\,,\tag 3.5.1$$
for the function $\lambda\to (1-\lambda )\exp(\lambda )$ is bounded by $
1$ for $\lambda\ge 0$.

It follows that to show the measure $\mu$ of (3.4.1) is finite, it
suffices to prove the following

\proclaim{(3.5.2)Proposition} We have
$$\exp\left(-\sum n(|\hat g(n)|^2-E|\hat g(n)|^2)\right)\in L^p(d
\nu^{*=*}_{\beta}),$$
for all $1\le p<\infty$.\endproclaim

\flushpar Proof.  To prove this we retrace the proof of (3.3.7).
We first observe that (3.5.2) is certainly true in the abelian
case; this is just the assertion that the number $\Cal E$ in (3.1.5) is
finite for any finite $k$.  It then follows as in step 1 of the
proof of (3.3.7) that
$$exp(-\lim_{\delta\to 0}\iint_{\vert t-\tau\vert >\delta}\frac {
\vert w(t)-w(\tau )\vert^2-E\vert w(t)-w(\tau )\vert^2}{\vert t-\tau
\vert^2}dtd\tau )\in L^p(d\nu_{\beta}^{0,*}),$$
for all finite $p$.

The main point is to check that this same statement is true
with $g(t)$ in place of $w(t)$, where $g=I(w)$.   To prove this it
suffices to show that
$$exp(-\iint I(w;t,\tau )dtd\tau )\in L^q(d\nu_{\beta}^{0,*}),\tag 3.5.3$$
for all finite $q$, where
$$I(w;t,\tau )=\frac {\vert g(t)-g(\tau )\vert^2-\vert w(t)-w(\tau
)\vert^2}{\vert t-\tau\vert^2}\tag 3.5.4$$
But
$$exp(-\iint I(w;t,\tau )dtd\tau )\le exp(\iint\vert I(w;t,\tau )
\vert dtd\tau )\tag 3.5.5$$
and $\iint\vert I(w;t,\tau )dtd\tau\in L^1(d\nu^{0,*}_{\beta})$ by (3.3.25).
Thus (3.5.3) follows
from Jensen's inequality.

Finally we localize this result to loops as in step 3. //

\bigskip

\centerline{\S 4. Existence of Invariant Measures.}

\bigskip

In this section we will regard the measures $\nu_{\beta ,k}$ as probability
measures on the formal loop space. We will use the
abbreviated notation
$$L=L_{formal}G=\Cal G/C_0.\tag 4.0.1$$
Our goal is to show that the measures $\nu_{\beta ,k}$ have weak limits
as $\beta\to 0$.  Our strategy will be the same as in Part II.  We will
first prove that the $\Cal N^{\pm}$ distributions have weak limits (\S 4.1 and
\S 4.2).  This will establish the existence of invariant measures
on flag spaces.  We will then consider the diagonal
distributions (\S 4.3).  The combination of these results will
imply weak convergence on the top stratum, and the existence
of invariant measures.

\bigskip

\flushpar\S 4.1. On approximate invariance of Wiener measures.

In chapter 1 of [MM2], the Malliavins state the following

\proclaim{(4.1.1)Proposition}Suppose that $g_L\in L_{W^1}K$.  Then
$$\int\vert 1-\frac {d\nu^{*=*}_{\beta}(g_Lg)}{d\nu^{*=*}_{\beta}
(g)}\vert^pd\nu^{*=*}_{\beta}(g)\to 0\quad as\quad\beta\to 0,$$
for every $1\le p<\infty$. \endproclaim

We will prove the following more precise result.

\proclaim{(4.1.2)Proposition}Suppose that $g_L\in Path^{1,1}_{W^1}
K$.  For all
$1\le p<\infty$ and $k\in K$,
$$\int\vert 1-\frac {d\nu_{\beta}^{1,k}(g_Lg)}{d\nu_{\beta}^{1,k}
(g)}\vert^pd\nu_{\beta}^{1,k}(g)\le 2\Gamma (\frac {p+1}2)\frac {
p_{T/2}(k)p_{T/2}(1)}{p_T(k)}(2\beta \Cal E(g_L))^{p/2},$$
where $p_T(k)$ denotes the heat kernel for $K$, and $T=1/\beta$.
\endproclaim

The basic idea, as usual, is to use Ito's map to recast the
analogous problem for paths in a linear setting, then
localize the result to conditioned paths.   We first consider the
linear problem, in the abstract framework advocated by Irving
Segal; see [S].

Suppose that $H$ is a separable Hilbert space, and let $\nu_{Gaus
s,\beta}$
denote the Gaussian measure associated to $H$ with inverse
temperature $\beta$.  It is well-known that $\nu_{Gauss,\beta}$ is
quasi-invariant with respect to the group of automorphisms
$O(H)\propto H$ of the affine space $H$.

\proclaim{(4.1.3)Proposition}For each $1\le p<\infty$, constant $
c_p$ such
that
$$\int\vert 1-\frac {d\nu_{Gauss,\beta}(g\cdot b)}{d\nu_{Gauss,\beta}
(b)}\vert^pd\nu_{Gauss,\beta}(b)\le 2\Gamma (\frac {p+1}2)(\beta\vert
h_0\vert^2)^{p/2}$$
for all $g\in O(H)\propto H$, where $h_0$ is the translational part of $
g$.

\endproclaim

\flushpar Proof.
$$\frac {d\nu_{Gauss,\beta}(g\cdot b)}{d\nu_{Gauss,\beta}(b)}=exp
(-\frac {\beta}2\vert h_0\vert^2-\beta (h_0,b)),$$
where $h_0$ is the translational part of $g$. This is a function of
one variable. Thus the integral in (4.1.3) is equal to
$$\int_{-\infty}^{+\infty}\vert 1-exp(-\frac {\beta}2\vert h_0\vert^
2+\beta\vert h_o\vert t)\vert^p(\frac {\beta}{2\pi})^{1/2}exp(-\frac {
\beta}2t^2)dt.\tag 4.1.4$$
Writing $s=\sqrt {\beta}\vert h_0\vert$ and changing $\sqrt {\beta}
t$ to $t$, we see that (4.1.4)
equals
$$\int_{-\infty}^{+\infty}\vert 1-exp(-\frac {s^2}2+st)\vert^pexp
(-\frac {t^2}2)dt\le c_ps^p,\tag 4.1.5$$
where
$$c_p=\sup_{0<s<\infty}\int\vert\frac {1-exp(-\frac {s^2}2+st)}s\vert^
pexp(-\frac {t^2}2)dt\tag 4.1.6$$
$$=\int\vert t\vert^pexp(-t^2/2)dt=2\Gamma (\frac {p+1}2)$$
This completes the proof.//

\smallskip

To apply this result to Brownian motion, let $H$ denote the
Hilbert space
$$H=\{x\in W^1([0,1],\frak k):x(0)=0\},\quad\vert x\vert_H^2=\int
\vert\dot {x}\vert^2dt.\tag 4.1.7$$
There is a natural map
$$Path^{1,*}_{W^1}K\to O(H)\propto H,\quad where\quad g\cdot dx=g
(dx)g^{-1}+gd(g^{-1}).\tag 4.1.8$$
If we concretely realize the Gaussian $\nu_{Gauss,\beta}$ as the Wiener
measure $\nu^{0,*}_{\beta}$, then the map
$$I:(Path^{0,*}\frak k,\nu_{\beta}^{0,*})\to (Path^{1,*}K,\nu_{\beta}^{
1,*}):x\to g,\tag 4.1.9$$
where $g$ is the solution to the Fisk-Stratonovich stochastic
differential equation
$$dg\circ g^{-1}=dx,\quad g(0)=1,\tag 4.1.10$$
induces a $Path^{1,*}_{W^1}K$-equivariant isomorphism of measure spaces.
We therefore have the following

\proclaim{(4.1.11)Corollary}For each $1\le p<\infty$,
$$\int\vert 1-\frac {d\nu_{\beta}^{1,*}(g_Lg)}{d\nu_{\beta}^{1,*}
(g)}\vert^pd\nu_{\beta}^{1,*}(g)\le 2\Gamma (\frac {p+1}2)(\beta
\Cal E(g_L))^{p/2}$$
for all $g_L\in Path_{W^1}^{1,*}K$.

\endproclaim

There is an analogous result for right translation, because
the Wiener measure $\nu_{\beta}^{1,*}$ is invariant with respect to
inversion, $g\to g^{-1}$.

The problem now is to localize this result to conditioned paths.
Before addressing this, we will first digress to show how we
can deduce that $\nu^{1,k}_{\beta}$ is quasi-invariant, for every $
k\in K$, from
the quasi-invariance of $\nu^{1,*}_{\beta}$, since our proof of (4.1.2) will
use
the same technique.

\smallskip

\flushpar(4.1.12) Proof that $\nu^{1,k}_{\beta}$ is $Path^{1,1}_{
W^1}K$-quasi-invariant.

Suppose that $g_L\in Path^{1,1}_{W^1}K$.  With respect to the fibration
$$Path^{1,*}K\to K:g\to g(1),$$
there is a essentially unique disintegration
$$(g_L)_{*}\nu^{1,*}_{\beta}=\int_K(g_L)_{*}\nu^{1,k}_{\beta}p_{1
/\beta}(k)dm(k),\tag 4.1.13$$
Essential uniqueness implies that $(g_L)_{*}\nu_{\beta}^{1,k}$ is equivalent to
$
\nu_{\beta}^{1,k}$
for almost every $k\in K$, hence that
$$\frac {d\nu^{1,*}_{\beta}(g_Lg)}{d\nu^{1,*}_{\beta}(g)}=\frac {
d\nu^{1,k}_{\beta}(g_Lg)}{d\nu_{\beta}^{1,k}(g)},\tag 4.1.14$$
where $k=g(1)$, for $a.e.$ $g$ $[\nu^{1,*}_{\beta}]$.  Using the ``principle of
localization'' from the Malliavin calculus for Wiener space (see
Theorem 3.1.3 of [MM1]), the Malliavins deduced that for $\underline {
every}$
$k\in K$, $(g_L)_{*}\nu^{1,k}_{\beta}$ is equivalent to $\nu_{\beta}^{
1,k}$, and therefore it is possible
to define the Radon-Nikodym derivatives so that (4.1.14) holds
for every $a.e.$ $g$ $[\nu^{1,k}_{\beta}]$, for every $k\in K$.  A more direct
argument can be given as follows.

It is very easy to directly compute that
$$\nu^{1,*}_{1/s}*\nu^{1,*}_{1/t}=\nu^{1,*}_{1/(s+t)},\tag 4.1.15$$
as a convolution equation in the group $Path^{1,*}K$. Using the
disintegration (4.1.13), with $g_L=1$, it follows that
$$\int_K(\nu_{1/s}^{1,l}*\nu_{1/t}^{1,l^{-1}k})p_s(l)p_t(l^{-1}k)
dm(l)=\nu_{1/(s+t)}^{1,k}p_{s+t}(k),\tag 4.1.16$$
for a.e.  $k\in K$.  This implies that for every $k\in K$, both sides
of (4.1.16) agree on cylinder functions, because when we use
evaluation at a finite number of points to push the measures
in (4.1.16) to a product of a finite number of $K$'s, both sides of
(4.1.16) are clearly continuous functions of $k$ (of course one
could also directly compute that the two sides agree as well).
Thus (4.1.16) is valid for every $k\in K$.  It follows that for
$g_L\in Path^{1,1}_{W^1}K$, as measure classes
$$\aligned
[\delta_{g_L}*\nu^{1,k}_T]&=[\int_K(\delta_{g_L}*\nu^{1,l}_{T/2}*
\nu_{T/2}^{1,l^{-1}k})dm(l)]\\&=[\int (\nu^{1,l}_{T/2}*\nu^{1,l^{
-1}k}_{T/2})dm(l)]=[\nu^{1,k}_T],\endaligned
\tag 4.1.17$$
for every $k\in K$. //

\smallskip

Now by (4.1.11) we know that
$$\int_KI_{\beta}(k)p_T(k)dm(k)\le 2\Gamma (\frac {p+1}2)(\beta \Cal E
(g_L))^{p/2}\tag 4.1.18$$
where
$$I_{\beta}(k)=\int\vert 1-\frac {d\nu_{\beta}^{1,k}(g_Lg)}{d\nu_{
\beta}^{1,k}(g)}\vert^pd\nu_{\beta}^{1,k}(g)\}\tag 4.1.19$$
We will use the formula (4.1.16) to localize (4.1.18), after
recalling some general facts.

Suppose that $\pi :X\to Y$ is a Borel map of standard Borel spaces.
The induced map on finite measures,
$$\pi_{*}:\Cal M(X)\to \Cal M(Y):\nu\to\pi_{*}\nu\tag 4.1.20$$
is a norm decreasing map.  If we fix one finite positive
measure $\nu\in \Cal M(X)$, then it follows that there is a norm
decreasing map
$$\pi^{\nu}_{*}:L^1(\nu )\to L^1(\pi_{*}\nu )\tag 4.1.21$$
where we view $L^1(\nu )$ as isometrically embedded in $\Cal M(X)$. More
generally there are induced norm decreasing maps
$$\pi^{\nu}_{*}:L^p(\nu )\to L^p(\pi_{*}\nu )\tag 4.1.22$$
for each $p$. To see this, suppose that $f\in L^p(\nu )$. By definition
$$\pi_{*}(fd\nu )=(\pi_{*}^{\nu}f)d\pi_{*}\nu\tag 4.1.23$$
Let $q$ denote the exponent conjugate to $p$. If $g\in L^q(\pi_{*}
\nu )$, then
$$\aligned
\vert\int_Yg\pi_{*}^{\nu}fd\pi_{*}\nu\vert&=\vert\int_Ygd\pi_{*}(
fd\nu )\vert =\vert\int_X(\pi^{*}g)fd\nu\vert\\&\le\vert\pi^{*}g\vert_{
L^q(\nu )}\vert f\vert_{L^p(\nu )}=\vert g\vert_{L^q(\pi_{*}\nu )}
\vert f\vert_{L^p(\nu )},\endaligned
\tag 4.1.24$$
implying that $\pi_{*}^{\nu}f\in L^p(\pi_{*}\nu )$ and $\vert\pi_{
*}^{\nu}f\vert_{L^p(\pi_{*}\nu )}\le\vert f\vert_{L^p(\nu )}$.
\smallskip

\flushpar Proof of (4.1.2).  We apply the remarks of the
preceding paragraph to the (pointwise) multiplication map
$$m:Path^{1,*}K\times Path^{1,*}K\to Path^{1,*}K\tag 4.1.25$$
and the measures
$$\nu =\int_K(\nu^{1,l}_{2\beta}\times\nu^{1,l^{-1}k}_{2\beta})\to
m_{*}\nu =\int_K(\nu_{2\beta}^{1,l}*\nu_{2\beta}^{1,l^{-1}k})=\nu_{
\beta}^{1,k},\tag 4.1.26$$
where the integration is with respect to the measure
$$\frac {p_{T/2}(l)p_{T/2}(l^{-1}k)}{p_T(k)}dm(l).\tag 4.1.27$$
As before, the crucial point is that (4.1.16) is true for every
$k\in K$.

Now the map $m$ is equivariant for the obvious left actions of
$Path_{W^1}^{1,*}K$; thus
$$m_{*}:d\nu (g_Lh_1,h_2)\to d\nu^{1,k}_{\beta}(g_Lh).\tag 4.1.28$$
These measures are absolutely continuous with respect to $\nu$
and $\nu^{1,k}_{\beta}$, respectively.  Since the measures $\nu^{
1,l}_{2/T}\times\nu^{1,l^{-1}k}_{2/T}$ have
disjoint support as $l$ varies over $K$, it follows that
$$m_{*}^{\nu}:\int_K\frac {d\nu^{1,l}(g_Lh_1)}{d\nu^{1,l}(h_1)}\to\frac {
d\nu^{1,k}(g_Lh)}{d\nu^{1,k}(h)},\tag 4.1.29$$
where the integral is with respect to the measure (4.1.27).
Since $m_{*}^{\nu}$ is norm decreasing on $L^p$, we conclude that
$$\int\vert 1-\frac {d\nu_{\beta}^{1,k}(g_Lh)}{d\nu_{\beta}^{1,k}
(h)}\vert^pd\nu^{1,k}_{\beta}(h)\le\int_K\{\int\int\vert 1-\frac {
d\nu_{2\beta}^{1,l}(g_Lh_1)}{d\nu_{2\beta}^{1,l}(h_1)}\vert^pd\nu^{
1,l}_{2\beta}(h_1)d\nu^{1,l^{-1}k}_{2\beta}(h_2)\}$$
$$\le\int\int\vert 1-\frac {d\nu_{2\beta}^{1,l}(g_Lh)}{d\nu_{2\beta}^{
1,l}(h)}\vert^pd\nu_{2\beta}^{1,l}(h)\frac {p_{T/2}(l)p_{T/2}(1)}{
p_T(k)}dm(l).$$
This together with (4.1.18) implies (4.1.2).//

\bigskip

\flushpar\S 4.2. Invariant measures on flag spaces.

Let $\mu_{\beta}^{\vert\tilde {L}\vert^{2k}}$ denote the measure on the formal
loop space with
values in the bundle $\vert\tilde {L}\vert^{2k}$ determined by
$$\vert\sigma_0\vert^{2k}d\mu_{\beta}^{\vert\tilde {L}\vert^{2k}}
=d\nu_{\beta ,k},\tag 4.2.1$$
where $\nu_{\beta ,k}$ is the probability measure constructed in the
previous section.  Let $\Lambda_0$ denote the basic dominant integral
functional, and let $\Cal L=\Cal L_{\Lambda_0}$.  We can project the measure
$\mu_{\beta}^{\vert\tilde {L}\vert^{2k}}$to $\Cal F$ to obtain a measure with
values in $
\vert \Cal L\vert^{2k}$, which we
will denote by $\mu_{\beta}^{\vert \Cal L\vert^{2k}}$.

To prove that there exists an $\tilde {L}_{pol}K$-invariant measure on $
\Cal F$
with values in $\vert \Cal L\vert^{2k}$, we need to prove that $\{
\nu_{\beta ,k}\}$ has weak
limit points.  It would suffice to show that the measures
$\mu_{\beta}^{\vert \Cal L\vert^{2k}}$ satisfy the condition (3.3.2) of \S 3.3
of Part I. This is
trivial for $k=0$, and it seems obvious for $k>0$, but we are
lacking a key step in our proof, as we will now explain.

Fix $g\in L_{pol}K$. We must show that
$$\int\Phi dg_{*}\mu_{\beta}^{\vert \Cal L\vert^{2k}}-\int\Phi d\mu_{
\beta}^{\vert \Cal L\vert^{2k}}\tag 4.2.3$$
tends to zero, for each density $\Phi$ of the form $\Phi =\phi$ $
\vert\sigma_0\vert^{2k}$,
where $\phi$ is bounded and, most importantly, depends upon
finitely many variables.  This difference equals
$$\int\phi (h)\frac {d\nu_{\beta ,k}(h)}{d\nu_{\beta}(h)}(\frac {
d\nu_{\beta}(g^{-1}h)}{d\nu_{\beta}(h)}-1)d\nu_{\beta}(h).\tag 4.2.4$$
If $k=0$, then this tends to zero as $\beta\to 0$ by (4.1.1), and this is
the case for any bounded $\phi$.

Now consider the case $k>0$.  Fix $p>2$ and let $q$ denote the
conjugate exponent.  Since
$$\vert\frac {d\nu_{\beta}(g^{-1}h)}{d\nu_{\beta}(h)}-1\vert_{L^p
(d\nu_{\beta}(h))}\le const\cdot\beta^{p/2},\tag 4.2.5$$
we would be done if we could prove that
$$\vert\frac {d\nu_{\beta ,k}}{d\nu_{\beta}}\vert_{L^q(d\nu_{\beta}
)}\le const/\beta^{power},\tag 4.2.6$$
where the power is less than $p/2$ (note that we do know that
the left side is in $L^q$, by (3.5.1) and (3.5.2)).  Whether there
exists such a bound is in doubt; it certainly is not true for the
abelian analogue; for by (3.1.5) the left side of (4.2.6) is
$$\frac {\Gamma (1+qk/\beta )^{1/q}}{\Gamma (1+k/\beta )}\tag 4.2.7$$
and this does not have a bound of the form (4.2.6). The point
apparently is that one must make essential use of the fact that $
\phi$
depends upon finitely many variables (in the abelian case, if
we project the measure (3.1.4) to the $N$ variables $x_1,..,x_N$,
then the left side of (4.2.6) is bounded by $const\cdot\beta^{1-q}$).

Therefore we have only proven the case $k=0$ of the following

\proclaim{(4.2.8)Conjecture}For each $k\ge 0$, there exists an
$\tilde {L}_{pol}K$-invariant measure on $\Cal F$ having values in $
\vert \Cal L_{\Lambda_0}\vert^{2k}$ such
that
$$\int_{\Cal F}\vert\sigma_0\vert^{2k}d\mu^{\vert \Cal L_{\Lambda_
0}\vert^{2k}}=1.$$
\endproclaim

\flushpar(4.2.9)Remark.  There is potentially a simple direct argument
for the existence of limit points for $\{\nu_{\beta ,k}\}$, basically
independent of (3.3.8) of Part I.   In rough outline the argument
goes as follows.  Write a generic $g\in L_{formal}G$ as
$$g=g_{-}\cdot g_0\cdot g_{+}$$
where $g_{\pm}\in G(\Bbb C[[z^{\pm}]])_1$.  Now fix $\beta >0$.  With $
\nu_{\beta}$ probability one,
$g_{\pm}$ is continuous up to the boundary of $D^{\pm}$, respectively.  This
implies that
$$\nu_{\beta ,k}\{\vert\hat {g}_{+}(n)\vert >R\}\to 0\quad as\quad
n\to\infty$$
for each $R>0$.  For sound reasons (related to the fact that
Wiener measure is formally weighted by the energy of a loop),
it is tempting to speculate that $\nu_{\beta ,k}\{\vert\hat {g}_{
+}(n)\vert >R\}$ is a
nonincreasing function of $n$, in particular that
$$\nu_{\beta ,k}\{\vert\hat {g}_{+}(n)\vert >R\}\le\nu_{\beta ,k}
\{\vert\hat {g}_{+}(1)\vert >R\}.$$
This would reduce the problem of proving tightness down to
the single case $n=1$. This case can be handled very simply
using approximate invariance (which we have unfortunately
established only for $k=0$).

We will consider this argument and its consequences in more
detail in  Part IV.

\bigskip

\flushpar\S 4.3. Existence of biinvariant measures on $L_{formal}
G$.

In this subsection and the next, we will need some detailed
structural information about affine algebras (see especially
chapter 7 of [Kac1], but remember that we are considering the
derived algebra, rather than the full Kac-Moody algebra).  We
will use the Kac convention that symbols associated with
$G=\dot {G}$ will be embellished with a dot, and symbols associated
with $\tilde {L}$ will generally have a tilde.  Let $\dot {A}$ denote the
Cartan
matrix for $\dot {\frak g}$, and let $A$ denote the extended Cartan matrix.
Then
$$\frak g=\frak g'(A)=\tilde {L}_{pol}\dot {\frak g}=L_{pol}\dot {
\frak g}\oplus \Bbb Cc.\tag 4.3.1$$
(To properly develop a theory of roots, one cannot avoid the
semidirect product
$$\Bbb Cd\propto \frak g'(A),\tag 4.3.2$$
where $d$ acts as the derivation $\frac 1i\frac d{d\theta}$, but we will leave
this in
the background).

The simple roots for $\frak g$ are $\{\alpha_j:0\le j\le rk\dot {
\frak g}\}$, where
$$\alpha_0=\delta -\theta ,\quad\alpha_j=\dot{\alpha}_j,j>0,\tag 4.3.3$$
$\delta (d)=1$, $\delta (\frak h)=0$, and $\theta$ is the highest root of $
\dot {\frak g}$.  The simple
coroots of $\frak g$ are $\{h_j:0\le j\le rk\dot {\frak g}\}$, where
$$h_0=c-\dot {h}_{\theta},\quad h_j=\dot {h}_j,j>0,\tag 4.3.4$$
and the $\{\dot {h}_j\}$ are the simple coroots of $\dot {\frak g}$.  There are
integers
$\check {a}_j>0$ such that
$$\dot {h}_{\theta}=\sum\check {a}_j\dot {h}_j.\tag 4.3.5$$
Given $\tilde {g}\in \Cal N^{-}\cdot H\cdot \Cal N^{+}\subset\tilde {
L}$, we have
$$\tilde {g}=l\cdot (diag\tilde {)}\cdot u,\quad where\quad (diag
\tilde {)}(\tilde {g})=\prod_0^{rk\dot {\frak g}}\sigma_j(\tilde {
g})^{h_j}.\tag 4.3.6$$
If $\tilde {g}$ projects to $g\in \Cal N^{-}\cdot\dot {H}\cdot \Cal N^{
+}\subset L$, then because $\sigma_0^{h_0}=\sigma_0^{c-\dot {h}_{
\theta}}$
projects to $\sigma_0^{-\dot {h}_{\theta}}$, we have
$$g=l\cdot diag\cdot u,\quad\tag 4.3.7$$
where
$$diag(g)=\sigma_0(\tilde {g})^{-\dot {h}_{\theta}}\prod_1^{rk\dot {
\frak g}}\sigma_j(\tilde {g})^{\dot {h}_j}=\prod_1^{rk\dot {\frak g}}\left
(\frac {\sigma_j(\tilde {g})}{\sigma_0(\tilde {g})^{\check {a}_j}}\right
)^{\dot {h}_j}\tag 4.3.8$$
Note that $diag$ can be viewed as a (infinite dimensional
analogue of a) rational map,
$$L\to \Bbb P^{rk\dot {\frak g}}:g\to \Bbb P(\sigma_j(\tilde {g})
\sigma_0(\tilde {g})^{\prod_{i\ne j}\check {a}_i}).\tag 4.3.8$$

Given $g\in L_{pol}G$, if
$$g=\sum_{-n}^n\hat {g}(j)z^j\tag 4.3.9$$
with respect to our representation $\pi$, then we will say that $
g$
has $degree\le n$.  For the map
$$i_{\alpha_0}=projection\circ\tilde{\imath}_{\alpha_0}:SL(2,\Bbb C
)\to\tilde {L}_{pol}G\to L_{pol}G,\tag 4.3.10$$
$$deg(i_{\alpha_0}\left(\matrix a&b\\c&d\endmatrix \right))\le n_
0\tag 4.3.11$$
where $n_0$ is a constant that depends upon the choice of
representation (we need $(d\pi (e_{\theta}))^{n_0+1}=0$).

In the case of $G=SL(2,\Bbb C)$,
$$diag=\left(\matrix \sigma_1/\sigma_0&0\\0&\sigma_0/\sigma_1\endmatrix \right
),\quad i_{\alpha_0}(\left(\matrix a&b\\c&d\endmatrix \right))=\left
(\matrix d&cz^{-1}\\bz&a\endmatrix \right)\tag 4.3.12$$
For the classical algebras and standard realizations as in \S 4 of
Part II, $n_0=1$, except for type $B$, where $n_0=2$.

The following technical result was cited in the introduction as
a crucial part of our intuition about the existence of measures.

\proclaim{(4.3.13)Lemma}Suppose that $M>0$. There exists $m>0$ such
that
$$\Cal N^{-}\cdot\dot {H}\cdot \Cal N^{+}\subset \Cal N^{-}_M\cdot
\{g\in L_{pol}K:deg(g)\le m\}\cdot \Cal N^{+}$$
\endproclaim\

The important point here is that
$$\{g\in L_{pol}K:deg(g)\le m\}\tag 4.3.14$$
is a compact set. Thus we are asserting that
$$\Cal N^{-}_M\backslash L/\Cal N^{+}\tag 4.3.15$$
is essentially compact.

\smallskip

\flushpar Proof. We first note that the map
$$\{g\in L_{pol}\dot {K}:deg(g)\le n_0,g\in \Cal N^{-}\cdot\dot {
H}\cdot \Cal N^{+}\}\to\dot {H}:g\to diag(g)\tag 4.3.16$$
is surjective. In fact for $k\in\dot {K}$, $\vert a\vert^2+\vert
b\vert^2=1$,
$$(i_{\alpha_0}\left(\matrix a&b\\-\bar {b}&\bar {a}\endmatrix \right
))k=exp(-\bar {b}a^{-1}e_{\theta}z^{-1})a^{-\dot {h}_{\theta}}exp
(a^{-1}be_{-\theta}z)k,\tag 4.3.17$$
provided $a\ne 0$, hence
$$diag(i_{\alpha_0}\left(\matrix a&b\\-\bar {b}&\bar {a}\endmatrix \right
)k)=a^{-\dot {h}_{\theta}}diag(k)=\prod\left(\frac {\sigma_j(k)}{
a^{\check {a}_j}}\right)^{\dot {h}_j}.\tag 4.3.18$$
Since $a$ is only constrained by $\vert a\vert\le 1$, and
$$diag(K)=\{\prod\eta_j^{\dot {h}_j}:\vert\eta_j\vert\le 1\},\tag 4.3.19$$
it follows from (4.3.18) that (4.3.16) is surjective.

The surjectivity of (4.3.16) implies that
$$\Cal N^{-}_M\cdot \{g\in L_{pol}K:deg(g)\le m\}\cdot \Cal B^{+}
\subset \Cal N^{-}_M\cdot \{g\in L_{pol}K:deg(g)\le m+n_0\}\cdot
\Cal N^{+}\tag 4.3.20$$
Now we also have
$$\{g\in L_{pol}K:deg(g)\le m\}\cdot \Cal B^{+}=\{g\in L_{pol}G:d
eg(g)\le m\}\cdot \Cal B^{+}\tag 4.3.21$$
Thus to complete the proof, it suffices to show that
$$\Cal N^{-}\cdot \Cal B^{+}\subset \Cal N^{-}_M\cdot \{g\in L_{p
ol}G:deg(g)\le m\}\cdot \Cal B^{+}\tag 4.3.22$$
for $m$ sufficiently large.  But (for Kac-Moody groups in
general)
$$\Cal N^{-}/\Cal N^{-}_M=N^{-}/N_M^{-}\tag 4.3.23$$
(where $N^{-}$ is the lower unipotent subgroup for $L_{pol}G$).  Hence
for sufficiently large $m$,
$$\{g\in N^{-}:deg(g)\le m\}\to \Cal N^{-}_M\backslash \Cal N^{-}\tag 4.3.24$$
is surjective. This implies (4.3.22) and completes the proof. //

\proclaim{(4.3.25)Lemma}Given an open subset $U$ of $\dot {H}$,
$$\Cal N^{-}\cdot\dot {H}\cdot \Cal N^{+}\subset\bigcup (g\cdot d
iag^{-1}U),$$
where the union is over all $g\in L_{pol}K$ of $degree\le n_0$.
\endproclaim

\flushpar Proof of (4.3.25).  Fix $\bold g\in \Cal N^{-}\cdot\dot {
H}\cdot \Cal N^{+}$.  We must show
that the map
$$\{g\in L_{pol}K:deg(g)\le n_0,g\in \Cal N^{-}\cdot\dot {H}\cdot
\Cal N^{+}\}\to\dot {H}:g\to diag(g\cdot \bold g)\tag 4.3.26$$
has dense image. We will prove more, namely that the
restriction
$$\{g\in i_{\alpha_r}(SU(2,\Bbb C))...i_{\alpha_0}(SU(2,\Bbb C)):
g\in \Cal N^{-}\cdot\dot {H}\cdot \Cal N^{+}\}\to\dot {H}:g\to di
ag(g\cdot \bold g)\tag 4.3.27$$
is surjective, where $r=rk\dot {\frak g}$.  Note that we can assume
$\bold g=exp(x)\in \Cal N^{-}$.  We establish this by showing the image of
(4.3.27) is closed and open (it is nonempty because 1 is in the
image).

The image of (4.3.27) is closed, because if $\{g_j\}$ is a sequence in the
domain of (4.3.27) and
$$diag(g_jexp(x))\to\dot {h}\in\dot {H},\tag 4.3.28$$
then we can assume that $g_j\to g\in\prod i_{\alpha_i}(SU(2,\Bbb C
)$, because this
set is compact, and because of the finiteness of the limit
(4.3.28), it follows that
$$gexp(x)\in \Cal N^{-}\cdot\dot {H}\cdot \Cal N^{+}\quad and\quad
diag(gexp(x))=\dot {h}.\tag 4.3.29$$

We now proceed to show the image is open. Let
$$g_j=\left(\matrix a_j&b_j\\-\bar {b}_j&\bar {a}_j\endmatrix \right
)\in SU(2,\Bbb C),\quad 0\le j\le r,\tag 4.3.30$$
$$\tilde {g}=\tilde{\imath}_{\alpha_r}(g_r)...\tilde{\imath}_{\alpha_
0}(g_0).\tag 4.3.31$$
We would like to calculate
$$(diag\tilde {)}(\tilde {g}exp(x))=\prod\sigma_j^{h_j}.\tag 4.3.32$$
It is not practical to do this in a completely explicit manner,
but the general form of the answer is given inductively by
$$\sigma_0=(a_0+b_0x_0),\sigma_1=(a_1+b_1\frac {x_1(g_0)}{\sigma_
0}),\sigma_2=(a_2+b_2\frac {x_2(g_0,g_1,\sigma_0)}{\sigma_1}),..\tag 4.3.33$$
where $x_0=\langle x,f_0\rangle$, and each $x_j$ is a polynomial.  Assuming
this
for the moment, we explain why this implies that (4.3.27) has
an open image.  If each $b_i\ne 0$, so that each $\vert a_i\vert
<1$, then
$$d(diag\tilde {)})\vert_{\tilde {g}}\tag 4.3.34$$
is surjective.  If $b_i=0$, then $\sigma_i=a_i$ has absolute value one,
and we have a boundary point for $(diag\tilde {)}$.  But $\sigma_
i$ can
nonetheless be varied in directions $\vert\sigma_i\vert\le 1$.  Thus the image
of
(4.3.27) will contain a set of the form
$$S_0\times S_1\times ..S_r\tag 4.3.35$$
where for each $i$, $S_i$ is either open in $\Bbb C^{*}\cong \Bbb C^{
h_i}$ or relatively
open in $\{0<\vert z\vert\le 1\}$.  The image of this in $\dot {H}$ will
contain an
open neighborhood of $diag(gexp(x))$.

To verify the formulas in (4.3.33), we first calculate
$$\tilde{\imath}_{\alpha_0}(g_0)exp(x)=\tilde{\imath}_{\alpha_0}(
g_0\left(\matrix 1&0\\x_0&1\endmatrix \right))exp(x').\tag 4.3.36$$
Now $G_{(0)}=\tilde{\imath}_{\alpha_0}(SL_2)$  normalizes
$$\sum_{i\ne 0}\Bbb Cf_i\oplus\sum_{height(\alpha )<-1}\frak g_{\alpha}
,\tag 4.3.37$$
(this is the radical for the parabolic $\frak g_{(0)}+\frak b^{-}$). Thus
(4.3.36) equals
$$exp(x^{\prime\prime})\tilde{\imath}_{\alpha_0}(\left(\matrix 1&
0\\\frac {-\bar {b}+\bar {a}x_0}{\sigma_0}&1\endmatrix \right)\left
(\matrix \sigma_0&b_0\\0&\sigma_0^{-1}\endmatrix \right))$$
$$=exp(x^{\prime\prime\prime})\sigma_0^{h_0}exp(\frac {b_0}{\sigma_
0}e_0),\tag 4.3.38$$
where $\sigma_0$ is given by (4.3.33), and $x^{\prime\prime\prime}
\in \frak n^{-}$ is a polynomial in $g_0$
and $\sigma_0^{-1}$.  Thus (4.3.32) equals
$$(diag\tilde {)}(\tilde{\imath}_{\alpha_r}(g_r)..\tilde{\imath}_{
\alpha_1}(g_1)exp(x^{\prime\prime\prime}))\sigma_0^{h_0},\tag 4.3.39$$
and we arrive at (4.3.33) by induction. This completes the
proof. //

\proclaim{(4.3.40)Corollary} Given an open subset $U$ of $\dot {H}$, there
exist a finite number of $g_i\in L_{pol}K$  such that
$$\Cal N^{-}\cdot\dot {H}\cdot \Cal N^{+}\subset\bigcup_{i\le n}(
g_i\cdot diag^{-1}(U)).$$
\endproclaim

\flushpar Proof of the Corollary. By (4.3.25) we have
$$L\subset\bigcup (g\cdot diag^{-1}U),\tag 4.3.41$$
where the union is over all $g\in L_{pol}K$ of $degree\le n_0$.

Now observe that for $M$ sufficiently large, if $deg(g)\le n_0$, then
$$\Cal N_M^{-}\cdot (g\cdot diag^{-1}U)\cdot \Cal N^{+}=g\cdot di
ag^{-1}U\tag 4.3.42$$
The invariance with respect to $\Cal N^{+}$ on the right is obvious, and
the invariance with respect to $\Cal N^{-}_M$ on the left follows from
$$g^{-1}(1+O(z^{-n}))g=1+O(z^{-n+2n_0}),\tag 4.3.43$$
and the left invariance of $diag^{-1}U$ with respect to $\Cal N^{
-}$.

By (4.3.13) there is a compact set $C\subset L$ such that
$$\Cal N^{-}\cdot\dot {H}\cdot \Cal N^{+}\subset \Cal N^{-}_M\cdot
C\cdot \Cal N^{+}.\tag 4.3.44$$
Since we can cover $C$ with a finite number of the sets
$g\cdot diag^{-1}U$, we are done.  //

\proclaim{(4.3.45)Theorem}There exists a $L_{pol}K$ biinvariant
probability measure on $L_{formal}G$.  \endproclaim

\flushpar Proof. Suppose that $U\subset\dot {H}$ is open. By way of
contradiction, suppose that there is a sequence $\beta_j\to 0$ such that
$$\nu_{\beta_j}diag^{-1}U\to 0.\tag 4.3.46$$
For $\{g_i:i\le n\}$ as in the Corollary, we would then have
$$\aligned
1=\nu_{\beta}(\bigcup_i(g_i\cdot diag^{-1}U)&\le\sum_i(\nu_{\beta}
diag^{-1}U+\nu_{\beta}(g_i\cdot diag^{-1}U)-\nu_{\beta}diag^{-1}U
)\\&\le n\nu_{\beta}diag^{-1}U-\sum_i4(\Cal E(g_i)\beta )^{1/2}\endaligned
\tag 4.3.47$$
by (4.1.2).  This is not compatible with (4.3.46).  We conclude that
for any open $U\subset\dot {H}$, there is an $\epsilon >0$ such that
$$\nu_{\beta}diag^{-1}U>\epsilon\tag 4.3.48$$
for all small $\beta$.  This means that $\{diag_{*}\nu_{\beta}\}$ has
nontrivial weak
limits with respect to continuous functions on $\dot {H}$ which vanish
at infinity (put another way, $\{diag_{*}\nu_{\beta}\}\subset Pro
b(\dot {H})$ necessarily
has weak limits (represented by probability measures) on the
one point compactification of $\dot {H}$; any such limit must
necessarily have positive measure on $\dot {H}$ by (4.3.48)).  Moreover, by
(4.1.2), if $diag_{*}\nu_{\beta_j}\to\nu$ in this sense, then
$$diag_{*}(L_g)_{*}(R_h)_{*}\nu_{\beta}\to\nu\tag 4.3.49$$
in this sense as well, for all $g,h\in L_{pol}K$, where $L_g$ and $
R_h$
denote left and right translation, respectively.

We know that the $\Cal N^{\pm}$ distributions of $\{\nu_{\beta}\}$ are tight.
It follows
that there exists a positive measure $\nu$ of total mass $\le 1$ on
$\Cal N^{-}\cdot\dot {H}\cdot \Cal N^{+}$, and there is a sequence $
\beta_j\to 0$ such that $g_{*}\nu_{\beta_j}\to\nu$
weakly with respect to all bounded continuous functions which
vanish at infinity in $\dot {H}$ directions, for all $g\in L_{pol}
K$. This means
that
$$\int fd\nu_{\beta_j}\to\int fd\nu\tag 4.3.50$$
for all $f$ in the $L_{pol}K$ invariant space of functions generated
by bounded continuous functions on $\Cal N^{-}\cdot\dot {H}\cdot
\Cal N^{+}$ which vanish at
infinity in $\dot {H}$ directions.  For such an $f$, and for $g\in
L_{pol}K$, we
then have
$$\aligned
\int fd(L_g)_{*}\nu&=\int f(g(\cdot ))d\nu =\lim\int f(g(\cdot ))
d\nu_{\beta_j}\\&=\lim\int fd\nu_{\beta_j}=\int fd\nu\endaligned
\tag 4.3.51$$
where the third equality follows from (4.1.2).  The class of
functions we are considering clearly separate measures on $L$,
hence $\nu$ is invariant with respect to left translation by $L_{
pol}K$.
Similarly $\nu$ is invariant with respect to right translation by
$L_{pol}K.$ This completes the proof.  //

\bigskip
\flushpar\S 4.4. On Harish-Chandra's {\bf c}-function, and diagonal
distributions.

Any linear function $\lambda$ on $\frak h$ can be written uniquely as
$$\lambda =\dot{\lambda }+\lambda (h_0)\Lambda_0,\tag 4.4.1$$
where $\dot{\lambda}$ can be identified with a linear function on $
\dot {\frak h}$.  A naive
analogue of the {\bf c}-function for the Kac-Moody algebra $\tilde {
L}\dot {\frak g}$ is
formally given by
$$\bold c(\rho -i\lambda )=\prod_{\alpha >0}\frac {\langle\rho -i
\lambda ,\alpha\rangle}{\langle\rho ,\alpha\rangle}=$$
$$\prod_{\dot{\alpha }>0}(1-i\frac {\langle\lambda ,\dot{\alpha}\rangle}{
\langle\dot{\rho }+2\dot {g}\Lambda_0,\dot{\alpha}\rangle})\prod_{
n>0}\left((\prod_{\dot{\alpha}}(1-i\frac {\langle\lambda ,\dot{\alpha }
+n\delta\rangle}{\langle\dot{\rho }+2\dot {g}\Lambda_0,\dot{\alpha }
+n\delta\rangle})^{-1})(1-i\frac {\langle\lambda ,n\delta\rangle}{
\langle\dot{\rho }+2\dot {g}\Lambda_0,n\delta\rangle})^{-rk\dot {
\frak g}}\right)$$
$$=\bold c_{\dot {\frak g}}(\dot{\rho }-i\dot{\lambda })\prod_{n>
0}\left((\prod_{\dot{\alpha}}(1-i\frac {\langle\dot{\lambda },\dot{
\alpha}\rangle +\lambda (h_0)n}{\langle\dot{\rho },\dot{\alpha}\rangle
+2\dot {g}n})^{-1})(1-i\frac {\lambda (h_0)n}{2\dot {g}n})^{-rk\dot {
\frak g}}\right)\tag 4.4.2$$
Here $\dot {g}$ denotes the dual Coxeter number for $\dot {\frak g}$, and $
\dot{\rho}$ and $\rho$ are
(formally) the sums of the positive roots for $\dot {\frak g}$ and $
\frak g$,
respectively.   By pairing the $\dot{\alpha}$ and $-\dot{\alpha}$ terms in the
infinite
product, it is easy to see that the product converges, provided
that $\lambda (h_0)=0$.  As it stands, this expression is problematical
when $\lambda (h_0)\ne 0$.

Now suppose that we had actually succeeded in establishing
the existence of a measure $d\mu^{\vert\tilde {L}\vert^{2k}}$, so that
$$d\mu_k=\vert\sigma_0\vert^{2k}d\mu^{\vert\tilde {L}\vert^{2k}}=\lim_{
\beta\downarrow 0}d\nu_{\beta ,k}.\tag 4.4.3$$
To compute the diagonal distribution, we proceed heuristically,
meaning that we pretend these measures are supported on $LK$,
and there is no problem in lifting them to the $\Bbb T$ extension $
\tilde {L}K$.
This might lead us to expect that
$$\int\dot {a}(g)^{-i\dot{\lambda}}d\mu_k=\int_{\tilde {L}K}\prod
\vert\sigma_j(\tilde {g})\vert^{-i\dot{\lambda }(h_j)}\vert\sigma_
0(\tilde {g})\vert^{(2k+i\dot{\lambda }(\dot {h}_{\theta}))\Lambda_
0(h_0)}/\int_{\tilde {L}K}\vert\sigma_0(\tilde {g})\vert^{2k}$$
$$=\bold c(\rho -i\{\dot{\lambda }+(-\dot{\lambda }(\dot {h}_{\theta}
)+i2k)\Lambda_0\})/\bold c(\rho +2k\Lambda_0)\tag 4.4.4$$
By (4.4.2) the numerator in (4.4.4) equals
$$\bold c_{\dot {\frak g}}(\dot{\rho }-i\dot{\lambda })\prod_{n>0}\left
((\prod_{\dot{\alpha}}(1-i\frac {\langle\dot{\lambda },\dot{\alpha}
\rangle -n\dot{\lambda }(\dot {h}_{\theta})+i2kn}{\langle\dot{\rho }
,\dot{\alpha}\rangle +2\dot {g}n})^{-1})(1-i\frac {-\dot{\lambda }
(\dot {h}_{\theta})n+i2kn}{2\dot {g}n})^{-rk\dot {\frak g}}\right
)\tag 4.4.5$$
The denominator in (4.4.4) does not regularize the divergence
caused by the presence of the $\dot{\lambda }(\dot {h}_{\theta})$ term, but it
does eliminate
other obstacles to convergence.  Without offering a rationale,
we will simply delete the troublesome term.

\proclaim{(4.4.6)Conjecture}Suppose that $\Lambda$ is a dominant integral
functional of level $k$.  Then for any $s\ge 0$,
$$\int\dot {a}(g)^{-i\dot{\lambda}}\vert\sigma_{\Lambda}\vert^{2s}
d\mu^{\vert\tilde {L}\vert^{2ks}}=\bold c_{\dot {\frak g}}(\dot{\rho }
+2s\dot{\Lambda }-i\dot{\lambda })\prod_{n>0}\prod_{\dot{\alpha}}\frac {
(1-i\frac {\langle\dot{\lambda }+i2s\dot{\Lambda },\dot{\alpha}\rangle
+i2skn}{\langle\dot{\rho },\dot{\alpha}\rangle +2\dot {g}n})^{-1}}{
(1-i\frac {i2skn}{\langle\dot{\rho },\dot{\alpha}\rangle +2\dot {
g}n})^{-1}}.$$
In particular,
$$\aligned
\lim_{\beta\downarrow 0}\int\dot {a}(g)^{-i\dot{\lambda}}d\nu_{\beta}
(g)&=\prod_{\dot{\alpha }>0}(1-i\frac {\langle\dot{\lambda },\dot{
\alpha}\rangle}{\langle\dot{\rho },\dot{\alpha}\rangle})^{-1}\prod_{
n>0}\prod_{\dot{\alpha}}(1-i\frac {\langle\dot{\lambda },\dot{\alpha}
\rangle}{\langle\dot{\rho },\dot{\alpha}\rangle +2\dot {g}n})^{-1}\endaligned
\tag 4.4.7$$
\endproclaim

\flushpar(4.4.8)Remarks (a) Note that
$$\lim_{\beta\uparrow\infty}\int\dot {a}(g)^{-i\dot{\lambda}}d\nu_{
\beta}(g)=\int_K\dot {a}(g)^{-i\dot{\lambda}}dm(g)=\prod_{\dot{\alpha }
>0}(1-i\frac {\langle\dot{\lambda },\dot{\alpha}\rangle}{\langle\dot{
\rho },\dot{\alpha}\rangle})^{-1}.$$
Could there possibly be an exact formula for the diagonal
distribution of $\nu_{\beta}$, for finite $\beta$?

(b) In the case of $G=Sl(2,\Bbb C)$, the right hand side of (4.4.7) is
equal to
$$\prod_{n>0}(1+(\frac {\lambda}{2n-1})^2)^{-1}=sech(\frac {\pi}2
\lambda ).\tag 4.4.9$$
where we identify $\lambda\in \Bbb R$ with the functional $\lambda
\dot{\alpha}_1$. If we identify
$$\Bbb R^{+}\to A:a\to\left(\matrix a&0\\0&a^{-1}\endmatrix \right
),\tag 4.4.10$$
then the corresponding distribution of $a(g)$ is given by
$$\{\frac 1{2\pi}\int_{-\infty}^{+\infty}sech(\frac {\pi}2\lambda
)a^{2i\lambda}d\lambda \}dm(a)=\frac 2{\pi}\frac 1{a^2+a^{-2}}\frac {
da}a.\tag 4.4.11$$

\smallskip

In the remainder of this subsection we will work exclusively
with our finite dimensional group $G$, so we drop the dots.

Let $P(\lambda )$ denote the function in (4..4.7).  By (4.4.9) of Part II,
we know that there is a unique $K$-biinvariant probability
measure on $G$, say $\phi (g)dm(g)$, with diagonal distribution
$P^{\check{}}(a)dm(a)$; in theory it can be calculated in the following
way.  From (4.4.27) of Part II, we know that
$$\aligned
\Cal H\phi (\lambda )\bold c(-\rho -i\lambda )=&P(\lambda -i\rho
)\endaligned
,\tag 4.4.12$$
where $\Cal H$ is the Harish transform.  We now apply
Harish-Chandra's inversion formula for $\Cal H$ (in the special case of
a complex group):  for $a\in A$,
$$\aligned
\phi (a)=&c\int \Cal H\phi (\lambda )\frac {\pi (\rho )}{\pi (i\lambda
)}\frac {\sum (det(w))a^{iw\cdot\lambda}}{\sum det(w)a^{w\cdot\rho}}
\vert \bold c(\lambda )\vert^{-2}d\lambda\endaligned
\tag 4.4.13$$
where $c$ is a constant (see Theorems 5.7 and 7.5 of [Helg]).

In the case of $SL(2,\Bbb C)$, if we again identify $\lambda\in \Bbb R$ and $
\lambda\alpha_1$, then
$$\Cal H\phi =-i\frac {\lambda -i}{cos(i\frac {\pi}2(\lambda -i))}
,\tag 4.4.14$$
$$\phi (\left(\matrix a&0\\0&a^{-1}\endmatrix \right))=c\int \Cal H
\phi (\lambda )\frac 1{2i}\frac {a^{i2\lambda}-a^{-i2\lambda}}{\lambda
sinh2loga}\lambda^2d\lambda ,\tag 4.4.15$$
where $c$ is a positive constant. Assuming that $a>1$,
$$\aligned
\int -i\frac {(\lambda -i)\lambda}{cos(i\frac {\pi}2(\lambda -i))}
a^{i2\lambda}d\lambda&=2\pi\sum_{n>0}Res(\frac {(\lambda -i)\lambda}{
cos(i\frac {\pi}2(\lambda -i))}a^{i2\lambda};2ni)\endaligned
$$
$$=2\pi\sum_{n>0}2ni(2n-1)ia^{-4n}\frac 1{-i\frac {\pi}2sin(i\frac {
\pi}2(2n-1)i)}$$
$$=-4ia^{-4}\sum 2n(2n-1)(a^{-2})^{2n-2}(-1)^{n-1}$$
$$=-4ia^{-4}(\frac d{da^{-2}})^2(\sum (-1)^{n-1}(a^{-2})^{2n})$$
$$=-4ia^{-4}\frac {2-6a^{-4}}{(1+a^{-4})^3}\tag 4.4.16$$
Also
$$\int -i\frac {(\lambda -i)}{cos(i\frac {\pi}2(\lambda -i))}a^{-
i2\lambda}d\lambda =2\pi\sum_{n>0}Res(\frac {(\lambda -i)\lambda}{
cos(i\frac {\pi}2(\lambda -i))}a^{-i2\lambda};-2ni)$$
$$=2\pi\sum (-2ni)(-2n-1)ia^{-4n}\frac 1{-i\frac {\pi}2sin(i\frac {
\pi}2(-2n-1)i)}$$
$$=-4ia^{-2}\sum 2n(2n+1)(a^{-2})^{2n-1}(-1)^n$$
$$=-4ia^{-2}(\frac d{da^{-2}})^2(\sum (-1)^n(a^{-2})^{2n+1})$$
$$=-4ia^{-2}(\frac d{da^{-2}})^2(a^{-2}\frac 1{1+(a^{-2})^2})$$
$$=-4ia^{-2}2\frac {a^{-6}-3a^{-2}}{(1+a^{-4})^3}.\tag 4.4.17$$
Inserting (4.4.16) and (4.4.17) into (4.4.15), we see that
$$\phi (\left(\matrix a&0\\0&a^{-1}\endmatrix \right))=const\frac
1{(a^2+a^{-2})^3}.\tag 4.4.18$$

Thus for $G=SL(2,\Bbb C)$, if we write a generic $g$ as
$$g=g_{-}\cdot g_0\cdot g_{+}$$
where $g_{\pm}\in G(\Bbb C[[z^{\pm}]])_1$ and $g_0\in G$, then we are
conjecturing that
$$\lim_{\beta\to 0}(g_0)_{*}\nu_{\beta}=\frac 1Etrace(g_0^{*}g_0)^{
-3}dm(g_0)\tag 4.4.19$$

\bigskip

\flushpar Acknowledgements.  I thank Sam Evens, Hermann
Flaschka, Howard Garland, Nasser Towghi and Gregg
Zuckermann for useful conversations.  I especially thank John
Palmer, who helped me in many ways with this work over
several years.

\bigskip

\centerline{References}

\bigskip

\flushpar[A,etal] S. Albeverio et al, Noncommutative
Distributions, Marcel Dekker, (1993).

\flushpar[BL] A. Beauville and Y. Laszlo, Conformal blocks and
generalized theta functions, CMP 164, 385-419 (1994)

\flushpar[Bill1] P.  Billingsley, Convergence of Probability
Measures, Wiley (1968)

\flushpar[Bill2] ------ , Probability and Measure, Wiley (1986)

\flushpar[CG] K. Clancey and I. Gohberg, Factorization of Matrix
Functions and Singular Integral Operators, Birkhauser (1981)

\flushpar[Gr1] L. Gross, Logarithmic Sobolev inequalities on loop
groups, JFA 101, (1991)

\flushpar[Gr2] -----, Uniqueness of ground states for
Schrodinger operators over loop spaces, JFA 112 (1993) 373-441

\flushpar[GV] R. Gangolli and V.S. Varadarajan, Harmonic
Analysis of Spherical Functions on Real Reductive Groups,
Springer-Verlag (1988)

\flushpar[Hart] R.  Hartshorne, Algebraic Geometry,
Springer-Verlag (1977)

\flushpar[Helg] S. Helgason, Groups and Geometric Analysis,
Academic Press (1983)

\flushpar[IW] N. Ikeda and S. Watanabe, Stochastic Differential
Equations and Diffusion Processes, North Holland (1981)

\flushpar[Kac1] V. Kac, Infinite Dimensional Lie Algebras,
Birkhauser (1984)

\flushpar[Kac2] -----, Constructing groups from infinite
dimensional Lie algebras, in Infinite Dimensional Groups with
Applications, edited by V. Kac, MSRI publication, Springer-Verlag,
(1985)

\flushpar[KP1] V. Kac and D. Peterson, Infinite flag varieties
and conjugacy theorems, Proc. Natl. Acad. Sci., USA, Vol. 80,
1778-1782 (1983)

\flushpar[KP2] ------, Unitary structure in representations of
infinite-dimensional groups and a convexity theorem, Invent.
Math.  76, 1-14 (1984)

\flushpar[KOV] S.  Kerov, G.I.  Olshanskii, A.M.  Vershik,
Harmonic analysis on the infinite symmetric group.  A
deformation of the regular representation, C.R.  Acad.  Sci.,  t.
316, Serie I, 773-778 (1993)

\flushpar[MM1] M. Malliavin and P. Malliavin, Integration
on loop groups, I. Quasi-invariant measures, JFA, Vol 93, No. 1,
207-237, (1990)

\flushpar[MM2] -----, Integration on loop groups, II.
Asymptotic Peter-Weyl orthogonality, preprint.

\flushpar[McKean] H. McKean, Stochastic Integrals, Academic
Press, (1969).

\flushpar[Ol] G. Olshanskii, Unitary representations of infinite
dimensional pairs (G,K) and the formalism of R. Howe, in
Representations of Lie Groups and Related Topics, Advanced
Studies in Contemp. Math., Vol 7, edited by A.M. Vershik and
D.P. Zhelobenko, Gordon and Breach, (1990).

\flushpar[OV] G. Olshanskii and A. Vershik, Ergodic unitarily
invariant measures on the space of infinite Hermitian matrices,
University of Tokyo preprint.

\flushpar[Pi1] D. Pickrell, Measures on infinite dimensional
Grassmann manifolds, JFA, 70, No.2, (1987), 323-356.

\flushpar[Pi2] ------, Separable representations of automorphism
groups of infinite dimensional symmetric spaces, JFA, 90, No.1,
(1990), 1-26.

\flushpar[Pi3] ------, Mackey analysis of infinite classical
motion groups, Pac. J. Math., Vol. 150, No.1, (1991), 139-166.

\flushpar[Pi4] ------, Extensions of loop groups, J.  Algebras,
Groups and Geometry, Volume in memory of Hanno Rund, 1-48,
June 1993.

\flushpar[Pi5] ------, Notes on harmonic analysis on infinite
dimensional symmetric spaces

\flushpar[PS] A. Pressley and G. Segal, Loop Groups, Oxford
University Press, (1986).

\flushpar[RC] E.  Rodriguez-Carrington, dissertation, Rutgers
University (1987).

\flushpar[Sar] D. Sarason, Function Theory on the Unit Circle,
Notes, Virginia Polytechnic Institute.

\flushpar[S] I.E. Segal, Algebraic integration theory, Bull. Amer.
Math. Soc. 71, 419-489 (1965).

\flushpar[W] H. Widom, Aymptotic behavior of block Toeplitz
matrices and determinants. II, Adv. Math. 21, 1-29 (1976)

\flushpar[Z] M. Zhou, dissertation, University of Arizona (1993)

\flushpar[Zygmund] A. Zygmund, Trigonometric Series, Cambridge
University Press (1959)

\end